\documentclass[12pt,reqno]{article}
\usepackage[english]{babel}
\usepackage{amsmath,amsthm}
\usepackage{amsfonts}
\usepackage{graphicx}
\usepackage{enumerate}
\usepackage[usenames,dvipsnames]{color}
\usepackage{subfigure}
\usepackage[right,pagewise,displaymath, mathlines]{lineno}
\usepackage{epstopdf}
\usepackage{color}
\usepackage{multirow}

\usepackage{tikz}
\usepackage{schemabloc}

\numberwithin{equation}{section}
\numberwithin{figure}{section}
\numberwithin{table}{section}

\newtheorem{theorem}{Theorem}[section]

\numberwithin{equation}{section}
\begin{document}
%\linenumbers

\begin{center}

{\Large\bf Assessing transfer functions in control systems}

\vspace*{7mm}

{\large Nadezhda Gribkova}

\medskip

\textit{Faculty of Mathematics and Mechanics, St.\,Petersburg State University,\\ St.\,Petersburg 199034, Russia}

\bigskip

{\large Ri\v cardas Zitikis}

\medskip

\textit{School of Mathematical and Statistical Sciences,
Western University, \break London, Ontario N6A 5B7, Canada}

\end{center}

\bigskip

\begin{quote}

{\bf Abstract.} When dealing with control systems, it is useful and even necessary to assess the performance of underlying transfer functions. The functions may or may not be linear, may or may not be even monotonic. In addition, they may have structural breaks and other abberations that require monitoring and quantification to aid decision making. The present paper develops such a methodology, which is based on an index of increase that naturally arises as the solution to an optimization problem. We show theoretically and illustrate numerically that the empirical counterpart of the index needs to be used with great care and in-depth knowledge of the problem at hand in order to achieve desired large-sample properties, such as consistency.

\medskip


\medskip


\end{quote}

\newpage

\section{Introduction}
\label{intro}

\begin{figure}[h!]
\centering
\[
\stackrel{\displaystyle x_{1}, \dots , x_{n}} {\xrightarrow{\hspace*{3cm}}}
\left \{
 \begin{array}{l}
\begin{tikzpicture}
\sbEntree{E1}
\sbBloc[0]{Bloc1}{Transfer function $h_1(x)$}{E1}
\sbSortie[9]{S1}{Bloc1}
\sbRelier[$y_{1}(1), \dots , y_{n}(1)$]{Bloc1}{S1}
\end{tikzpicture}
\\
\hspace*{10mm}\vdots
\\
\begin{tikzpicture}
\sbEntree{E1}
\sbBloc[0]{Bloc1}{Transfer function $h_d(x)$}{E1}
\sbSortie[9]{S1}{Bloc1}
\sbRelier[$y_{1}(d), \dots , y_{n}(d)$]{Bloc1}{S1}
\end{tikzpicture}
\\
\hspace*{10mm}\vdots
\\
\begin{tikzpicture}
\sbEntree{E1}
\sbBloc[0]{Bloc1}{Transfer function $h_D(x)$}{E1}
\sbSortie[9]{S1}{Bloc1}
\sbRelier[$y_{1}(D), \dots , y_{n}(D)$]{Bloc1}{S1}
\end{tikzpicture}
  \end{array}
\right.
\]
\caption{The dynamical system.}
	\label{fig-1}
\end{figure}
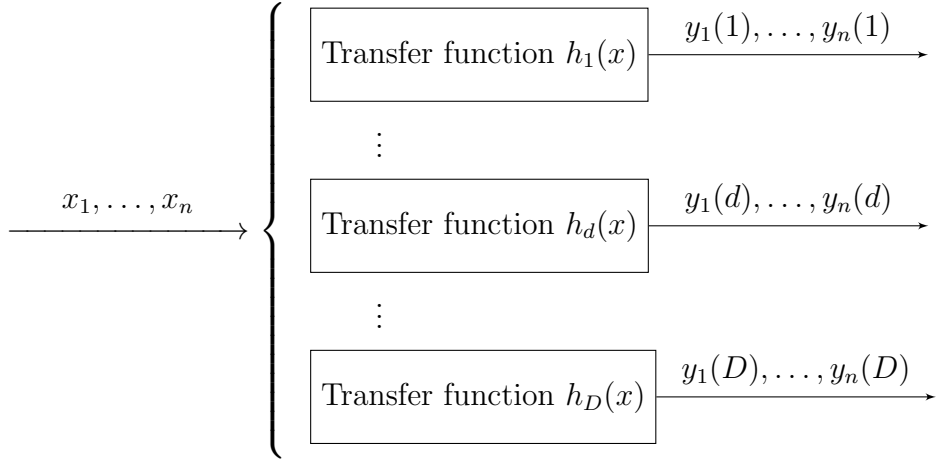

We are interested in assessing the performance of a control system, which could, for example, be one of those arising in time series (e.g., Tong, 1990; Box et al., 2015) or in an architecture such as SISO, SIMO, etc.  of wireless communications (e.g., Tse and Viswanath,~2005; Kshetrimayum,~2017). Let the system be comprised of $D\ge 1$ filters, and let each of the filters $d\in \{1,\dots , D\}$ be represented by a transfer function $h_d(x)$ (Figure~\ref{fig-1}).
The inputs are realized values $x_{1}, \dots , x_{n}$ of random variables $X_{1}, \dots , X_{n}$, and the corresponding outputs are the transformations $y_{1}(d), \dots , y_{n}(d)$ of the $x_i$'s according to the formula $y_{i}(d)=h_d(x_i)$, $1\le i \le n$. That is, the random outputs $Y(d)$ arise via the formula $Y(d)=h_d(X)$ for every $d\in \{1,\dots , D\}$. We denote the cumulative distribution function (cdf) of $X$ by $F(x)$; with $X\sim F$ as the shorthand.

Due to structural breaks or other abberations, the actual transfer functions $h_1(x),\dots , h_D(x)$ may not be known, or just partially known. Naturally, some information about the transfer functions can be gleaned from the scatterplots
\[
\mathcal{S}_d=\big \{(x_i,h_d(x_i))~:~ 1\le i \le n \big \} \quad \textrm{for} \quad d=1,\dots , D,
\]
but as we shall soon see, unanticipated and perhaps even surprising issues arise when estimating, quantifying, and comparing certain features of the transfer functions.

As an example, choose any transfer function $h(x)$ among $h_1(x),\dots , h_D(x)$, and suppose that we wish to asses its monotonicity on a transfer window $(a,b]$. Several methods have been suggested in the literature for the task (e.g., Davydov and Zitikis, 2017;  Chen et al.,  2018; and references therein). Assume that the function $h(x)$ is absolutely continuous (e.g., differentiable). Following Davydov and Zitikis (2017), we can assess and quantify monotonicity of $h(x)$ using the index of increase
\begin{equation}
\mathrm{I}(h)={\int_a^b (h'(x))_{+}\text{d}x \over \int_a^b |h'(x)|\text{d}x },
\label{index-norm}
\end{equation}
where $y_{+}$ denotes the positive part of $y\in \mathbf{R}$, that is, $y_{+}$ is equal to $y$ if $y>0$ and $0$ otherwise. With the notation $g(t)=h(a+t(b-a))$ for $t\in [0,1]$, the index of increase $\mathrm{I}(h)$ reduces to
\begin{equation}
\mathrm{I}(g)={\int_0^1 (g'(t))_{+}\text{d}t \over \int_0^1 |g'(t)|\text{d}t }.
\label{index-norm-g}
\end{equation}
Consider now calculating the index from the practical perspective.

When inputs can freely be chosen by the researcher, they could, for example, be set to the equidistant points $x_{i,n}=a+t_{i,n}(b-a)$ with $t_{i,n}=(i-1)/(n-1)$ for $i=1,\dots , n$. The outputs in this case are $g(t_{i,n})$, $1\le i \le n$, and the index $\mathrm{I}(g)$ can be approximated by
\begin{equation}\label{approx-numeric}
\mathrm{I}_n^0(g)={\sum_{i=2}^n (g(t_{i,n})-g(t_{i-1,n}))_{+} \over \sum_{i=2}^n |g(t_{i,n})-g(t_{i-1,n})| }
\end{equation}
as closely as desired by choosing a sufficiently large $n$. Namely, we have the statement
\begin{equation}\label{conv-I}
\mathrm{I}_n^0(g) \to \mathrm{I}(g)  \quad \textrm{when} \quad n\to \infty ,
\end{equation}
whose numerical explorations can be found in Chen and Zitikis (2017), and Chen et al.~(2018).

Quite often, however, the outputs $g(t_{i,n})$ are contaminated by random errors $\varepsilon_i$, thus making the actual outputs to become
\begin{equation}\label{data-I}
Y_{i,n}=g(t_{i,n})+\varepsilon_i  \quad \textrm{for} \quad  i=1,\dots, n.
\end{equation}
In this case, the index $\mathrm{I}_n^0(g)$ turns into the random index $\mathrm{I}_n^{\varepsilon}(g)$ defined by
\begin{equation}\label{approx-random}
\mathrm{I}_n^{\varepsilon}(g)={\sum_{i=2}^n (Y_{i,n}-Y_{i-1,n})_{+} \over \sum_{i=2}^n |Y_{i,n}-Y_{i-1,n}| } ,
\end{equation}
which, as Chen et al.~(2018) have shown, happens to be asymptotically degenerate, that is,
\begin{equation}\label{conv-0.5}
\mathrm{I}_n^{\varepsilon}(g) \stackrel{\mathbf{P}}{\to}  {1\over 2}  \quad \textrm{when} \quad  n\to \infty ,
\end{equation}
with $\stackrel{\mathbf{P}}{\to}$ denoting convergence in probability. To circumvent the issue, Chen et al.~(2018) have proposed a method based on averaging a certain number of outputs to mitigate the influence of random errors on the non-contaminated outputs $g(t_{i,n})$. In the current paper we further advance the understanding of this research area.

We have organized the rest of the paper as follows. In Section~\ref{main} we place the above observations within a general framework and specify assumptions under which the empirical index of increase converges to a quantity that can be used to quantify and compare the (lack of) increase of various transfer functions; three theorems are devoted to the issue. In Section~\ref{illustration} we report findings of extensive numerical and graphical explorations of the empirical index of increase, thus providing an illustration of how the herein developed theory works in practice. Section~\ref{proofs} contains proofs of the main theorems of Section~\ref{main}. A brief summary and concluding remarks make up Section~\ref{conclude}.

\section{Main results}
\label{main}

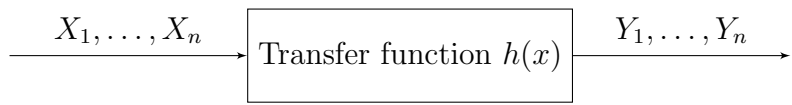
\begin{figure}[h!]\bigskip
\centering
\begin{tikzpicture}
\sbEntree{E1}
\sbBloc[8]{Bloc1}{Transfer function $h(x)$}{E1}
\sbRelier[$X_{1},\dots , X_{n}$]{E1}{Bloc1}
\sbSortie[7]{S1}{Bloc1}
\sbRelier[$Y_{1}, \dots , Y_{n}$]{Bloc1}{S1}
\end{tikzpicture}
\caption{Is the transfer function $h(x)$ still increasing or has sustained a breakdown?}
	\label{fig-2}
\end{figure}

Let $X_{1},\dots , X_{n} \sim F$ be independent and identically distributed (iid) inputs, with the corresponding outputs $Y_{i}=h(X_{i})$ arising via an (unknown) transfer function $h(x)$ (Figure~\ref{fig-2}).
We wish to assess monotonicity of the function either on its entire domain of definition or only on its sub-domain. From the practical point of view, we can do so only via the available for us information,  which contains the realized values of the input variables as well as the corresponding outputs. Under this set-up, the empirical index of increase is
\begin{equation}\label{random-index-10}
\mathrm{I}_n={\sum_{i=2}^n (h(X_{i:n})-h(X_{i-1:n}))_{+} \over \sum_{i=2}^n |h(X_{i:n})-h(X_{i-1:n})|}.
\end{equation}

We have already alluded to the fact that the index may or may not be asymptotically degenerate, and so it is natural to ask whether or not -- and if yes, then when -- the index can be used for the intended purpose, which is to assess and quantify the (lack of) monotonocity of the transfer function $h(x)$. In the following two complimenting each other theorems, we resolve this issue.

We make the following assumptions:
\begin{enumerate}[(C1)]
\item\label{cond-1}
The cdf $F(x)$ of $X$ is absolutely continuous, and we denote its probability density function (pdf) by $f(x)$.
\item\label{cond-2}
There are points $\tau_0:=0<\tau_1<\dots<\tau_m<\tau_{m+1}:=1$, for some $m\geq 0$, such that the function
\begin{equation}\label{funct-g}
H(u)={h'\circ F^{-1}(u) \over f\circ F^{-1}(u)}
\end{equation}
is well-defined and continuous on every interval $(\tau_{k},\tau_{k+1})$, $k=0,\dots , m$, and has finite left- and right-hand limits at every point $\tau_1,\dots,\tau_m \in (0,1)$.
\item\label{cond-3}
The function $h\circ F^{-1}(u)$ is left-continuous at the points $\tau_1,\dots,\tau_m \in (0,1)$, which are defined in condition (C\ref{cond-2}), and has finite right-hand limits at these points. Furthermore, either $h\circ F^{-1}(u)$ has at least one non-zero jump among the points $\tau_1,\dots,\tau_m \in (0,1)$, or the function $h'\circ F^{-1}(u)$ is not identically equal to $0$ on the interval $(0,1)$.
\end{enumerate}

The following theorem, formulated under the additional (to be subsequently relaxed) assumption that the function $H(u)$ has finite right- and left-hand limits at the endpoints of its domain of definition $(0,1)$, introduces what we call the index of increase, denoted by $\mathrm{I}_F(h)$, and whose several counterparts have already appeared above.

\begin{theorem}
\label{Thm-2}
Let conditions (C\ref{cond-1})--(C\ref{cond-3}) be satisfied, and let the function $H(u)$ have finite right- and left-hand limits at $\tau_0=0$ and $\tau_{m+1}=1$, respectively. Then
\begin{equation}\label{ii-main-2}
\mathrm{I}_n \stackrel{\mathbf{P}}{\longrightarrow}
\mathrm{I}_F(h):=\frac{\sum_{k=1}^m (\Delta h_k)_{+}+\int_{0}^{1} H_{+}(u)\mathrm{d}u}
{\sum_{k=1}^m|\Delta h_k|+\int_{0}^{1}|H|(u)\mathrm{d}u}
\end{equation}
where
\begin{itemize}
  \item $\Delta h_k:=h\circ F^{-1}(\tau_k+0)-h\circ F^{-1}(\tau_k)$, $1\le k \le n$, are the jump sizes of the function $h\circ F^{-1}(u)$ at the points $\tau_k$ defined in condition (C\ref{cond-2}), with the two sums $\sum_{k=1}^m$ vanishing when $m=0$;
  \item $H_{+}(u)=\max \{H(u), 0\} $, $H_{-}(u)=\max \{-H(u), 0\} $, and $|H|(u)=H_{+}(t)+H_{-}(u)$.
\end{itemize}
\end{theorem}

The next theorem relaxes the aforementioned additional assumption of Theorem~\ref{Thm-2} by allowing $H(u)$ to grow indefinitely -- though not too fast -- when it approaches the endpoints of its domain of definition $(0,1)$. This we achieve at the expense of assuming differentiability of $H(u)$ and adding certain assumptions on the growth of $H(u)$ and its derivative $H'(u)$ when they are approaching the endpoints of $(0,1)$.

\begin{theorem}
\label{Thm-3}
Let conditions (C\ref{cond-1})--(C\ref{cond-3}) be satisfied, and let there be constants $c>0$ and $b>0$ such that
\begin{equation}
\label{iib}
|H(u)|\leq c  \big ( u(1-u) \big )^{-1+b} \quad \textrm{for} \quad u\in (0,1)\setminus \{\tau_1,\dots, \tau_k\}.
\end{equation}
Furthermore, let the function $H(u)$ be differentiable on the set
$\Theta:=(0,\theta)\cup(1-\theta,1)$ for some (small) $\theta >0$, and let the bound
\begin{equation}
\label{iibb}
|H'(u)|\leq c \big (u(1-u) \big )^{-2+b} \quad \textrm{for} \quad u\in \Theta
\end{equation}
be satisfied with the same constant $b>0$ as in condition \eqref{iib}, but the constant $c<\infty $ can be different. Then statement \eqref{ii-main-2} holds.
\end{theorem}

In the next section, we shall numerically illustrate the above two theorems under various scenarios. The rest of the current section is devoted to sheading more light on why and when the asymptotic degeneracy of the empirical index $\mathrm{I}_n$ occurs. This can be done in a very general setup. Hence, unless explicitly stated otherwise, for the rest of this section we deal with generic random pairs $(X_1,Y_1), \dots, (X_n,Y_n)$, which, for the sake of transparency, we assume to be independent copies of $(X,Y)$. We assume that the marginal cdf $F(x)$ of $X$ is continuous, in which case, almost surely, we can unambiguously order the first coordinates $X_1,\dots , X_n$ in the strictly ascending fashion. This gives rise to the  order statistics $X_{1:n}< \cdots < X_{n:n}$ and their corresponding concomitants $Y_{1,n}, \dots, Y_{n,n}$ (e.g., David and Nagaraja,~2003). Based on the latter ones, we define the index of increase
\begin{equation}\label{random-index}
\mathrm{I}_n={\sum_{i=2}^n (Y_{i,n}-Y_{i-1,n})_{+} \over \sum_{i=2}^n |Y_{i,n}-Y_{i-1,n}|},
\end{equation}
whose definition mimics those introduced earlier. The following is a very general result that, as we shall see below, introduces an almost optimal condition which describes the set of those pairs $(X_i,Y_i)$ that result in the asymptotic degeneracy of $\mathrm{I}_n$ when $n\to \infty $.

\begin{theorem}
\label{Thm-1}
Let $Y$ have finite second moment. If the degeneracy condition
\begin{equation}\label{ii-main}
 B_n:={1\over \sqrt{n}} \sum_{i=2}^n \big | Y_{i,n}-Y_{i-1,n} \big |
 \stackrel{\mathbf{P}}{\to} \infty
\end{equation}
is satisfied, then $\mathrm{I}_n\stackrel{\mathbf{P}}{\to} 1/2$ when $n\to \infty $.
\end{theorem}

With the help of Theorem \ref{Thm-1}, the remainder of the current section is devoted to an explanation of why statement (\ref{conv-I}) holds when there are no measurement errors, and statement (\ref{conv-0.5}) holds when there are (non-degenerate) errors. Naturally, our explanation, which consists of two parts, relies on the asymptotic behaviour of the quantity $B_n$.

\paragraph{\it Part 1.}

Consider the case when there are no measurement errors, and let the function $g(t)$ be  $\gamma$-H\"{o}lder continuous for some $\gamma \ge 1/2$. (Note, for example, that $g(t)$ with a uniformly on $[0,1]$ bounded derivative $g'(t)$ would correspond to the case $\gamma=1$.) Hence, we have
\begin{align*}
B_n&={1\over \sqrt{n}} \sum_{i=2}^n \big | g(t_{i,n})-g(t_{i-1,n}) \big |
\\
&\le {c\over \sqrt{n}} \sum_{i=2}^n \big | t_{i,n}-t_{i-1,n} \big |^{\gamma }
\\
&= {c\over n^{\gamma -1/2}} ,
\end{align*}
which implies $B_n=O(1)$ when $n\to \infty $, which is the opposite of degeneracy condition (\ref{ii-main}). In this case, as noted by Davydov and Zitikis (2017) and further explored by Chen and Zitikis (2017), we have statement (\ref{conv-I}).

\paragraph{\it Part 2.}

We now look at model (\ref{data-I}) with non-degenerate measurement errors $\varepsilon_{i}$, which for the sake of argument are assumed to be iid with finite second moments. As before, we let $g(t)$ be $\gamma$-H\"{o}lder continuous for some $\gamma \ge 1/2$. We have
\begin{align*}
B_n&\ge {1\over \sqrt{n}} \sum_{i=2}^n | \varepsilon_{i}-\varepsilon_{i-1} |
-{1\over \sqrt{n}} \sum_{i=2}^n \big | g(t_{i,n})-g(t_{i-1,n}) \big |
\\
&\ge {1\over \sqrt{n}} \sum_{i=2}^n | \varepsilon_{i}-\varepsilon_{i-1} |
-{c\over n^{\gamma -1/2}}
\\
&= {1\over \sqrt{n}} \sum_{i=2}^n \Big (| \varepsilon_{i}-\varepsilon_{i-1} |
- \mathbf{E}\big [| \varepsilon_{i}-\varepsilon_{i-1} | \big ]\Big )
+ {n-1\over \sqrt{n}} \mathbf{E}\big [| \varepsilon_{2}-\varepsilon_{1} | \big ]
-{c\over n^{\gamma -1/2}}
\\
&= \sqrt{n}~ \mathbf{E}\big [| \varepsilon_{2}-\varepsilon_{1} | \big ]
+O_{\mathbf{P}}(1)
\end{align*}
when $n\to \infty $. The right-hand side tends to infinity in probability because $\mathbf{E} [| \varepsilon_{2}-\varepsilon_{1} |  ]$ is positive whenever $\varepsilon_{i}$'s are non-degenerate, which we assume. Hence,  degeneracy condition (\ref{ii-main}) is satisfied. Recall that in this case statement (\ref{conv-0.5}) holds, which is corroborated by Theorem \ref{Thm-1}. We refer to Chen et al.~(2018) for additional details on this topic, as well as to the earlier work of Davydov and Zitikis (2007) on distinguishing deterministic and random noises.

\section{A numerical illustration}
\label{illustration}

We have subdivided this section into two major parts, which are Sections \ref{section-31} and \ref{section-32}. In the first part, we specify an underlying signal distribution, from which we then generate input random variables. In the same section, we also specify the parameter values that we use in our following numerical explorations, and we also specify a transfer function. In the second part, which is Section~\ref{section-32}, we report findings of our numerical and graphical explorations.

\subsection{Model specification}
\label{section-31}

For the sake of illustrative simplicity, we use iid inputs, which are copies of a random variable $T$ to be specified in a moment, to assess monotonicity of the transfer function $h(x)$ on certain subintervals $(a,b]$, called transfer windows, of its domain of definition. These windows $(a,b]$ must also be subsets of the support of the random variable $T$, since only in this case can we guarantee some data to at least hope to assess monotonicity of $h(x)$ on $(a,b]$. In the specification of $T$ that follows, the support of this random variable is the entire real-half line $[0,\infty )$, and thus, in principle, we can choose any window $(a,b]\subset [0,\infty )$ and have some data in it: sometimes plentiful but sometimes not. The crucial relationship between $(a,b]$ and the availability of data will be discussed below.

\subsubsection{The signal distribution}

In numerous problems, quite notably in those related to the Internet traffic, inputs are heavy tailed; yet modelling small inputs is also important. In such cases, the Lomax distribution has turned out to be a good model (e.g., Holland et al., 2006; Weigle, 2006; Mattos et al., 2014; Chen et al., 2015; Wilson et al., 2018). The distribution is Pareto-like, and we refer to Arnold (2015) for an authoritative treatment of such distributions. The Lomax cdf is given by the formula
\begin{equation}
F_T(t)=1-\bigg ( {\beta \over t+\beta} \bigg )^{\alpha }  \quad \textrm{for} \quad  t\ge 0,
\label{cdf-t}
\end{equation}
where $\alpha >0$ and $\beta >0$ are the shape and scale parameters, respectively. The corresponding pdf and quantile functions are
\begin{equation}
f_T(t)={\alpha\beta^{\alpha} \over (t+\beta)^{\alpha+1}}  \quad \textrm{for} \quad   t\ge 0,
\label{pdf-t}
\end{equation}
and
\begin{equation}
F_T^{-1}(u)=\beta \bigg ( {1\over (1-u)^{1/\alpha}} -1\bigg )  \quad \textrm{for} \quad   0<u<1,
\label{qf-t}
\end{equation}
respectively. In what follows, we shall also need the density-quantile function of $T$, which in view of equations (\ref{pdf-t}) and (\ref{qf-t}) has the expression
 \begin{equation}
f_T\circ F_T^{-1}(u)={\alpha \over \beta } (1-u)^{1+1/\alpha}   \quad \textrm{for} \quad   0<u<1.
\label{dqf-t}
\end{equation}

\subsubsection{The input distribution and window choices}

Since we are concerned with the transfer function $h(x)$ on the window $(a,b]$, we input into the filter only those values of $T$ that fall into the window $(a,b]$. In other words, $X$ is the random variable $T$ conditioned on $T\in (a,b]$. Concisely, we write
\[
X=\langle T \mid T\in (a,b] \rangle .
\]
The cdf of $X$ is given by the formula
\begin{equation}
F(x)=
\left\{
  \begin{array}{ll}
    0 & \hbox{ when } x\le a,\\
    \displaystyle {F_T(x)-F_T(a)\over F_T(b)-F_T(a)} & \hbox{ when } x\in (a,b],\\
    1 & \hbox{ when } x>b, \\
  \end{array}
\right.
\label{cdf-x}
\end{equation}
and the quantile function is
\begin{equation}
 F^{-1}(t)= F_T^{-1}\big (F_T(a)+t ( F_T(b)-F_T(a) )\big )  \quad \textrm{for} \quad   0<t<1.
\label{qf-x}
\end{equation}
The latter equation and the closed-form expressions for $F_T(x)$ and $F_T^{-1}(t)$ allow us to conveniently simulate $X_1,\dots , X_n$ by first simulating uniform on $[0,1]$ random variables $U_1,\dots , U_n$ and then setting  $X_i=F^{-1}(U_i)$.

The Lomax distribution is absolutely continuous, and thus the cdf $F(x)$ of $X$ is also absolutely continuous, as required by condition (C\ref{cond-1}). Its pdf is
\begin{equation}
f(x)=
\left\{
  \begin{array}{ll}
    \displaystyle {f_T(x)\over F_T(b)-F_T(a)} & \hbox{ when } x\in (a,b], \\
0 & \hbox{ otherwise, }
  \end{array}
\right.
\label{gen-d}
\end{equation}
and the density-quantile function is
\begin{equation}
f\circ F^{-1}(t)={1\over F_T(b)-F_T(a) }
f_T\circ F_T^{-1}\big (F_T(a)+t ( F_T(b)-F_T(a) )\big ) \quad \textrm{for} \quad  0<t<1.
\label{gen-dq}
\end{equation}
In the following numerical explorations, we work with the three windows
\begin{equation}
(0,2],\quad (8,12],  \quad \text{and} \quad (0,20].
\label{par-4}
\end{equation}

\subsubsection{Lomax parameter choices}

The mean of the Lomax distribution exists when $\alpha >1$ and is equal to $\beta/(\alpha-1)$. Quite frequently, especially in the aforementioned Internet-related applications, the mean exists but the variance does not. Guided by this observation, in our numerical explorations we set
\begin{equation}
\alpha=1.5 \quad \text{and} \quad \beta=1.
\label{par-1}
\end{equation}
Hence, the mean of $T$ is $2$. For the variance to exist, which is equal to $\beta^2\alpha /((\alpha-1)^2(\alpha-2))$, we must have $\alpha>2$. To illustrate this case, we set
\begin{equation}
\alpha=5 \quad \text{and} \quad \beta=1.
\label{par-2}
\end{equation}
The mean of the latter $T$ is equal to $0.25$, and the variance is $5/48\approx 0.1042$. Hence, most of the observations of $T$ congregate near $x=0$, with only few ones beyond, say, $x=15$. This, as we shall see in Section~\ref{section-323}, will be an impediment when estimating the true index of increase in the window $(0,20]$.

\begin{figure}[h!]
\centering
\subfigure[$\alpha=1.5$]{%
\includegraphics[height=0.35\textwidth]{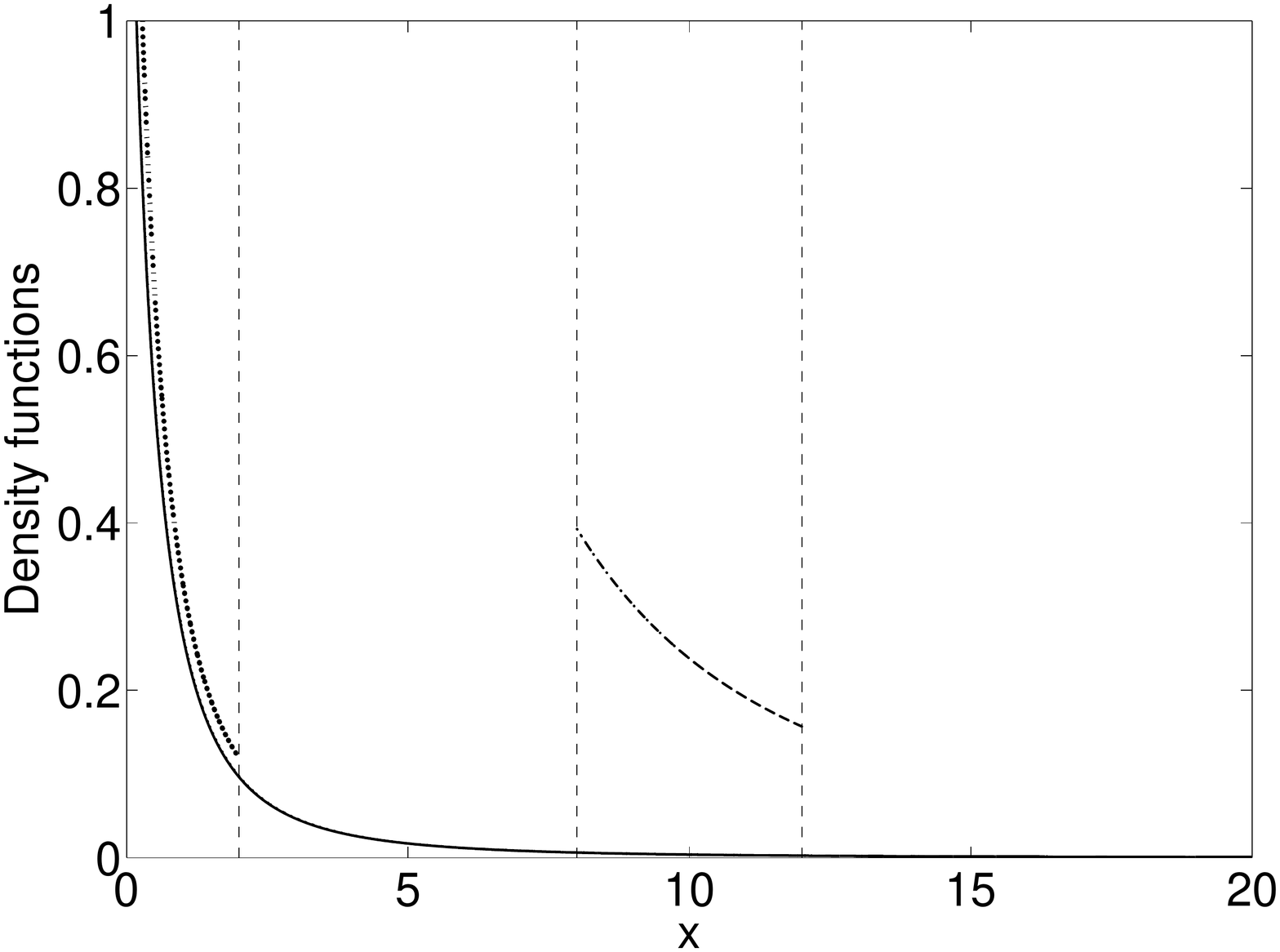}}
\hfill
\subfigure[$\alpha=5$]{%
\includegraphics[height=0.35\textwidth]{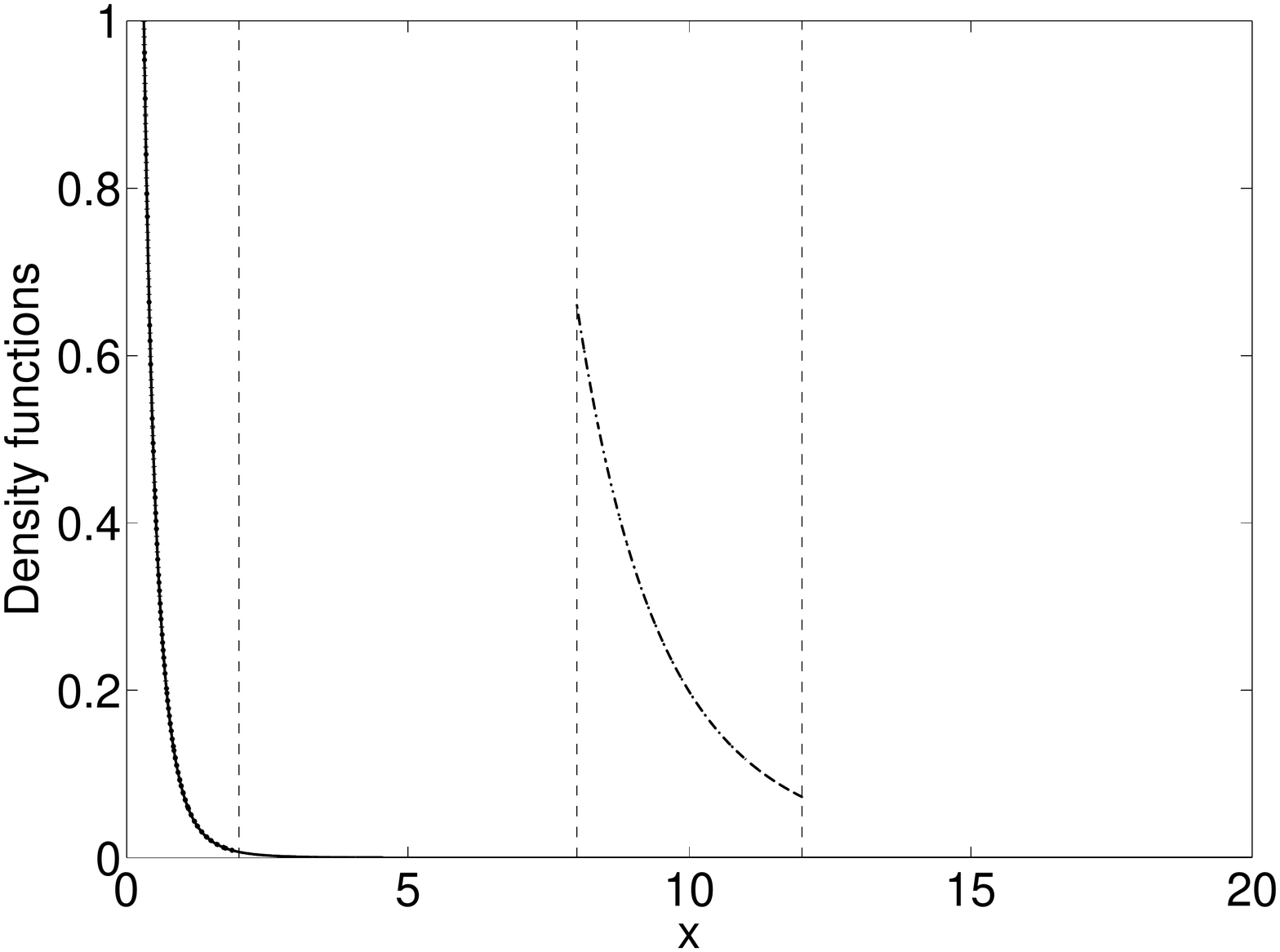}}
\caption{The Lomax pdf's: the unconditional (solid) and conditional ones in the windows $(0,2] $ (dotted) and $(8,12]$ (dash-dotted).}
\label{fig-lomax}
\end{figure}

To better understand our following numerical results,  in Figure~\ref{fig-lomax}
we have depicted the unconditional and conditional pdf's in windows (\ref{par-4}). Namely, we have drawn the Lomax pdf $f_T(t)$ (equation (\ref{pdf-t})) under the above two sets of parameters, where we have also depicted the conditional pdf $f(x)$ (equation (\ref{gen-d})) in windows (\ref{par-4}).  Not all of these conditional pdf's are clearly visible in the figures: the one in the window $(0,20]$ virtually coincides with the unconditional pdf under both sets of parameters (\ref{par-1}) and (\ref{par-2}), and so is the one in the window $(0,2]$ under the set of parameters (\ref{par-2}). The conditional pdf in the window $(8,12]$ is quite distinct from the unconditional pdf under both sets of parameters (\ref{par-1}) and (\ref{par-2}), and thus clearly visible in the two panels of Figure \ref{fig-lomax}. We shall depict more minute details of these conditional pdf's in Figures~\ref{pdf-02} and \ref{pdf-812} below.

\subsubsection{Transfer function}

\begin{figure}[h!]
\centering
\includegraphics[height=0.35\textwidth]{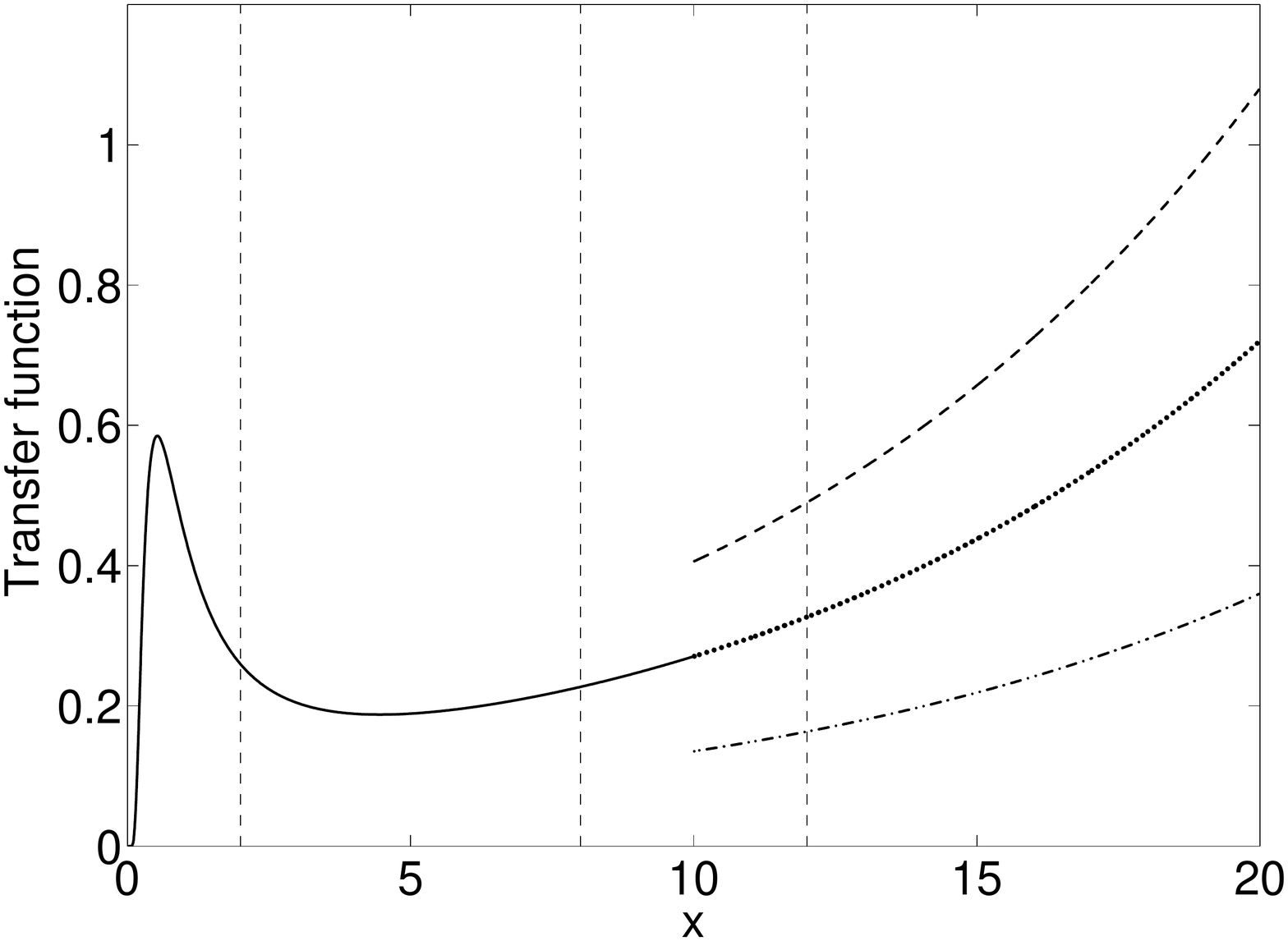}
\caption{The transfer function $h(x)$ below the threshold $x_0=10$ (solid) and above it when $\varrho = -0.5$ (dashed), $\varrho =  0$ (dotted), and $\varrho =  0.5$ (dash-dotted), with the transfer windows delineated by vertical lines at $x=2$, $8$, and $12$.}
\label{fig-h}
\end{figure}

The four-parameter transfer function that we use in our numerical explorations is given by
\begin{equation}
h(x)=\bigg( \gamma + {\delta \over x^2}\bigg) e^{\gamma x -\delta/x} \big ( 1+\varrho \mathbf{1}_{(0,\infty)}(x-x_0) \big )  \quad \textrm{for} \quad   x\ge 0.
\label{tf-1}
\end{equation}
Unless $\varrho =0$, the function has a regime switching at the point $x=x_0$. We choose to work with the following parameter values:
\begin{equation}
\gamma=0.1,~ \delta =1,~ x_0 =10, \quad \text{and} \quad \varrho \in \{0,\pm 0.5\}.
\label{par-3}
\end{equation}
Hence, in the window $(8,12]$, we may have a regime switching, which could be either a drop or a jump depending on the sign of $\varrho $. Specifically, when $\varrho = -0.5$, the transfer function has a sudden drop at the point $x_0=10$, and when $\varrho = 0.5$, it has a sudden jump at the point. When $\varrho=0$, the function is continuous on the entire real half-line; it actually coincides with the hazard rate function of Bebbington et al.~(2007). We have depicted the function $h(x)$ in all these three cases in Figure~\ref{fig-h}.
In the following section we shall explore the performance of the empirical index $\mathrm{I}_n$ when the inputs are fed into the filter through windows (\ref{par-4}).

\subsection{Numerical explorations}
\label{section-32}

\begin{figure}[h!]
\centering
\subfigure[$\alpha =  1.5$, window \textrm{(0,2]}, $\mathrm{I}_F(h)\approx 0.6424$.]{%
\includegraphics[height=0.35\textwidth]{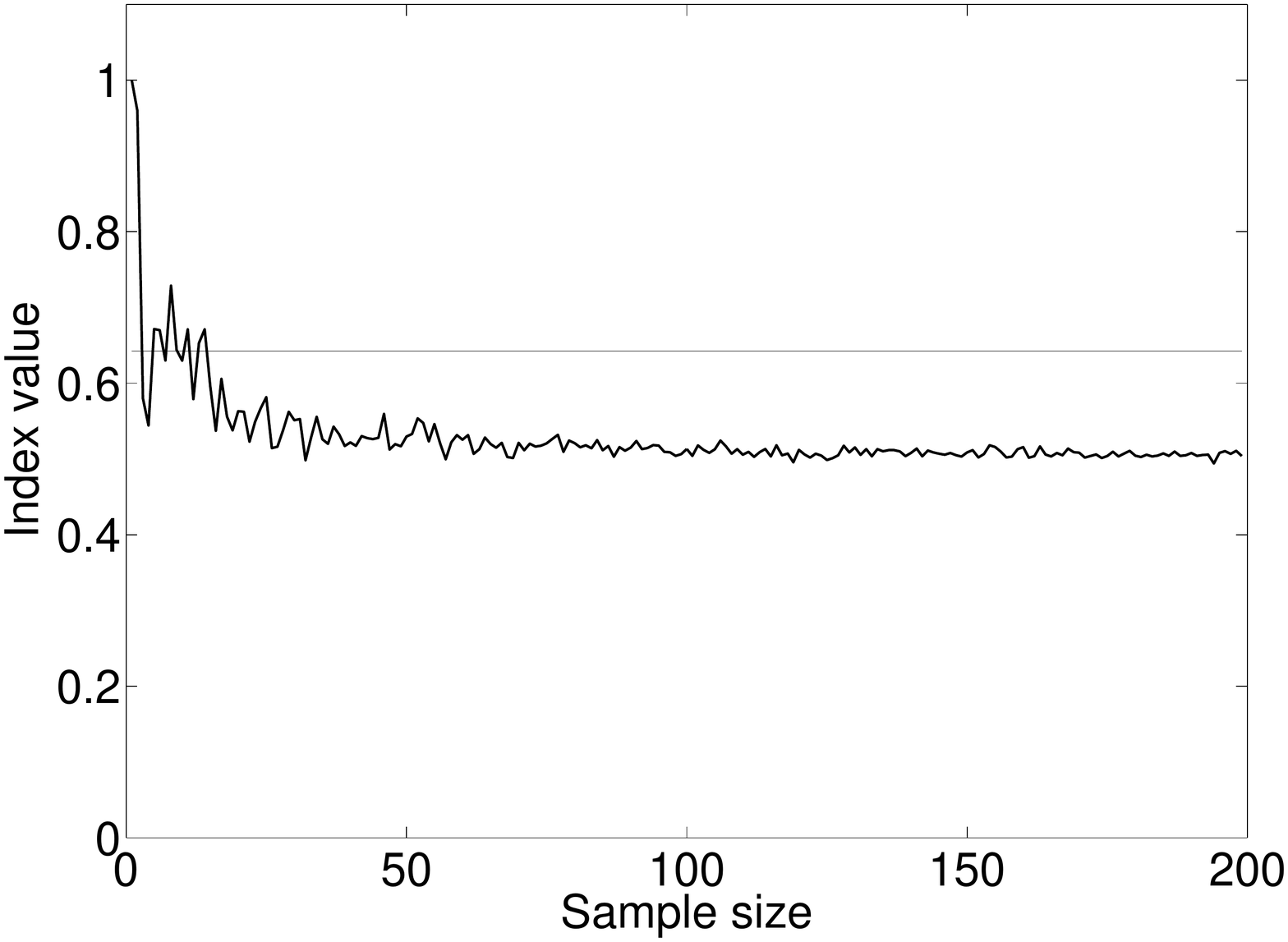}}
\hfill
\subfigure[$\alpha =  5$, window \textrm{(0,2]}, $\mathrm{I}_F(h)\approx 0.6424$.]{%
\includegraphics[height=0.35\textwidth]{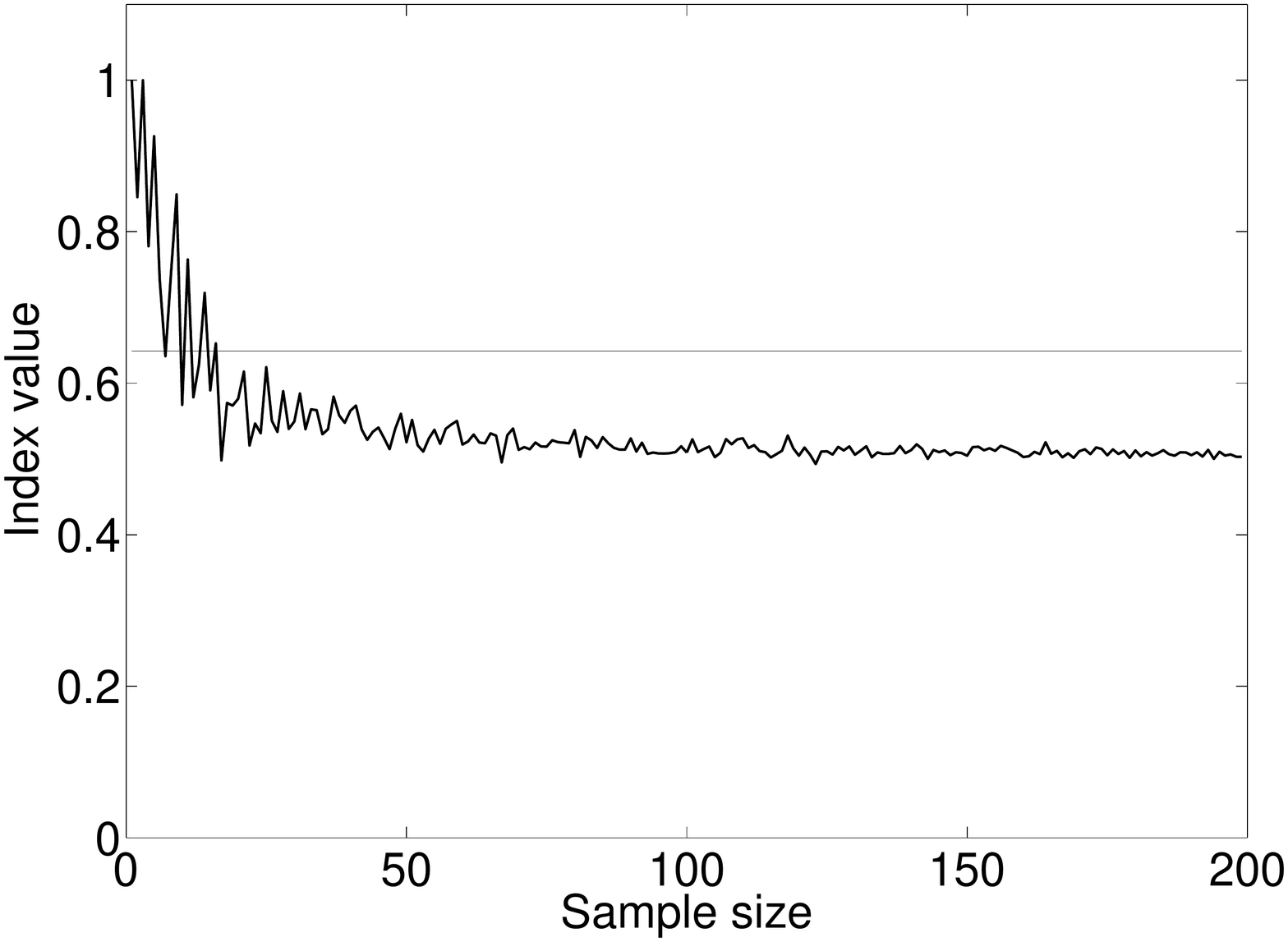}}
\hfill
\subfigure[$\alpha =  1.5$, window \textrm{(8,12]}, $\mathrm{I}_F(h)=1$.]{%
\includegraphics[height=0.35\textwidth]{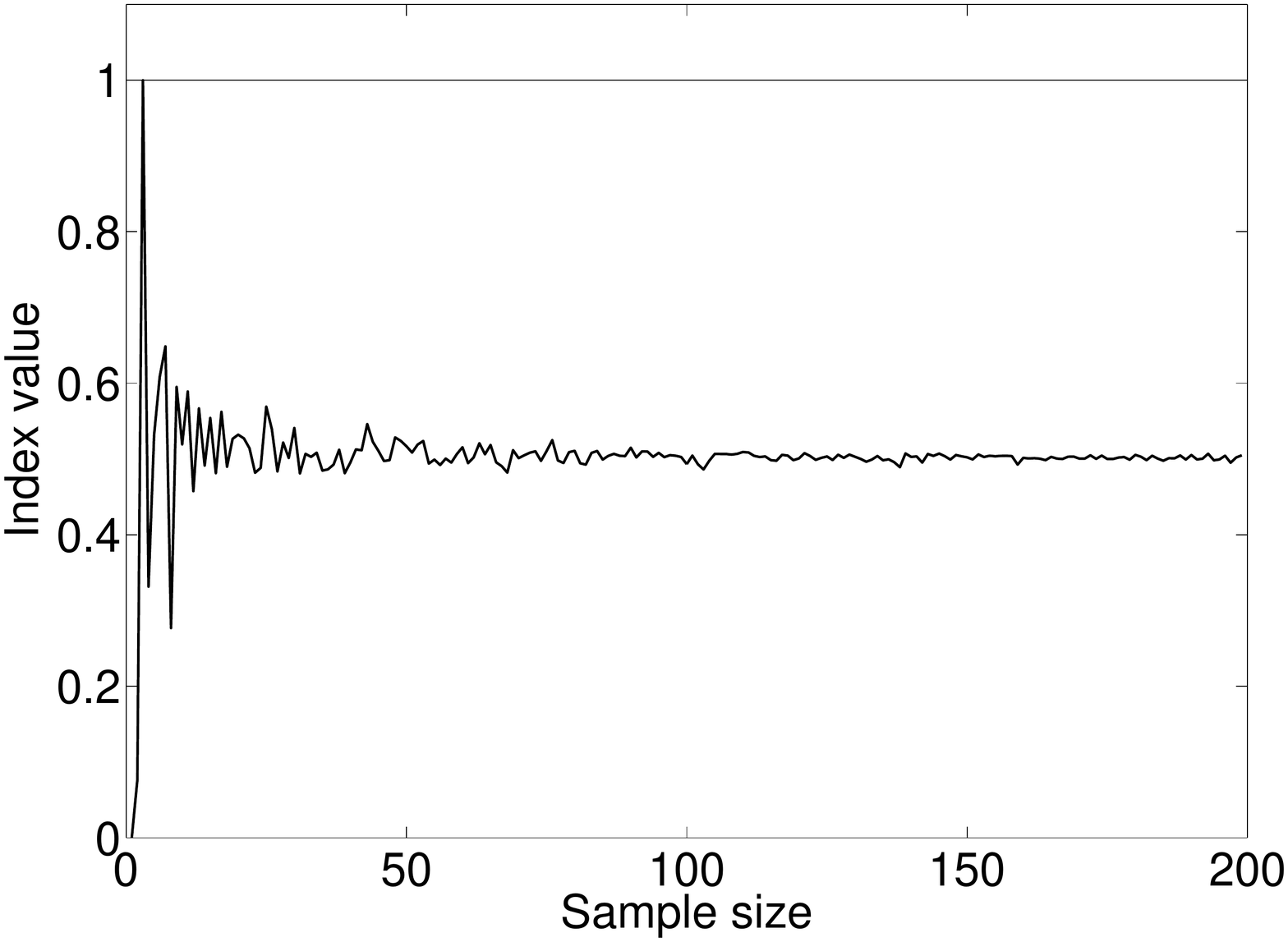}}
\hfill
\subfigure[$\alpha =  5$, window \textrm{(8,12]}, $\mathrm{I}_F(h)=1$.]{%
\includegraphics[height=0.35\textwidth]{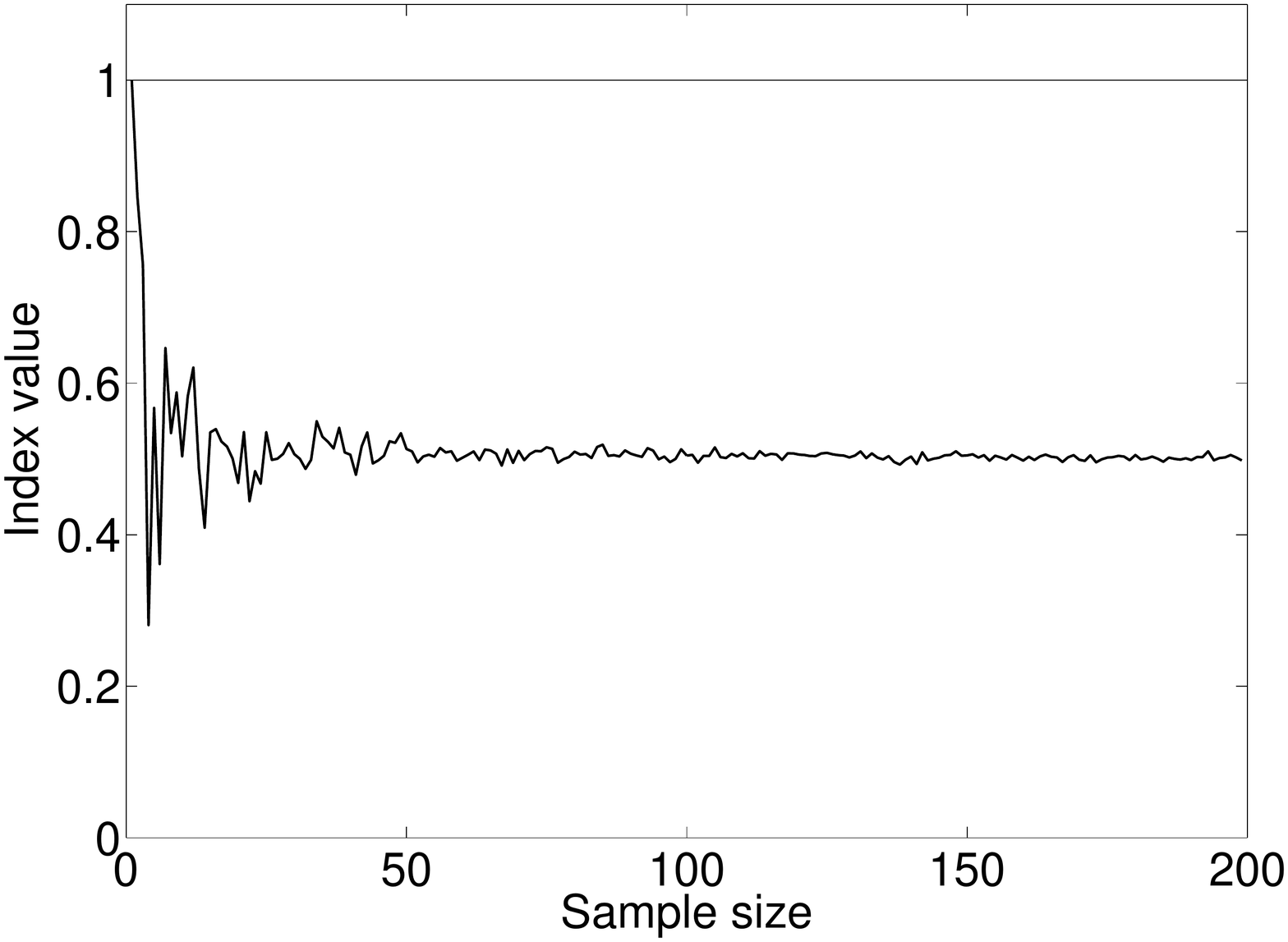}}
\caption{The index $\mathrm{I}_n$ for the transfer function $h(x)$ with $\varrho =  0$ and the added Gaussian noise with mean $0$ and standard deviation $0.1$, with the horizontal lines at the heights equal to the respective values of $\mathrm{I}_F(h)$.}
\label{fig-noise}
\end{figure}

We start numerical explorations with the case when a noise is added to the outputs. In this case the empirical index $\mathrm{I}_n$, which is defined by equation (\ref{random-index}), tends to $1/2$ instead of the true value $\mathrm{I}_F(h)$. We see this from Figure~\ref{fig-noise},
where the outputs $Y_{i,n}=h(X_{i:n})+\varepsilon_{i} $ have been generated using:
\begin{itemize}
\item
iid inputs $X_1,\dots , X_n$ that follow the Lomax distribution with the specified in the panels $\alpha $'s and the scale parameter $\beta=1$;
\item
the transfer function $h(x)$ given by formula (\ref{tf-1}) with
$\gamma=0.1$, $\delta =1$, and $\varrho=0$;
\item
additive noise made up of independent Gaussian variables $\varepsilon_{i}$ with means $0$ and standard deviations $\sigma=0.1$.
\end{itemize}
The horizontal lines in the panels of Figure~\ref{fig-noise} are at the heights of the corresponding values of $\mathrm{I}_F(h)$, which is defined by equation (\ref{ii-main-2}) and whose calculations we provide next.

Since in the current example we assume $\varrho=0$, the transfer function $h(x)$ is continuous on its entire domain of definition. Furthermore, as seen from equations (\ref{qf-t}) and (\ref{qf-x}), the quantile function of $X$ is continuous on the interval $(0,1)$. Hence, the number $m$ of jumps that show up in condition (C\ref{cond-2}) is zero, and so the index of increase given by equation (\ref{ii-main-2}) reduces to
\begin{align}
\mathrm{I}_F(h)
={\int_{0}^{1} H_{+}(u)\mathrm{d}u \over \int_{0}^{1}|H|(u)\mathrm{d}u}
&= {\int_{0}^{1} (h')_{+}\circ F^{-1}(u)\mathrm{d} F^{-1}(u) \over \int_{0}^{1}|h'|\circ F^{-1}(u)\mathrm{d} F^{-1}(u)}
\notag
\\
&= {\int_{a}^{b} (h')_{+}(x)\mathrm{d} x \over \int_{a}^{b}|h'|(x)\mathrm{d} x}.
\label{I-cont}
\end{align}
To calculate the index $\mathrm{I}_F(h)$, therefore, we need to know the derivative of the transfer function $h(x)$, which, when $\varrho=0$, is equal to
\begin{equation}
h'(x)={ (\gamma x^2 + \delta)^2-2\delta x \over x^4} e^{\gamma x -\delta/x}  .
\label{tf-deriv}
\end{equation}
(This derivative has played a major role in reliability-engineering modelling by Bebbington et al.~(2007).) Numerical integration on the right-hand side of equation (\ref{I-cont}) yields $\mathrm{I}_F(h)\approx 0.6424$ when the window is $(0,2]$ (the top two panels of Figure~\ref{fig-noise}) and $\mathrm{I}_F(h)=1$ when the window is $(8,12]$ (the bottom two panels of Figure~\ref{fig-noise}). This concludes our illustration of the case with the additive noise.

Next we explore statement \eqref{ii-main-2} within the contexts of Theorems \ref{Thm-2} and \ref{Thm-3}, and in each window (\ref{par-4}).

\subsubsection{Window $(0,2]$}

\begin{figure}[h!]
\centering
\subfigure[$\alpha=1.5$, $\varrho=0.5$]{%
\includegraphics[height=0.35\textwidth]{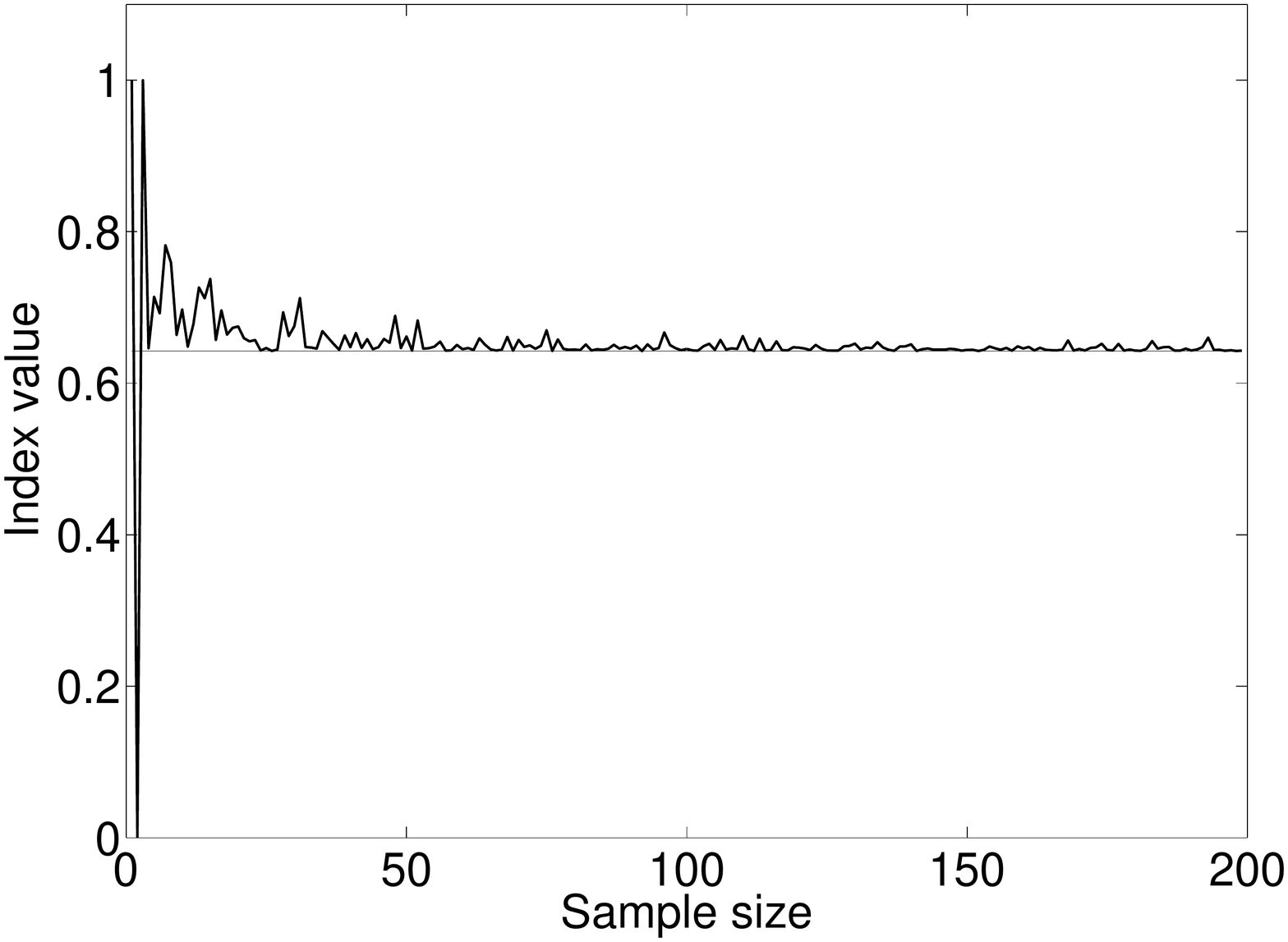}}
\hfill
\subfigure[$\alpha=5$, $\varrho=0.5$]{%
\includegraphics[height=0.35\textwidth]{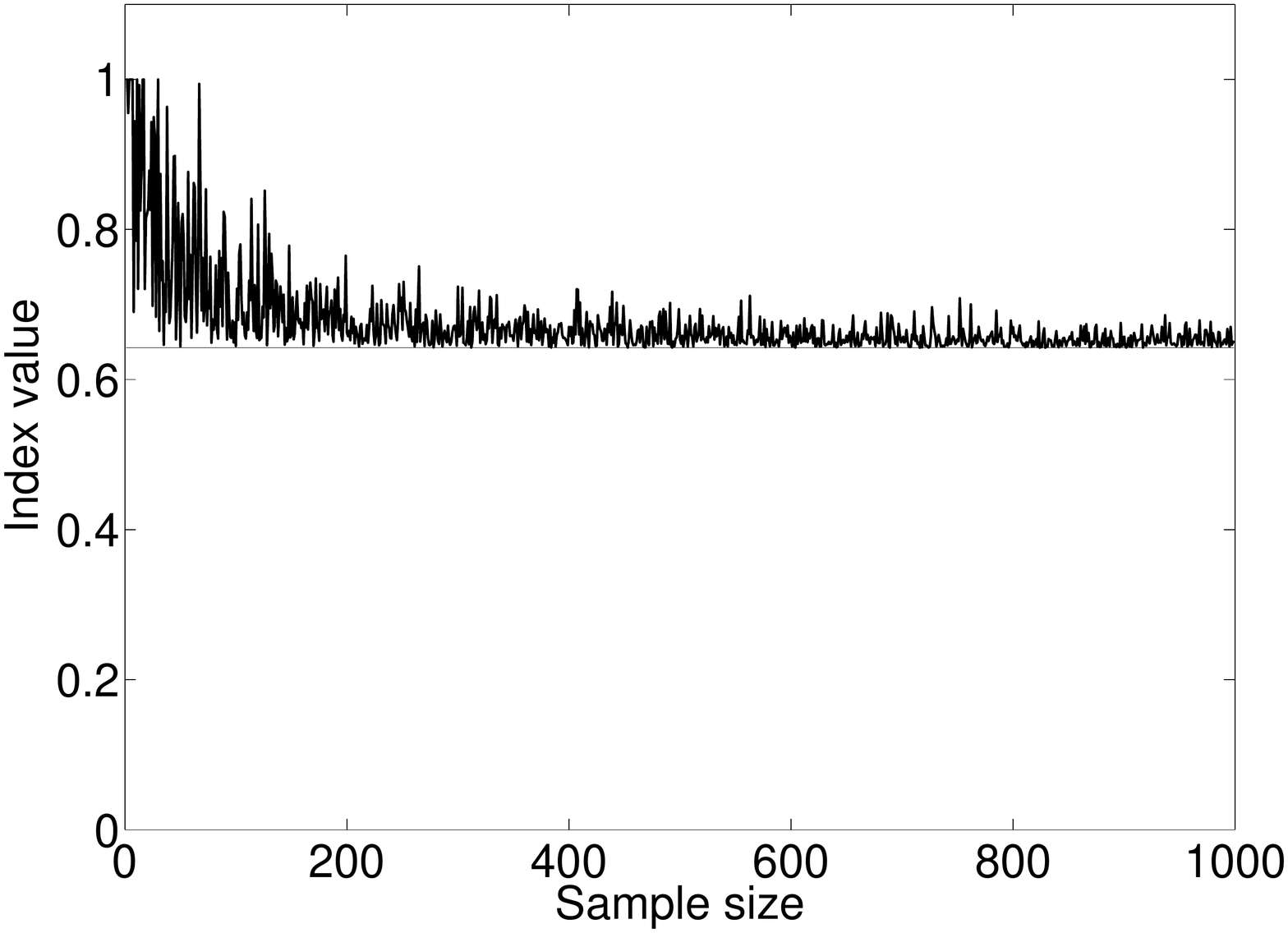}}
\\
\subfigure[$\alpha=1.5$, $\varrho=0$]{%
\includegraphics[height=0.35\textwidth]{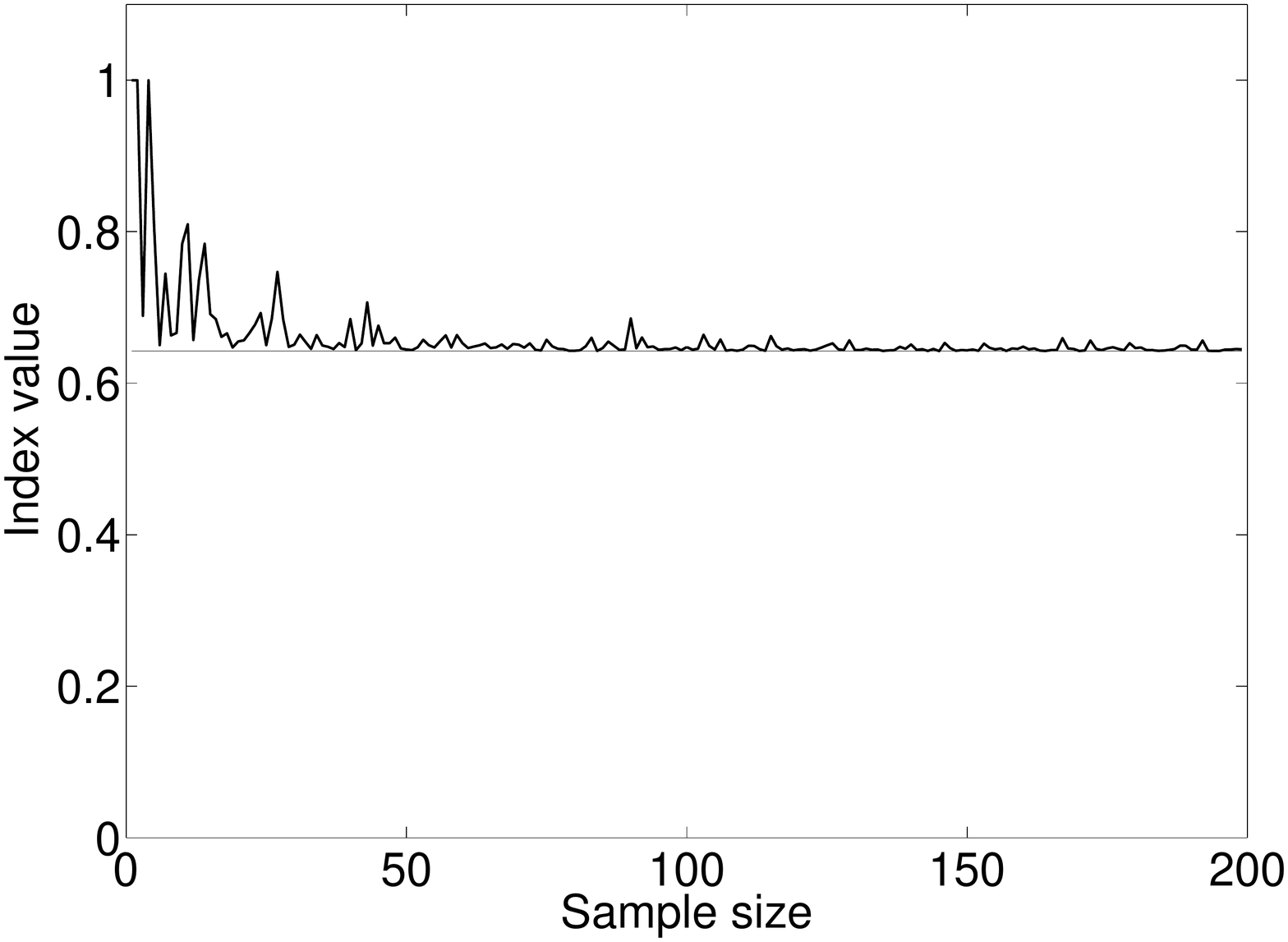}}
\hfill
\subfigure[$\alpha=5$, $\varrho=0$]{%
\includegraphics[height=0.35\textwidth]{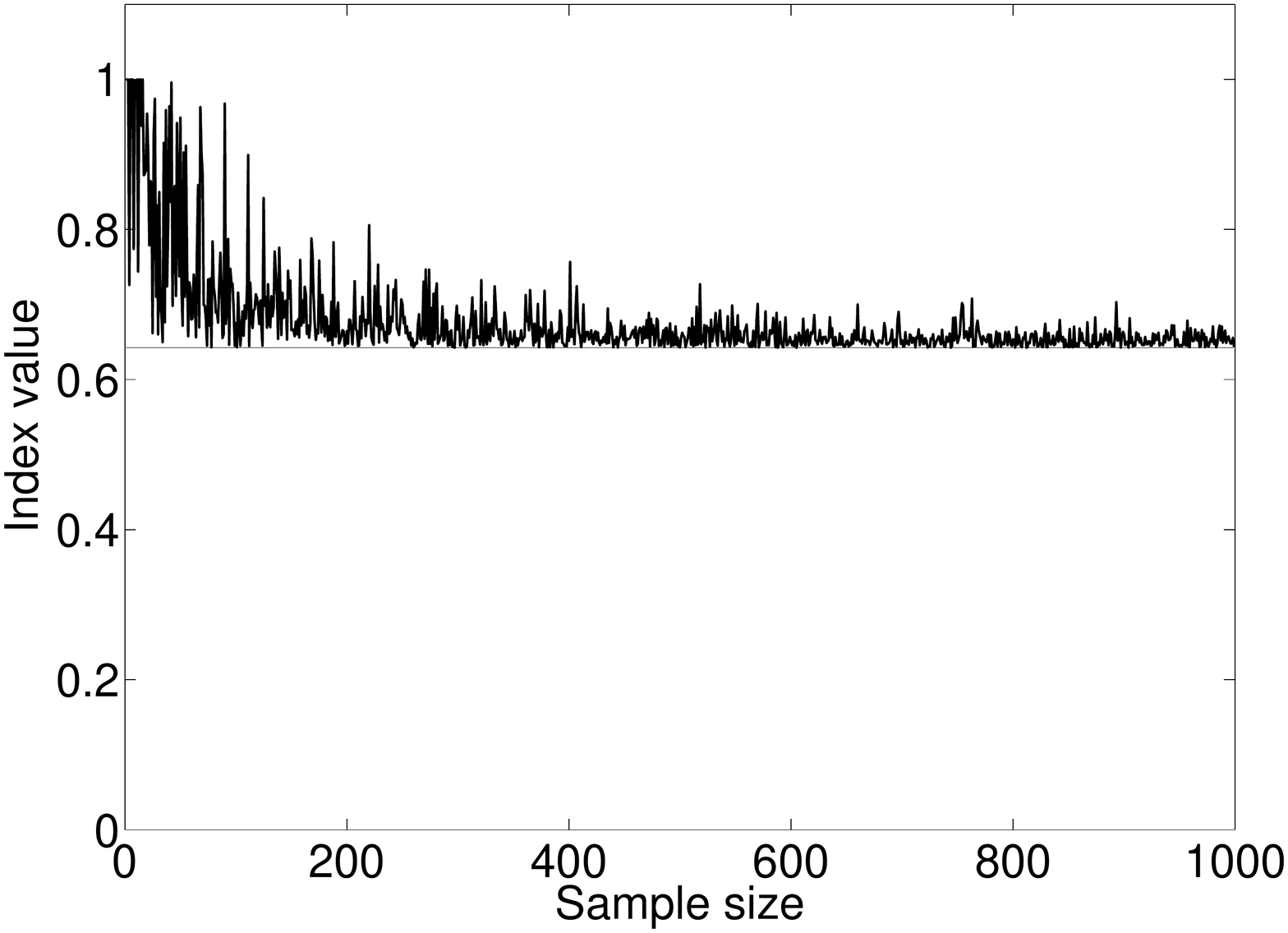}}
\\
\subfigure[$\alpha=1.5$, $\varrho=-0.5$]{%
\includegraphics[height=0.35\textwidth]{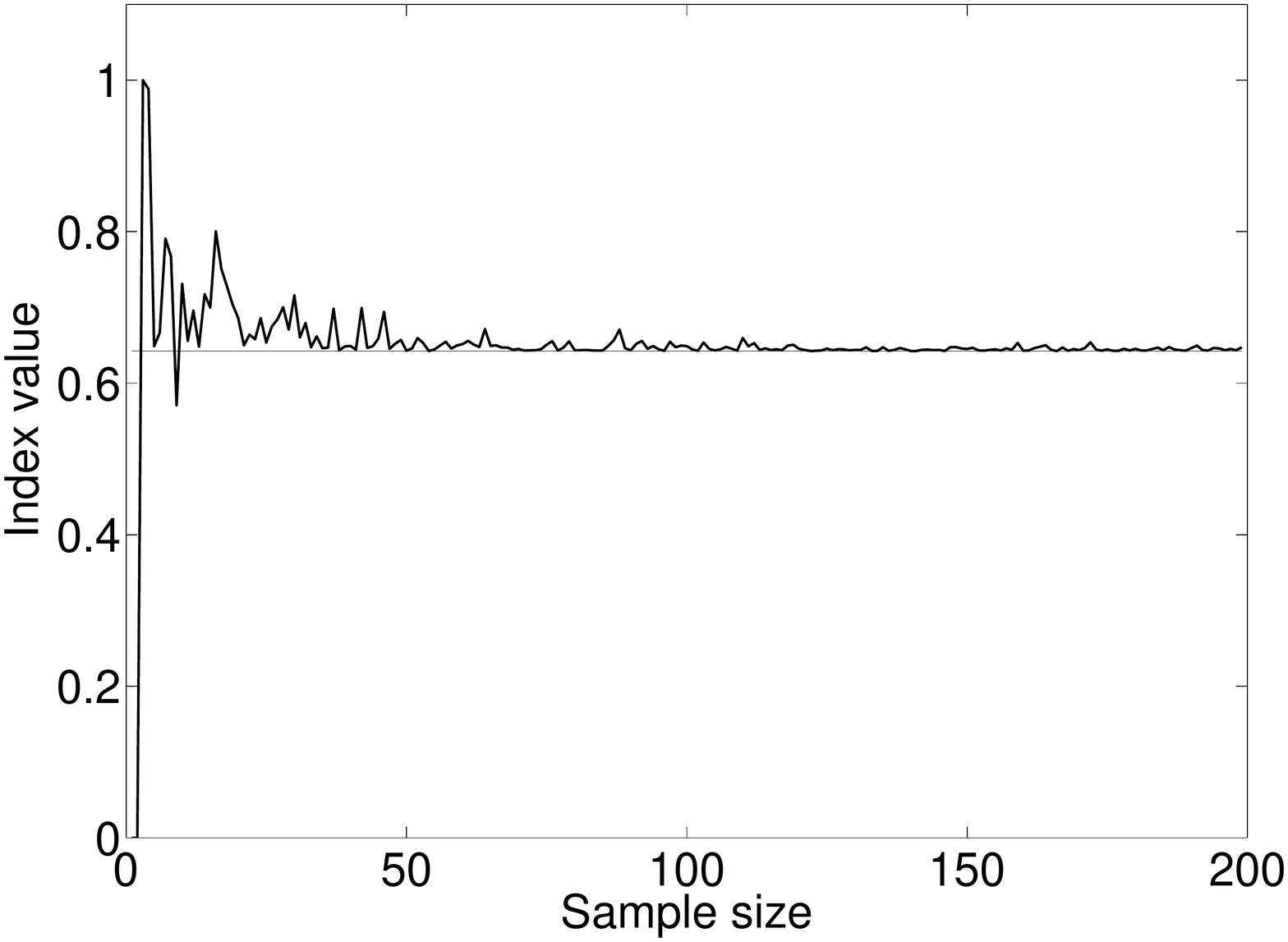}}
\hfill
\subfigure[$\alpha=5$, $\varrho=-0.5$]{%
\includegraphics[height=0.35\textwidth]{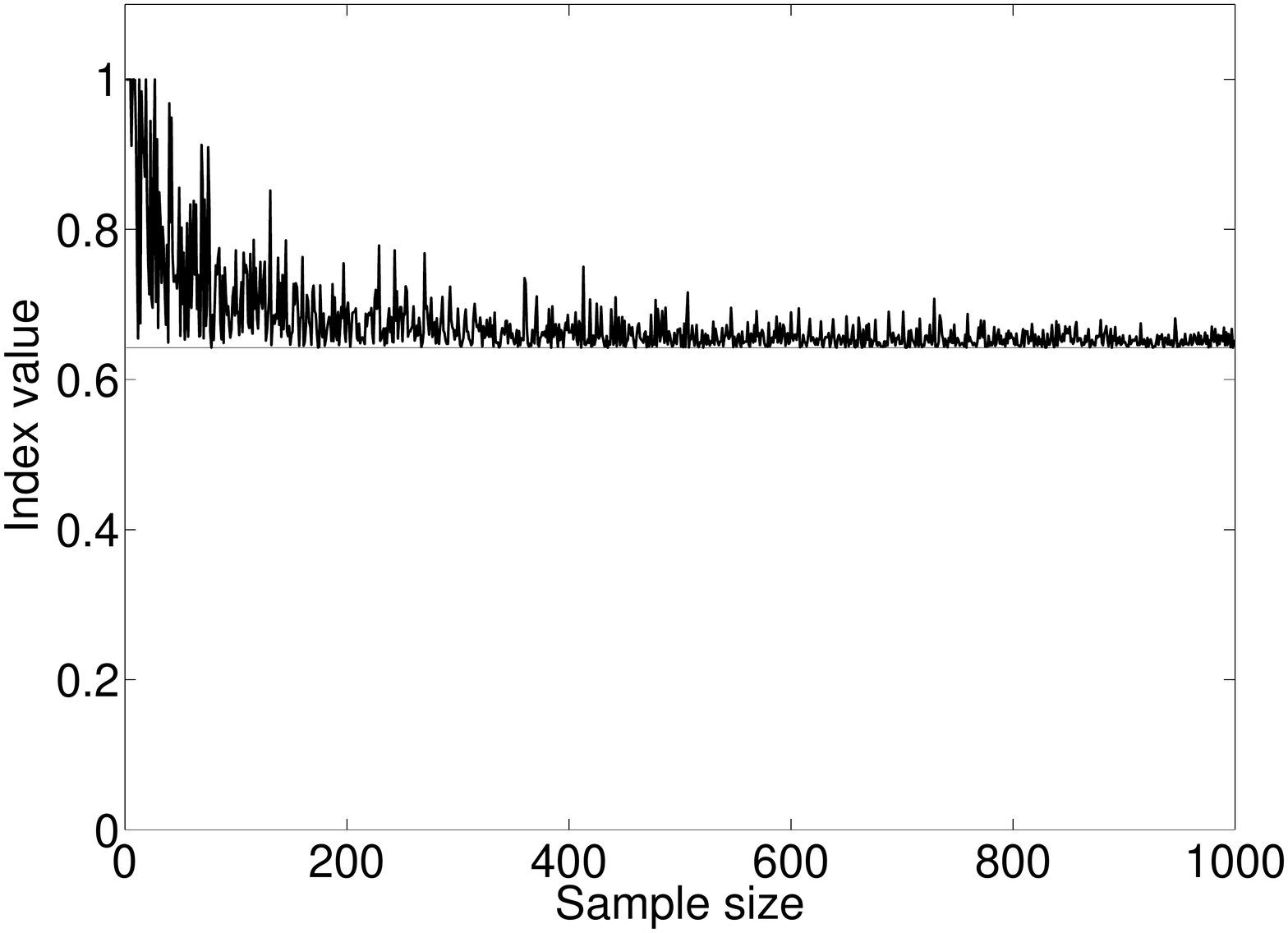}}
\caption{Performance of $\mathrm{I}_n$ in the window $(0,2]$, with the horizontal lines at the height of $\mathrm{I}_F(h)=0.6424$.}
\label{fig-win02}
\end{figure}

In the window $(0,2]$, the function $h(x)$ and thus the indices $\mathrm{I}_n$ and $\mathrm{I}_F(h)\approx 0.6424$ are not sensitive to the parameter $\varrho $ values, which regulate the sizes of jumps at $x=10$. Hence, in all the six panels of Figure~\ref{fig-win02}
the theoretical index is (cf.~equation (\ref{I-cont}))
\begin{equation}
\mathrm{I}_F(h)
={\int_{0}^{2} (h')_{+}(x)\mathrm{d} x \over \int_{0}^{2}|h'|(x)\mathrm{d} x}\approx 0.6424,
\label{index-i-00}
\end{equation}
which is the height of the horizontal line in each panel. We see that the empirical index $\mathrm{I}_n$ converges to the true value $0.6424$ in all the panels, albeit much slower when $\alpha=5$ (the three right-hand panels) than when $\alpha=1.5$ (the three left-hand panels); notice the different scales in the left- and right-hand panels. At first sight, these convergence results may look surprising because lighter-tailed distributions usually exhibit better statistical properties, but the current context is quite different: we are feeding into the filter not the values of the original random variable $T$ but those of its truncated version $X=\langle T \mid T\in (0,2] \rangle $. We next explain, aiming primarily at intuition, how this affects the performance of the estimator $\mathrm{I}_n$.

\begin{figure}[h!]
\centering
\subfigure[Pdf $f(x)$]{%
\includegraphics[height=0.35\textwidth]{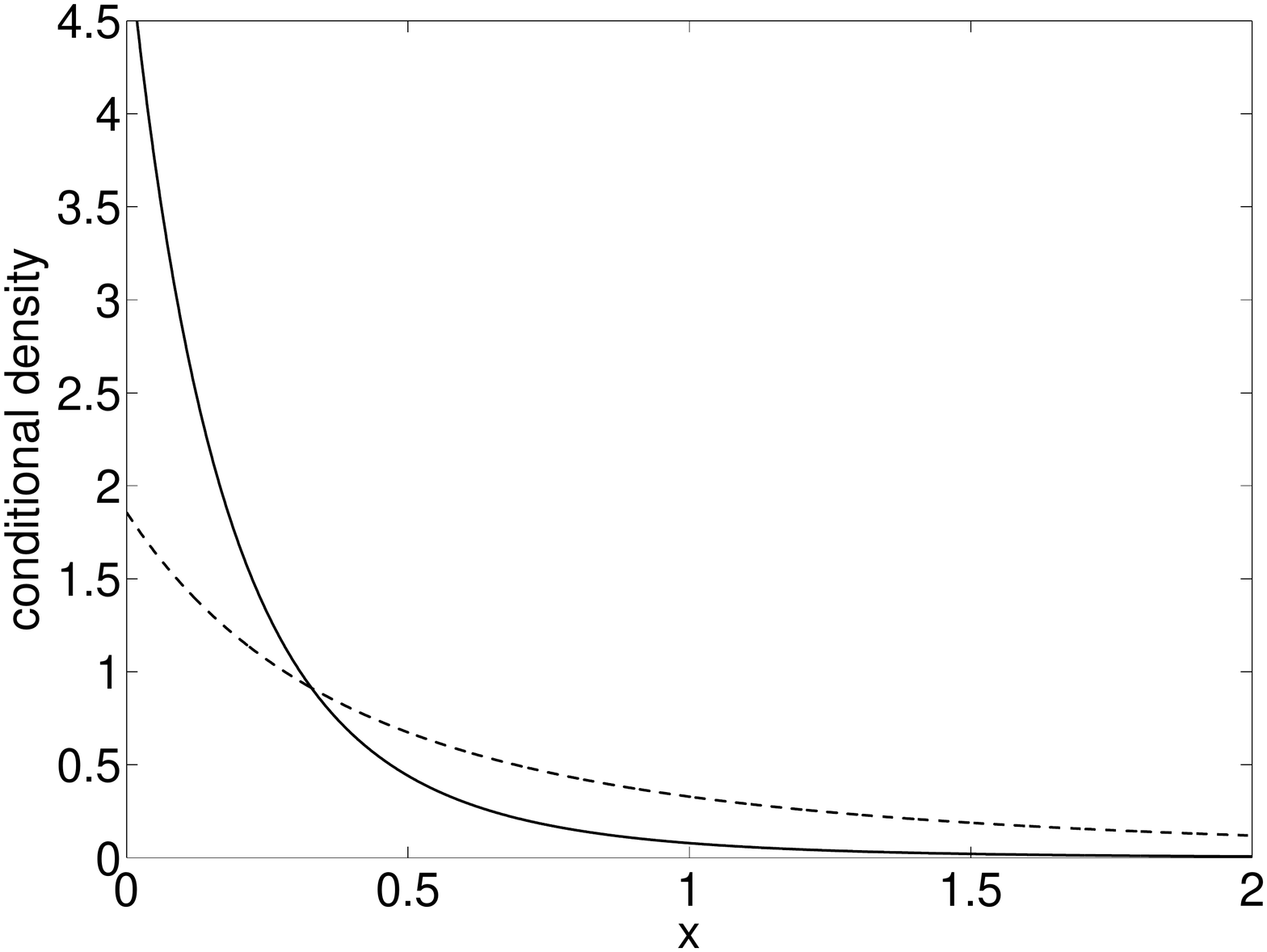}}
\hfill
\subfigure[Function $H(u)$]{%
\includegraphics[height=0.35\textwidth]{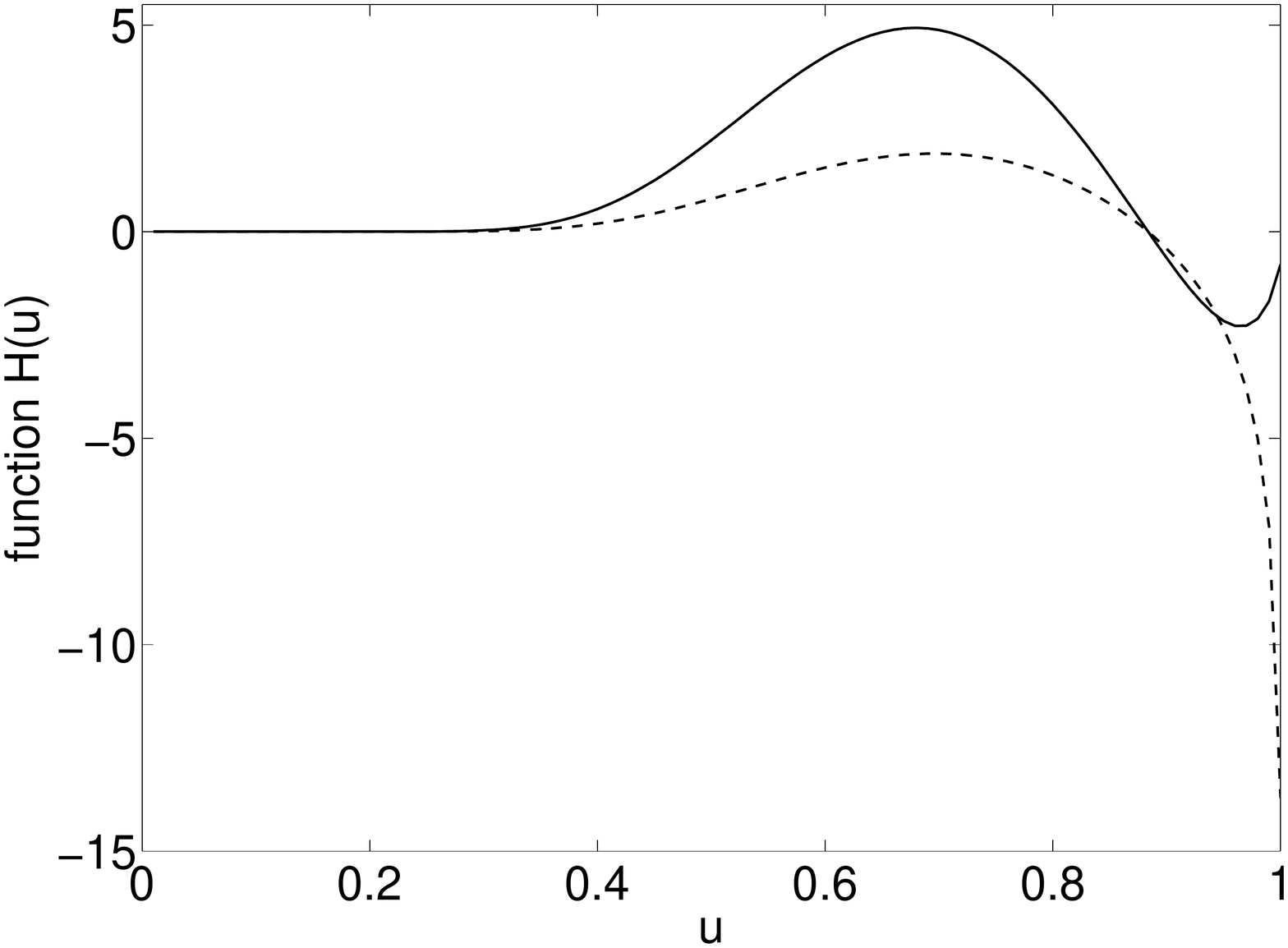}}
\caption{Two functions corresponding to $X=\langle T \mid T\in (0,2]\rangle $ when $\alpha=1.5$ (solid) and $\alpha=5$ (dashed).}
\label{pdf-02}
\end{figure}

Namely, loosely speaking, the less the function $H(u)$ jolts on the interval $(0,1)$, the faster the rate of convergence of $\mathrm{I}_n$ to $\mathrm{I}_F(h)$. For example, if the transfer function $h(x)$ does not change much in the window $(a,b]$, whatever it might be, but the pdf $f(x)$ on one part of the window is very large but on another part is nearly zero, then the function $H(u)$ has a considerable jolt. Going back to the currently explored window $(0,2]$, we see from the right-hand panel of Figure~\ref{pdf-02}
that there is some wiggling in the case $\alpha=1.5$, but it is less pronounced than in this case $\alpha=5$, when the function $H(u)$ dips well below zero when it approaches the endpoint $u=1$. Hence, we see a slower rate of convergence in the latter case. Yet, in neither case is the convergence of $\mathrm{I}_n$ to $\mathrm{I}_F(h)$ completely derailed, unlike in some other windows that we shall explore below. In those cases the index $\mathrm{I}_n$ would not even get close to $\mathrm{I}_F(h)$ due to the low density $f(x)$ values and thus scarcity of data in certain regions of the input interval, thus making fluctuations of the function $H(u)$ large in some regions and low in other ones.

\subsubsection{Window $(8,12]$}

\begin{figure}[h!]
\centering
\subfigure[Pdf $f(x)$]{%
\includegraphics[height=0.35\textwidth]{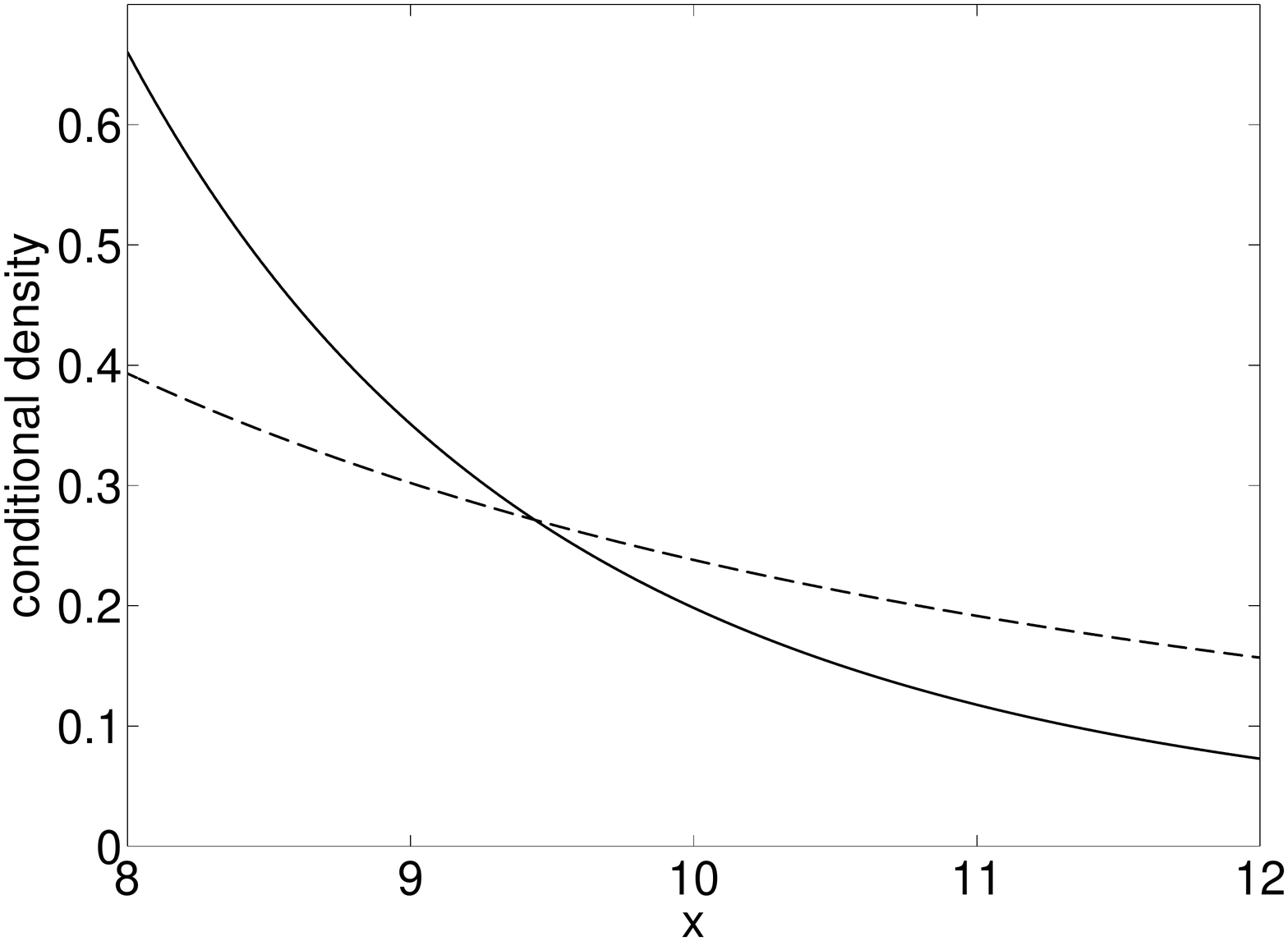}}
\hfill
\subfigure[Function $H(u)$ when $\varrho = - 0.5$]{%
\includegraphics[height=0.35\textwidth]{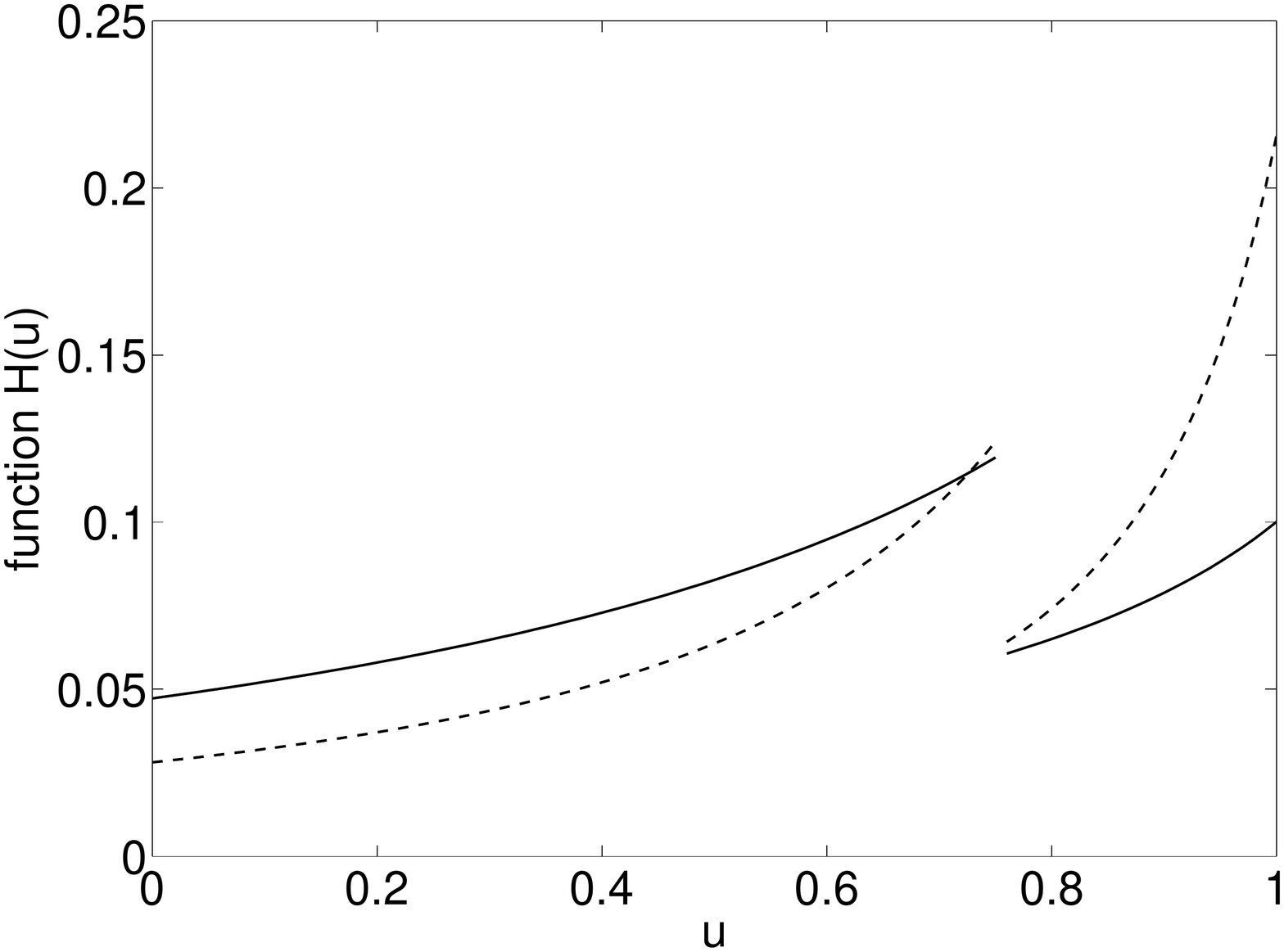}}
\\
\subfigure[Function $H(u)$ when $\varrho = 0$]{%
\includegraphics[height=0.35\textwidth]{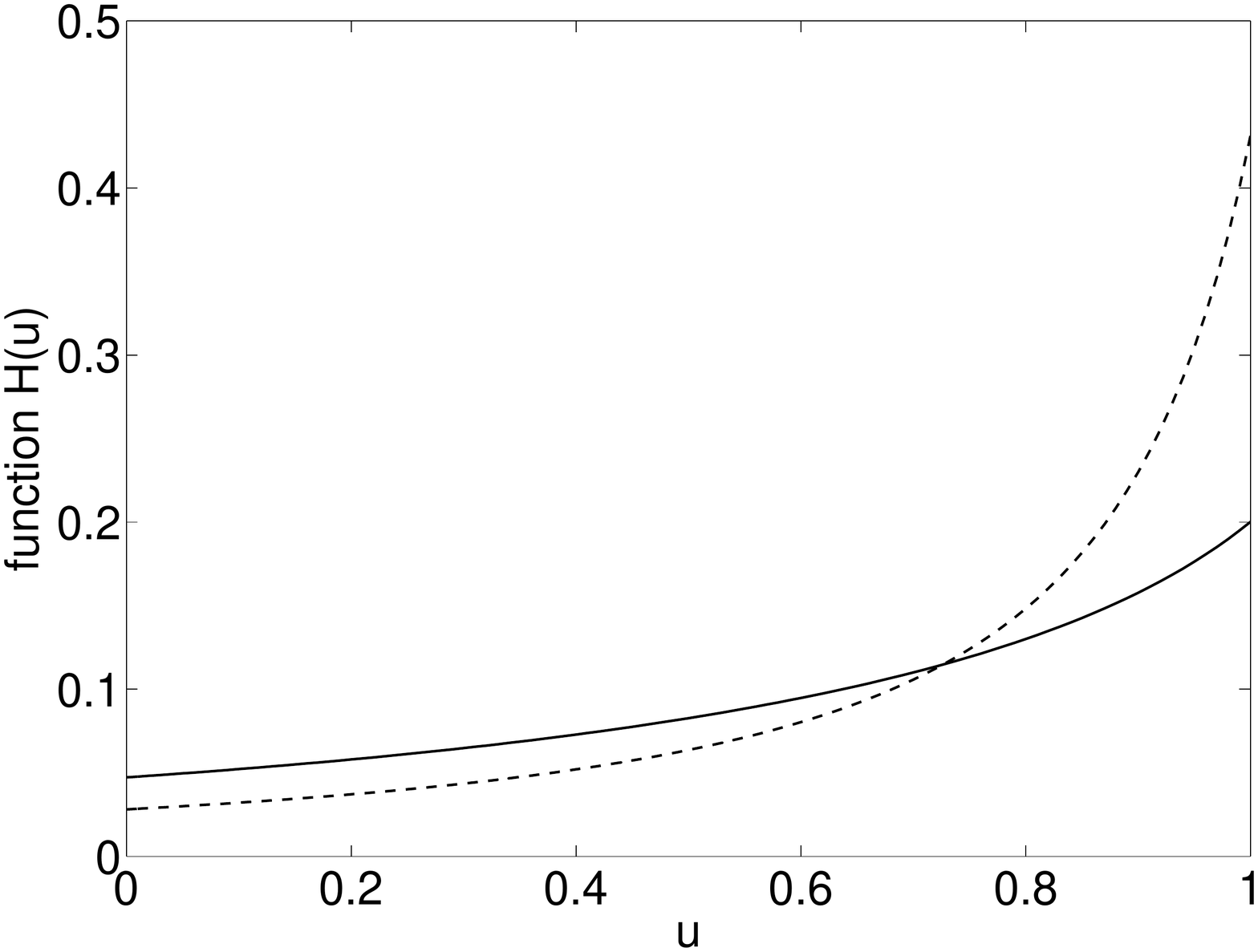}}
\hfill
\subfigure[Function $H(u)$ when $\varrho = 0.5$]{%
\includegraphics[height=0.35\textwidth]{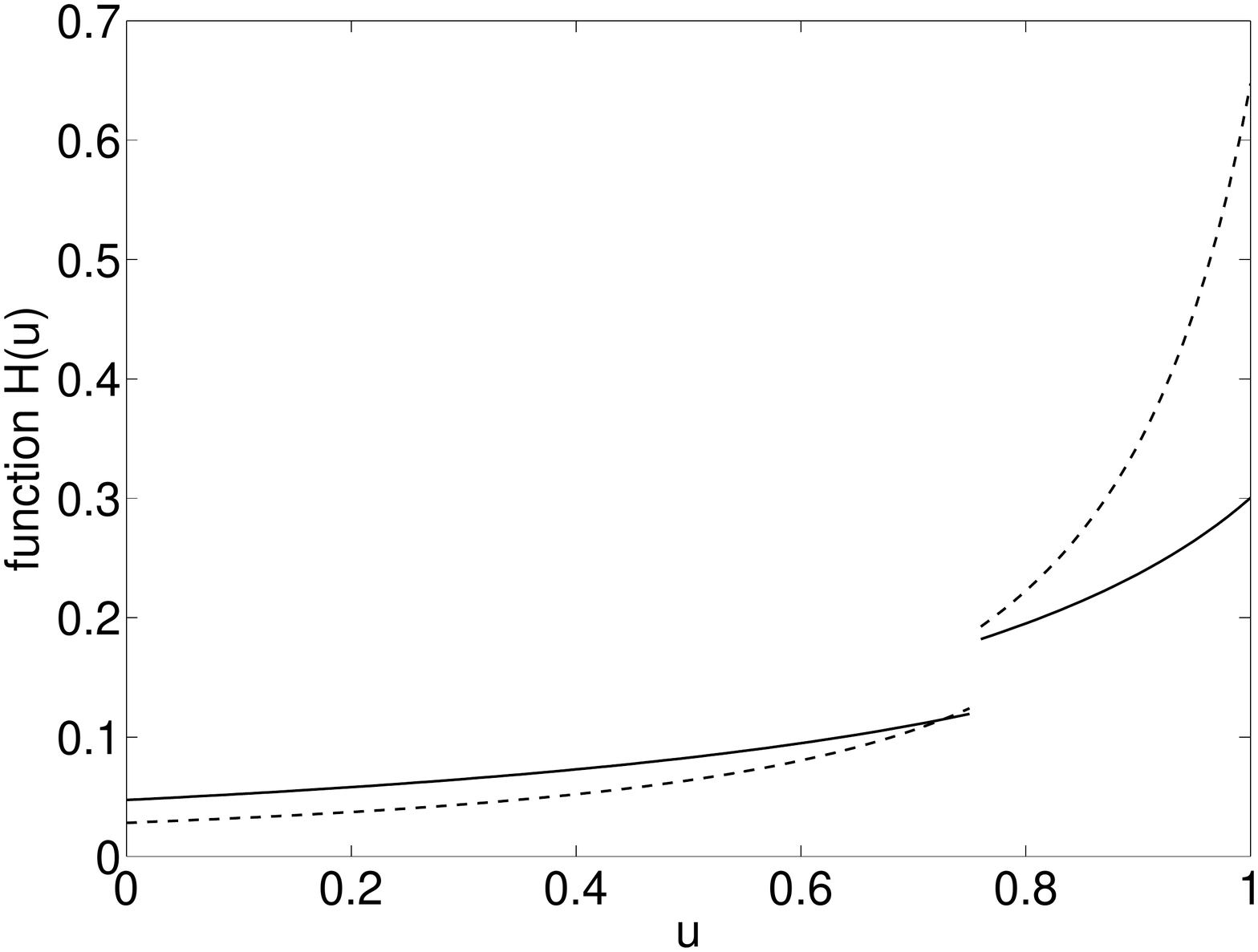}}
\caption{Pdf $f(x)$ and the function $H(u)$ corresponding to $X= \langle T \mid T\in (8,12] \rangle $ when $\alpha=1.5$ (solid) and $\alpha=5$ (dashed).}
\label{pdf-812}
\end{figure}

\begin{figure}[h!]
\centering
\subfigure[$\alpha=1.5$, $\varrho=-0.5$, $\mathrm{I}_F(h)\approx 0.3459$]{%
\includegraphics[height=0.35\textwidth]{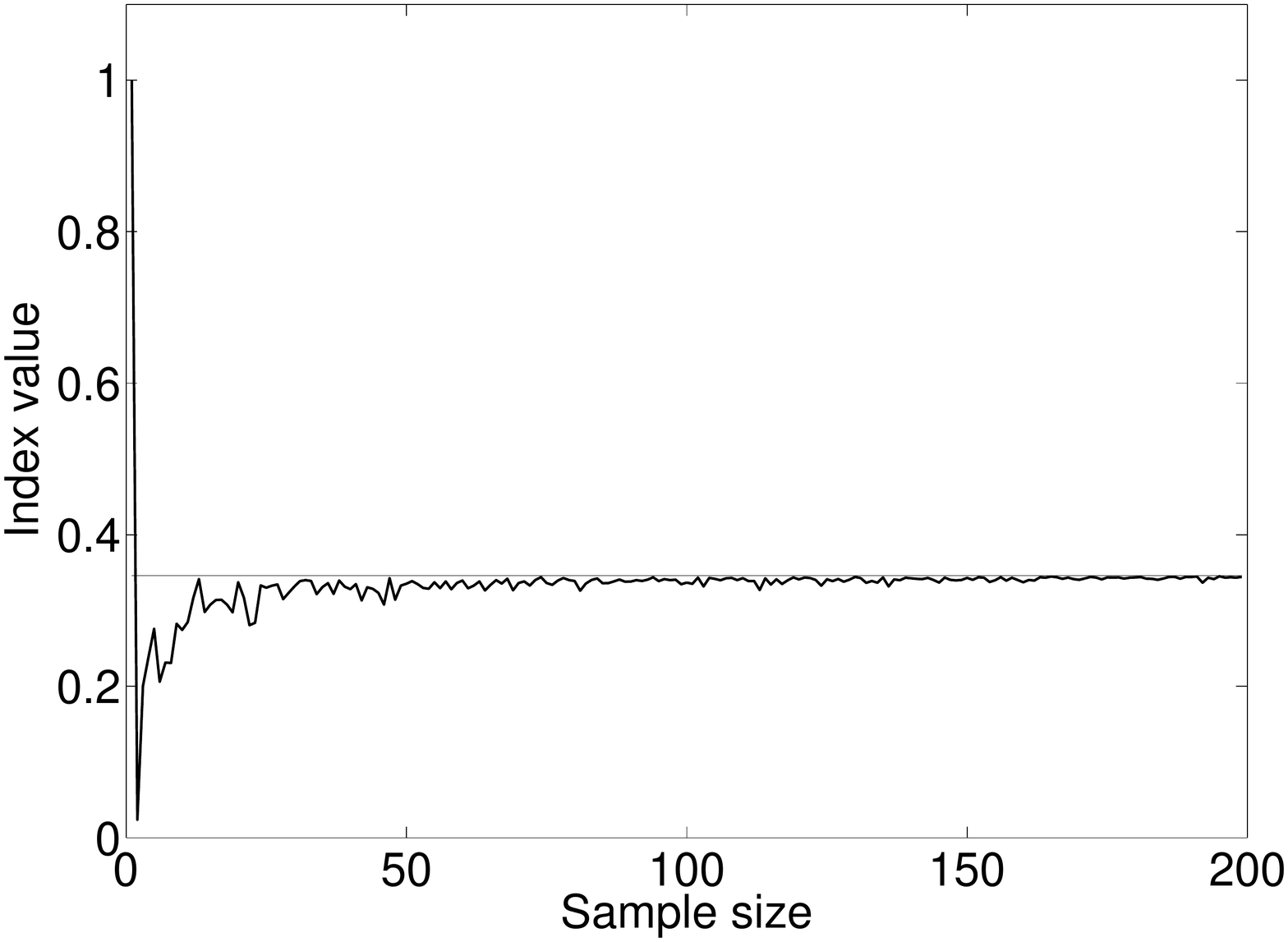}}
\hfill
\subfigure[$\alpha=5$, $\varrho=-0.5$, $\mathrm{I}_F(h)\approx 0.3459$]{%
\includegraphics[height=0.35\textwidth]{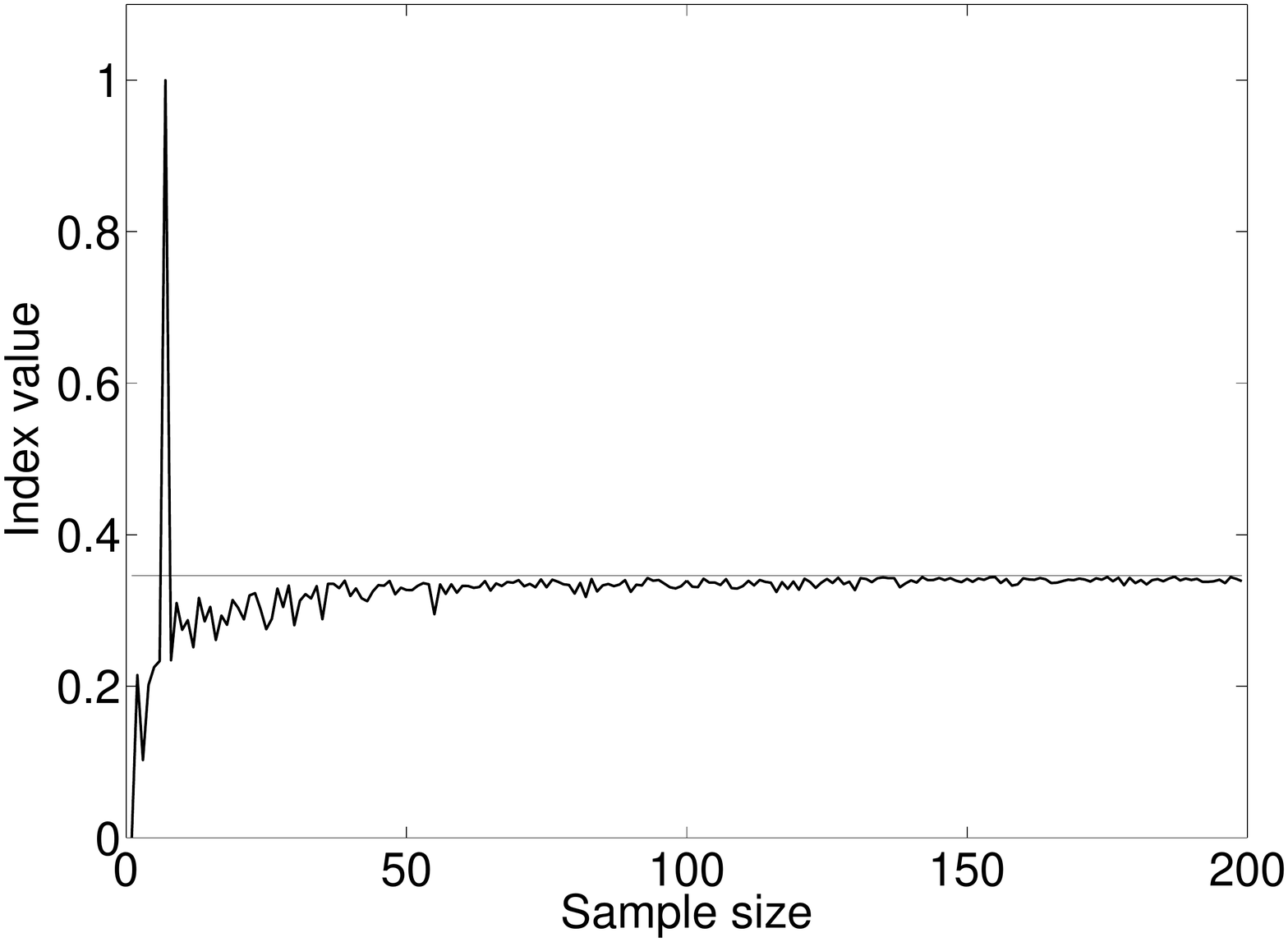}}
\caption{Performance of $\mathrm{I}_n$ in the window $(8,12]$, with the horizontal lines at the height of $\mathrm{I}_F(h)=0.3459$.}
\label{fig-win812}
\end{figure}

The window $(8,12]$ contains a regime switching whenever $\varrho \neq 0$, and we choose the values $\varrho = \pm 0.5$ to illustrate this case, in addition to the continuous case $\varrho = 0$. Note at the outset that the window $(8,12]$ is deeper into the Lomax tail than in the previously considered window $(0,2]$, and thus, as seen from Figure~\ref{pdf-812},
the function $H(u)$ exhibits, arguably, more comparable deviations from constancy when $\alpha=1.5$ and $\alpha=5$ than those that we earlier saw in the window $(0,2]$. As a consequence, in Figure~\ref{fig-win812},
which concerns only the case $\varrho = -0.5$, we see very similar convergence patters of the estimator $\mathrm{I}_n$ to $\mathrm{I}_F(h)$ when $\alpha=1.5$ and $\alpha=5$.

As to calculating the theoretical value $\mathrm{I}_F(h)$, we first note that since the transfer function $h(x)$ has a drop at the point $x=10$, we therefore have $\tau_1=F(10)$ and thus $\Delta h_1=h\circ F^{-1}(\tau_1+0)-h\circ F^{-1}(\tau_1)=\varrho e^{9/10} 11/100\approx -0.1353$. Consequently,
\[
\mathrm{I}_F(h)
={\int_{8}^{12} (h')_{+}(x)\mathrm{d} x \over |\Delta h_1|+\int_{8}^{12}|h'|(x)\mathrm{d} x}
\approx {0.0715 \over 0.1353+0.0715}\approx 0.3459.
\]
The values of the above integrals in the numerator and denominator are equal, because $h'(x)$ is non-negative to the left and to the right of the drop point $x=10$.

The reason we have depicted only the case $\varrho=-0.5$ in Figure~\ref{fig-win812} is that when $\varrho=0$ and $\varrho=0.5$, the empirical and theoretical indices of increase are equal to $1$, due to the fact that the transfer function $h(x)$ is increasing on the window $(8,12]$. In more detail, when $\varrho=0$, the transfer function $h(x)$ is increasing on the window $(8,12]$, and thus
\[
\mathrm{I}_F(h)
={\int_{8}^{12} (h')_{+}(x)\mathrm{d} x \over \int_{8}^{12}|h'|(x)\mathrm{d} x}=1.
\]
When $\varrho=0.5$, the transfer function is increasing on the window $(8,12]$ but has a jump at $x=10$. Hence, according to definition (\ref{ii-main-2}) of $\mathrm{I}_F(h)$, we have
\[
\mathrm{I}_F(h)=\frac{(\Delta h_1)_{+}+\int_{8}^{12} (h')_{+}(x)\mathrm{d} x }
{|\Delta h_1|+\int_{8}^{12}|h'|(x)\mathrm{d} x}=1,
\]
because of two facts: first, $\Delta h_1=\varrho e^{9/10} 11/100\approx 0.1353$ is a positive quantity, and second, the two integrals have the same values.

\subsubsection{Window $(0,20]$}
\label{section-323}

\begin{figure}[h!]
\centering
\subfigure[$\alpha=1.5$, $\varrho=0.5$, $H(1) = 144.1699$]{%
\includegraphics[height=0.35\textwidth]{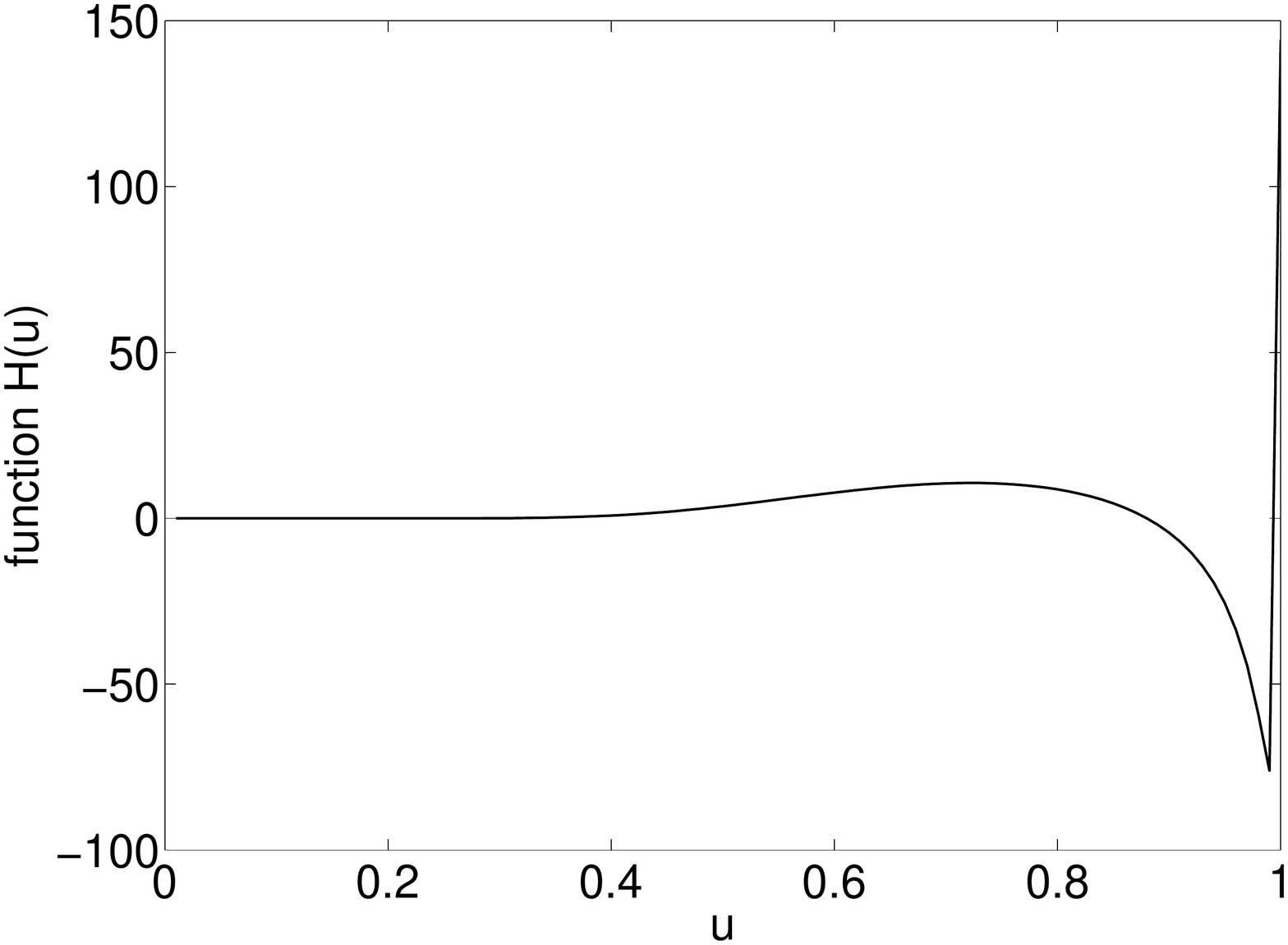}}
\hfill
\subfigure[$\alpha=5$, $\varrho=0.5$, $H(1) =  18.548  \times 10^5$]{%
\includegraphics[height=0.35\textwidth]{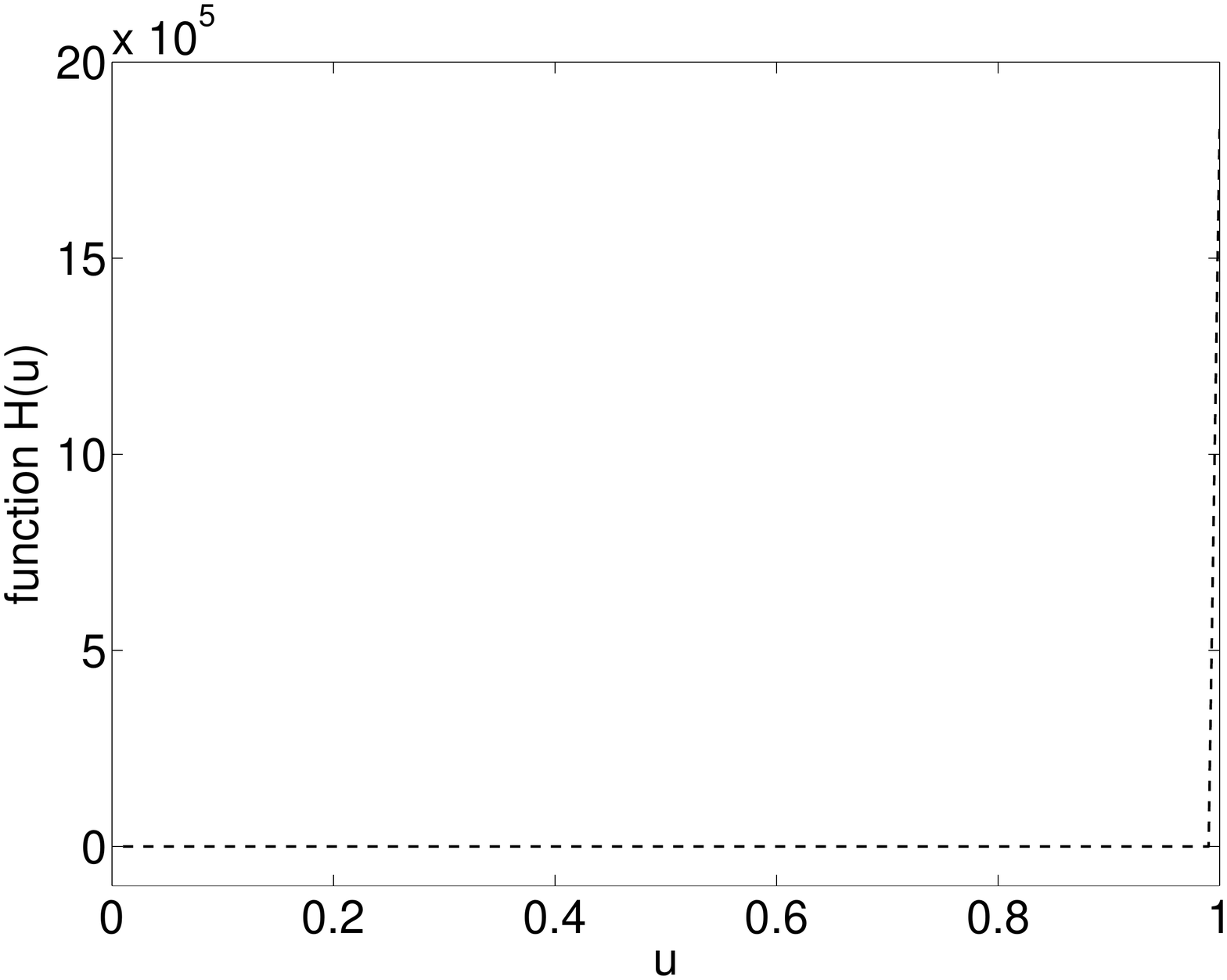}}
\\
\subfigure[$\alpha=1.5$, $\varrho=0$, $H(1) =   96.1133$]{%
\includegraphics[height=0.35\textwidth]{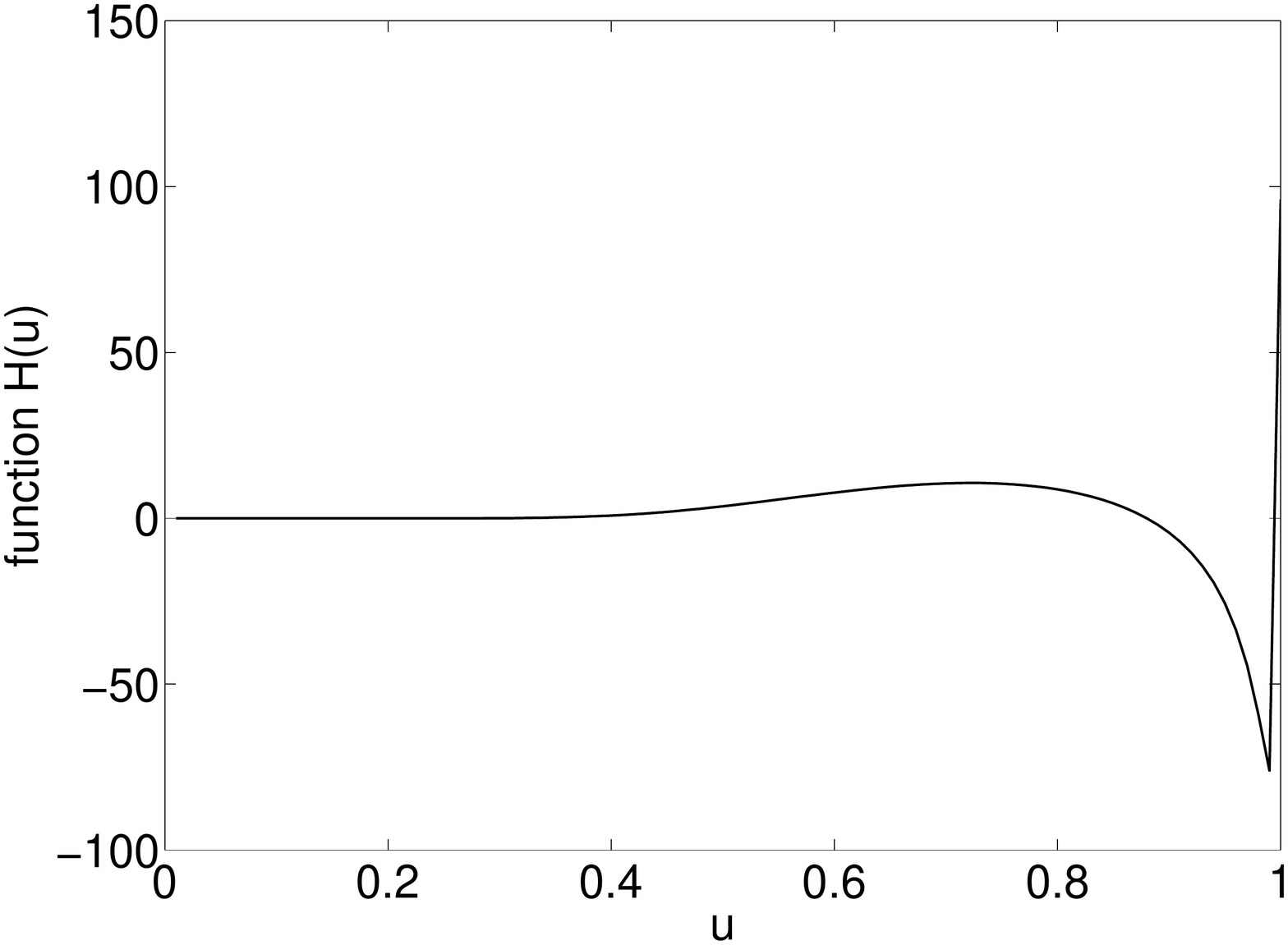}}
\hfill
\subfigure[$\alpha=5$, $\varrho=0$, $H(1) =  12.365  \times 10^5$]{%
\includegraphics[height=0.35\textwidth]{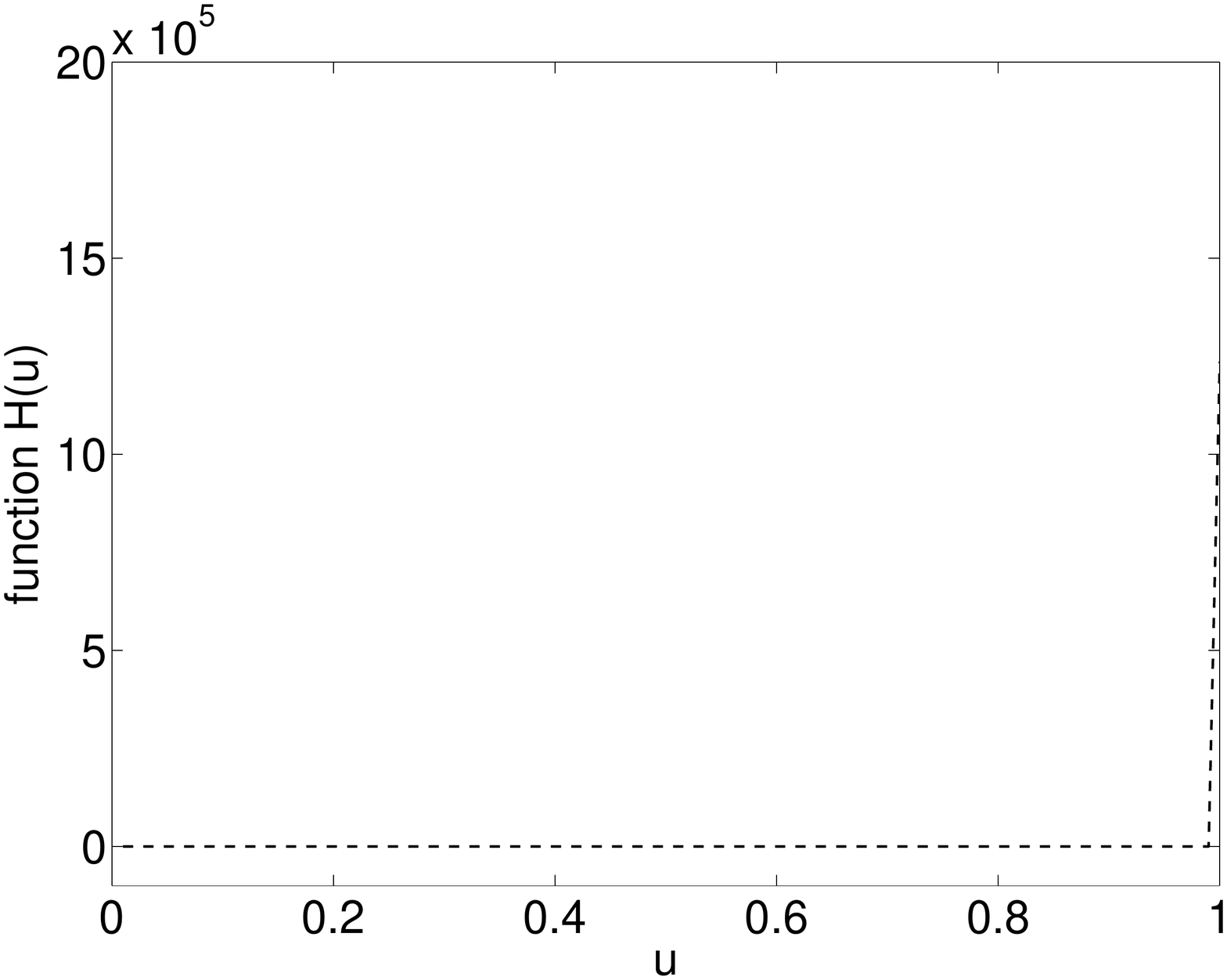}}
\\
\subfigure[$\alpha=1.5$, $\varrho=-0.5$, $H(1) = 48.0566$]{%
\includegraphics[height=0.35\textwidth]{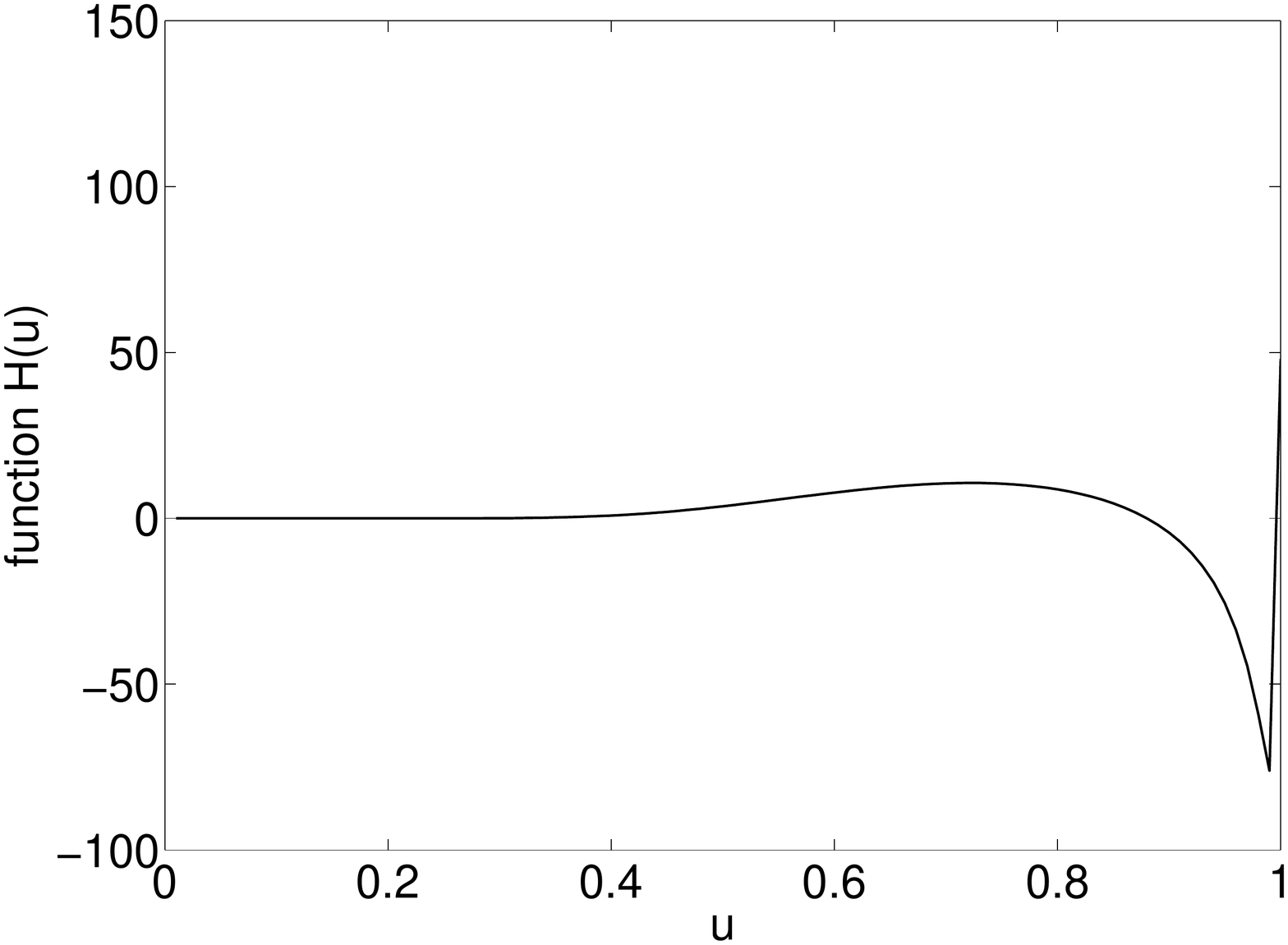}}
\hfill
\subfigure[$\alpha=5$, $\varrho=-0.5$, $H(1) =  6.183 \times 10^5$]{%
\includegraphics[height=0.35\textwidth]{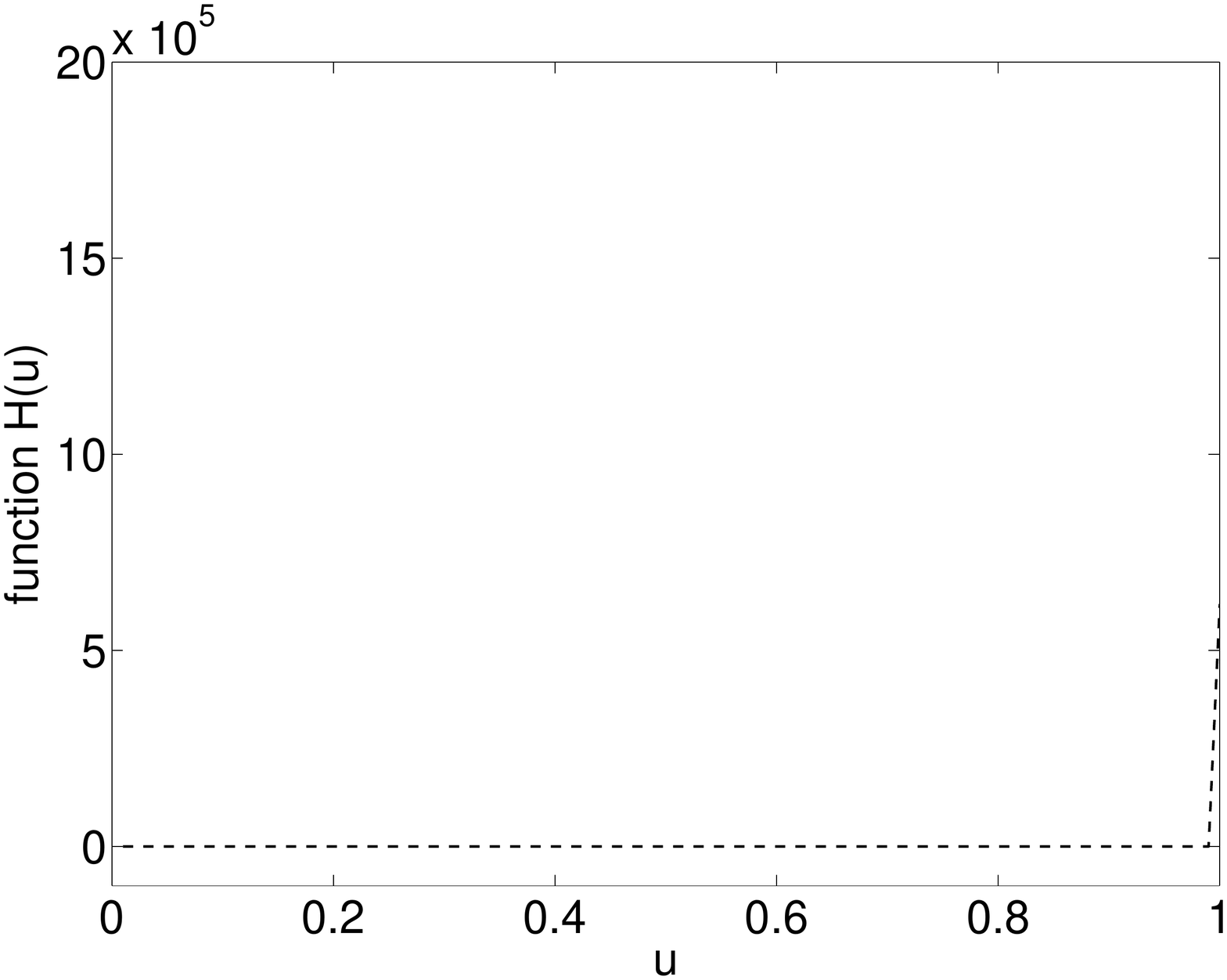}}
\caption{Function $H(u)$ corresponding to $X= \langle T \mid T\in (0,20] \rangle $; note the decreasing (from top to bottom) values of $H(1)$ and thus the decreasing sizes of the jolts near $u=1$.}
\label{pdf-020}
\end{figure}

The window $(0,20]$ is much winder than the previous two ones, and thus brings up the issue of sufficient data in every region of the window. The scarcity of data results in bigger jolts of the function $H(u)$, which we have drawn in Figure \ref{pdf-020}
under the two parameter values $\alpha=1.5$ and $\alpha=5$. Note the decreasing sizes of the jolts near $u=1$ when going from top to bottom. Following our interpretations already used above, we would expect to see performance improvement of the estimator $\mathrm{I}_n$ when going from $\varrho=0.5$, to $\varrho=0$, and then to $\varrho=-0.5$. This we shall indeed see in our following explorations.

\begin{figure}[h!]
\centering
\subfigure[$\alpha=1.5$, $\varrho=0.5$, $\mathrm{I}_F(h)\approx 0.7881$]{%
\includegraphics[height=0.35\textwidth]{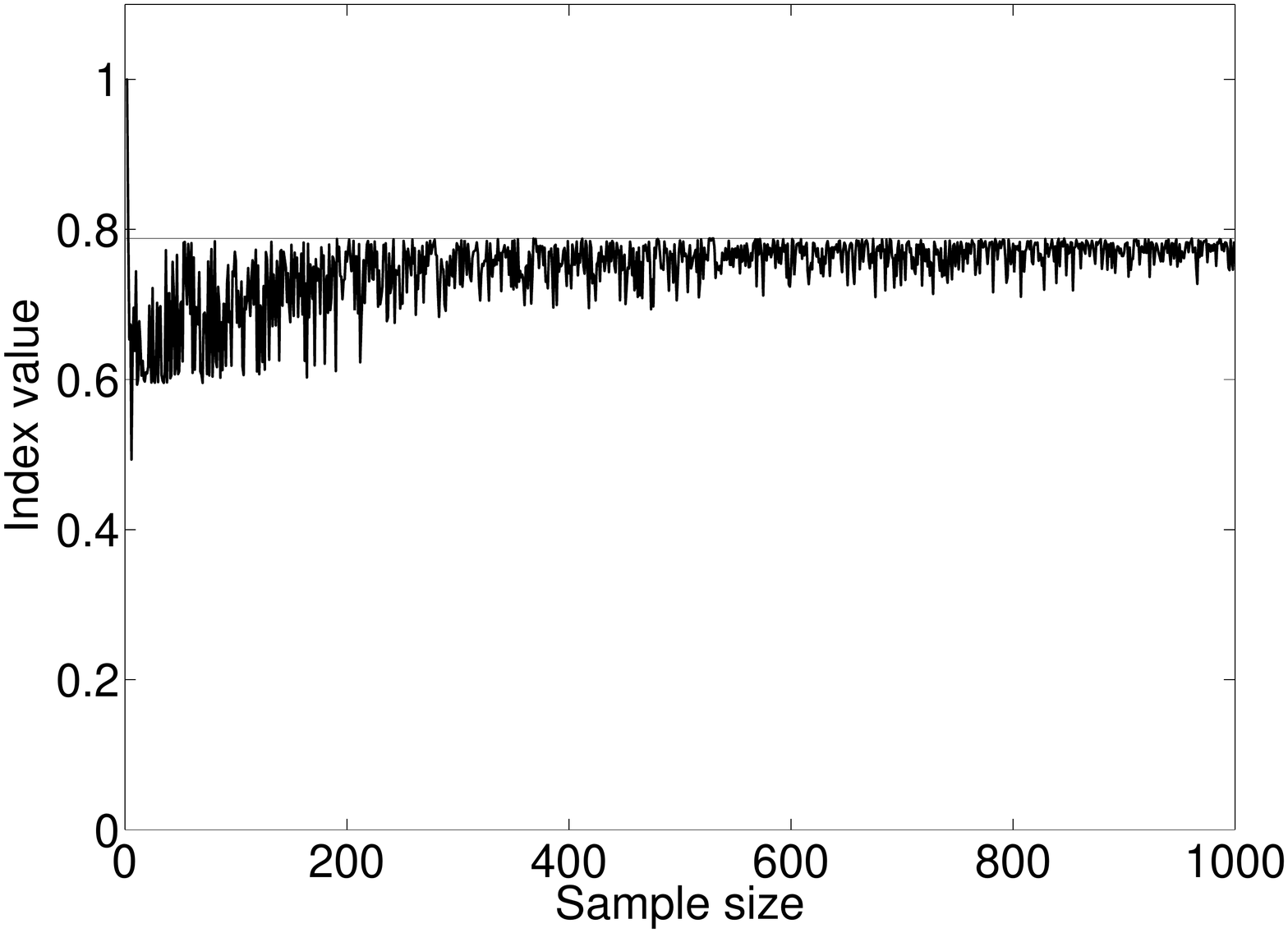}}
\hfill
\subfigure[$\alpha=5$, $\varrho=0.5$, $\mathrm{I}_F(h)\approx 0.7881$]{%
\includegraphics[height=0.35\textwidth]{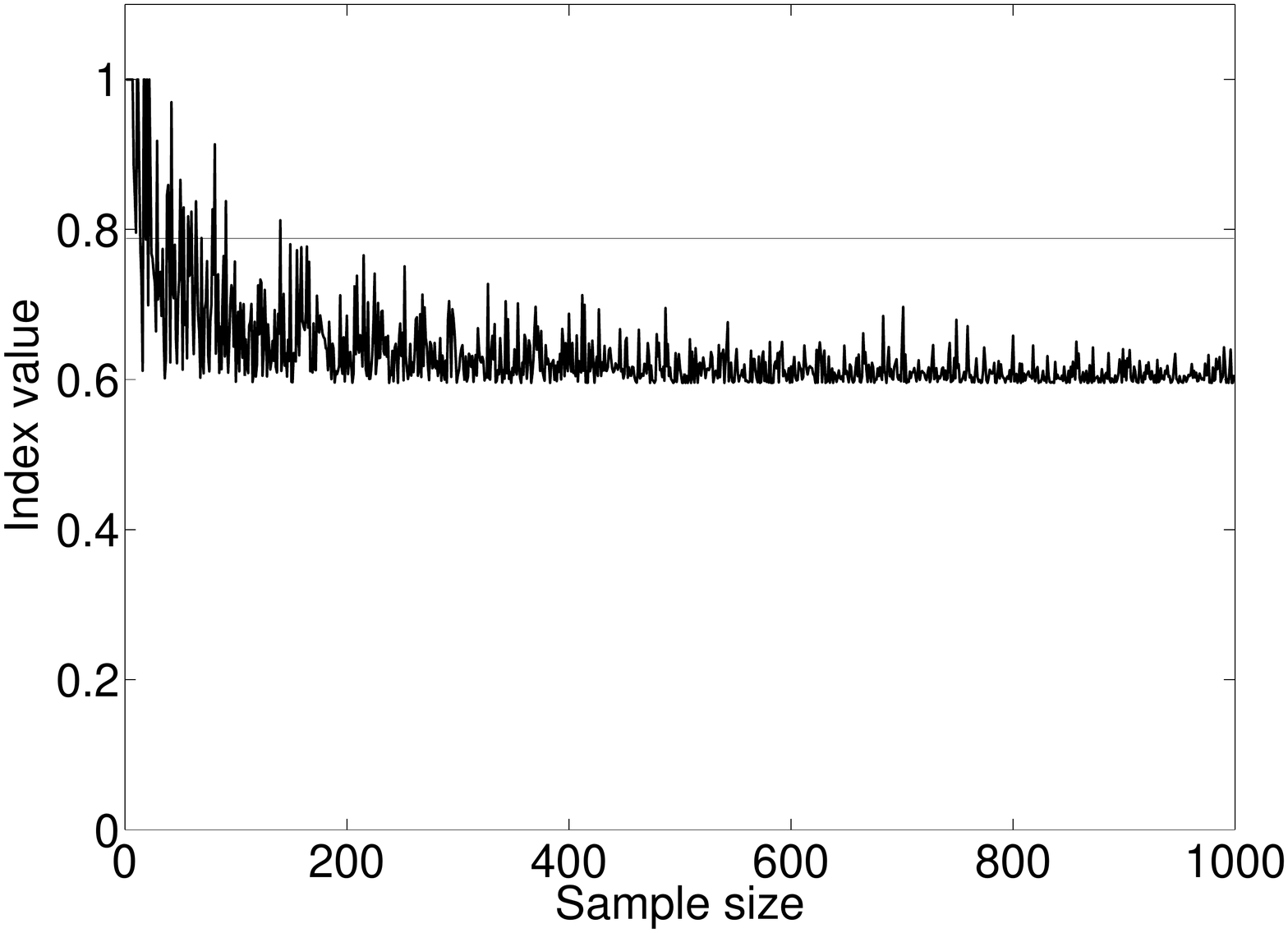}}
\\
\subfigure[$\alpha=1.5$, $\varrho=0$, $\mathrm{I}_F(h)\approx 0.7378$]{%
\includegraphics[height=0.35\textwidth]{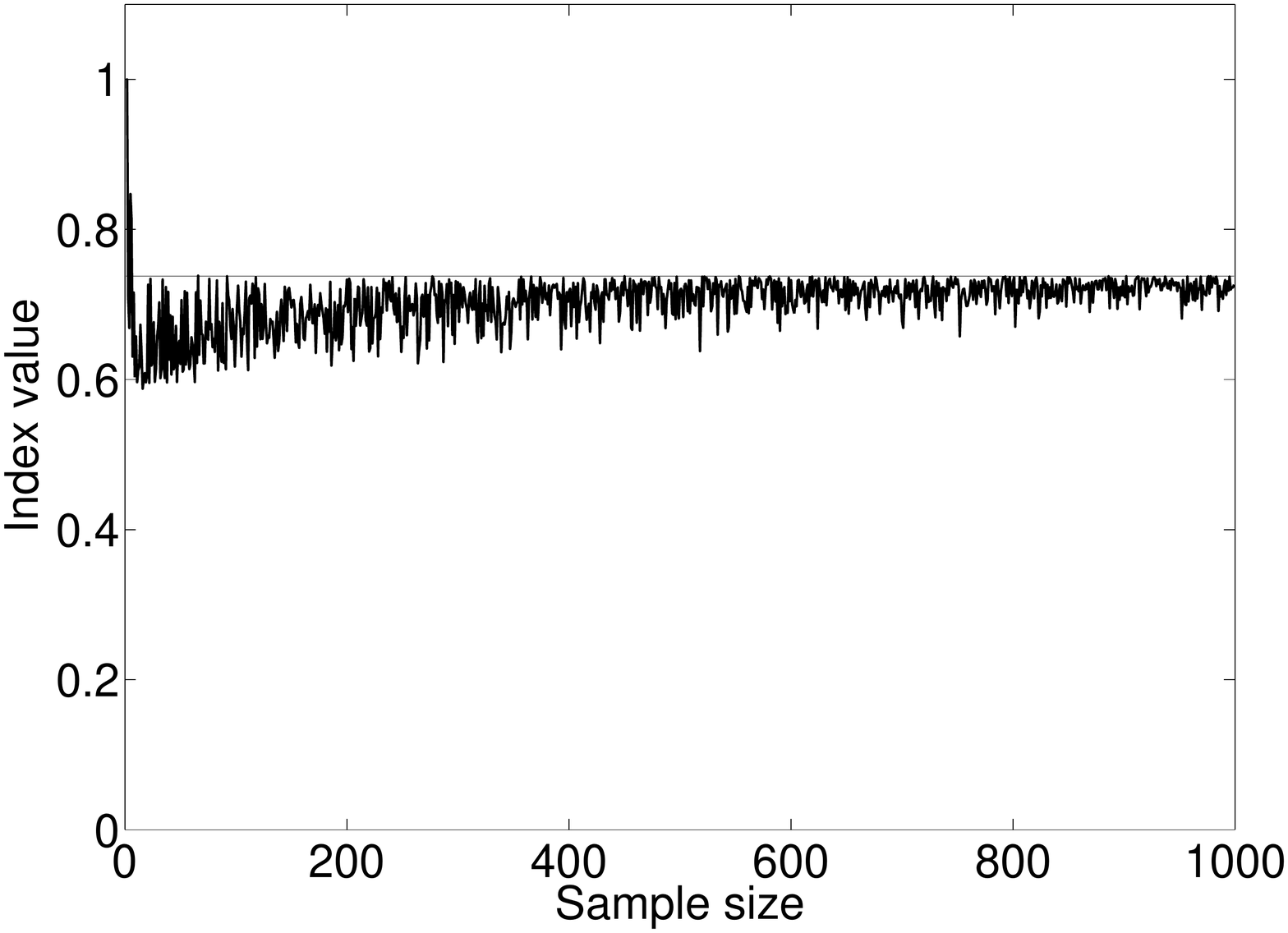}}
\hfill
\subfigure[$\alpha=5$, $\varrho=0$, $\mathrm{I}_F(h)\approx 0.7378$]{%
\includegraphics[height=0.35\textwidth]{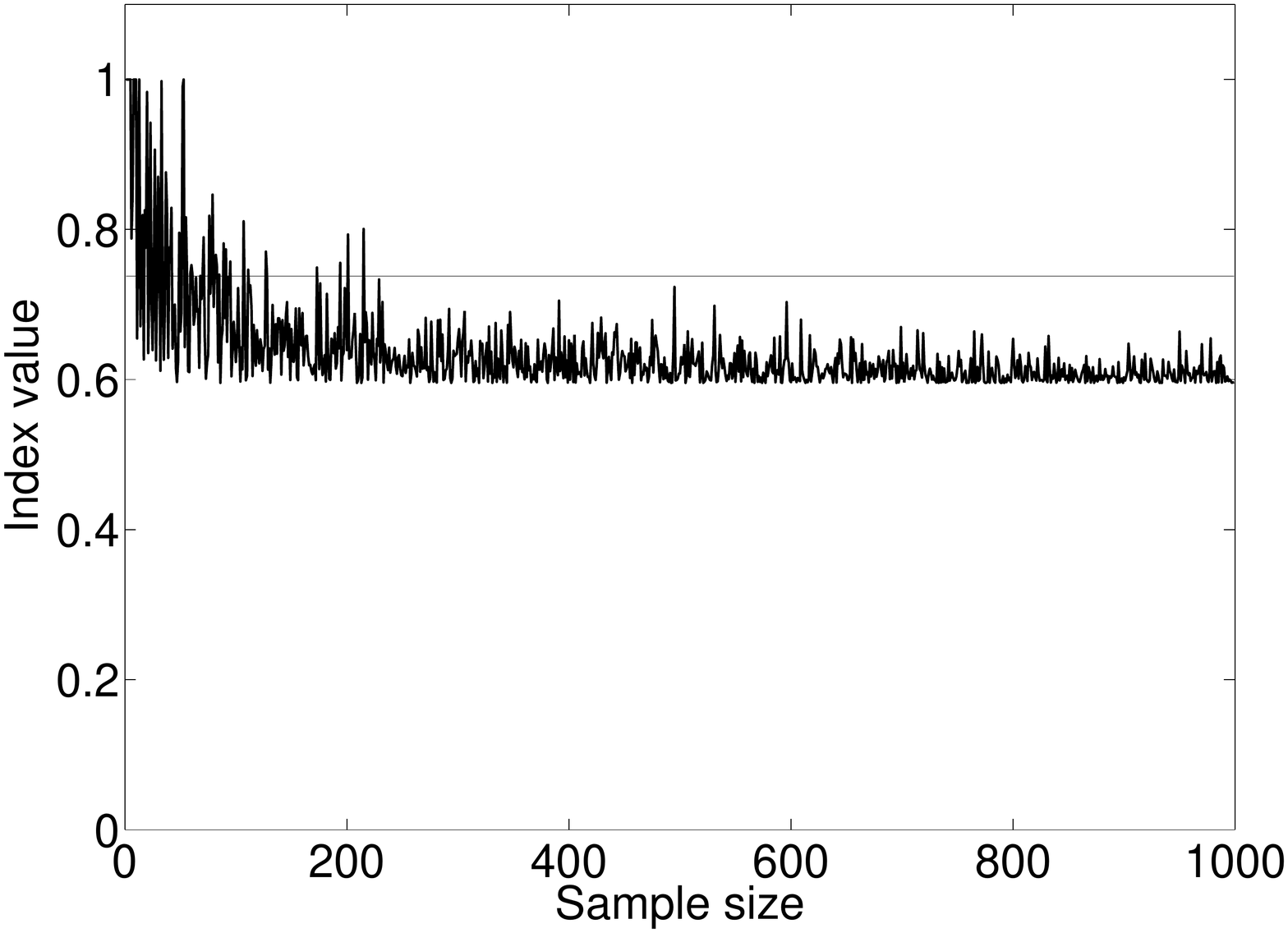}}
\\
\subfigure[$\alpha=1.5$, $\varrho=-0.5$, $\mathrm{I}_F(h)\approx 0.3459$]{%
\includegraphics[height=0.35\textwidth]{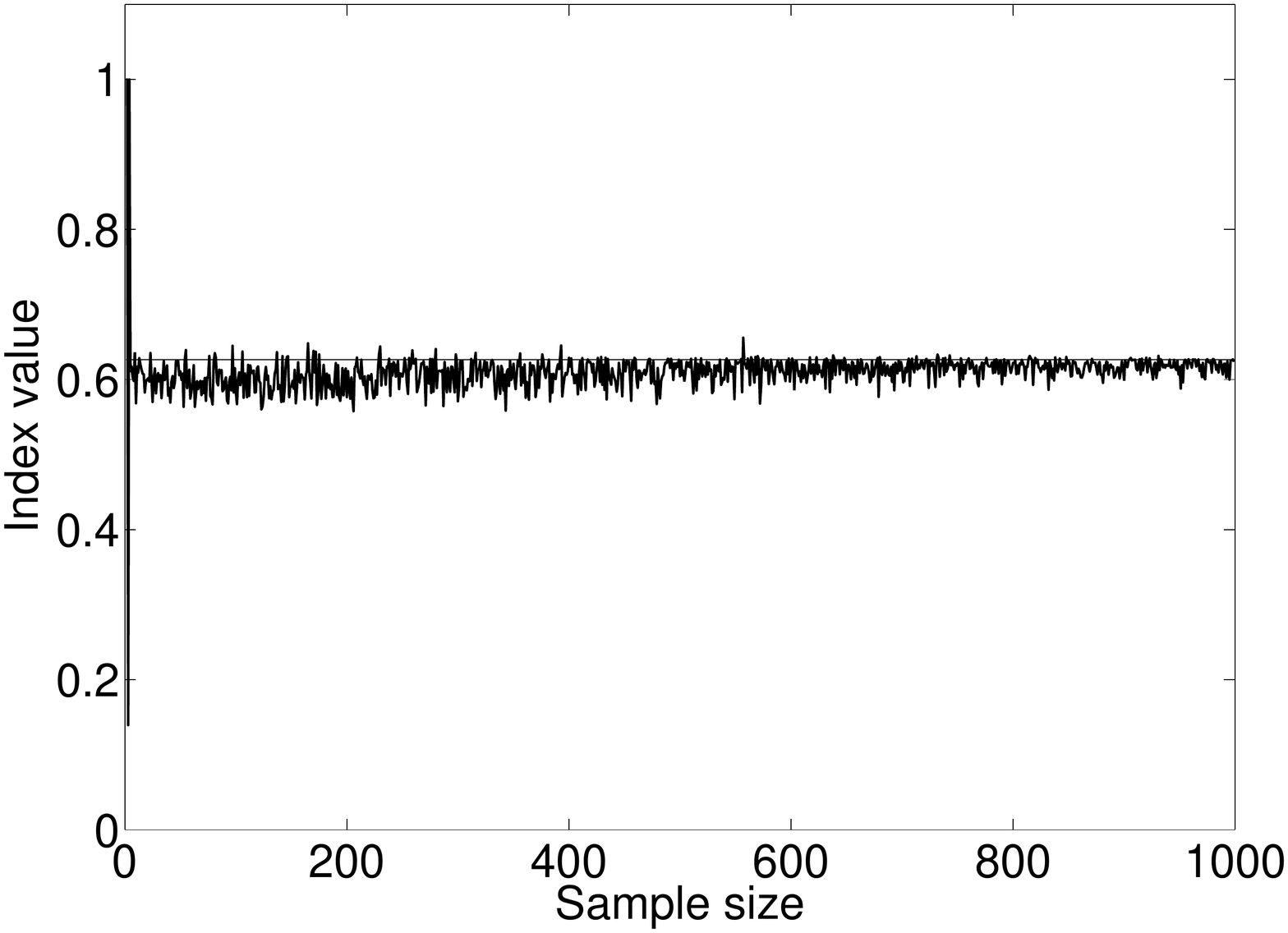}}
\hfill
\subfigure[$\alpha=5$, $\varrho=-0.5$, $\mathrm{I}_F(h)\approx 0.3459$]{%
\includegraphics[height=0.35\textwidth]{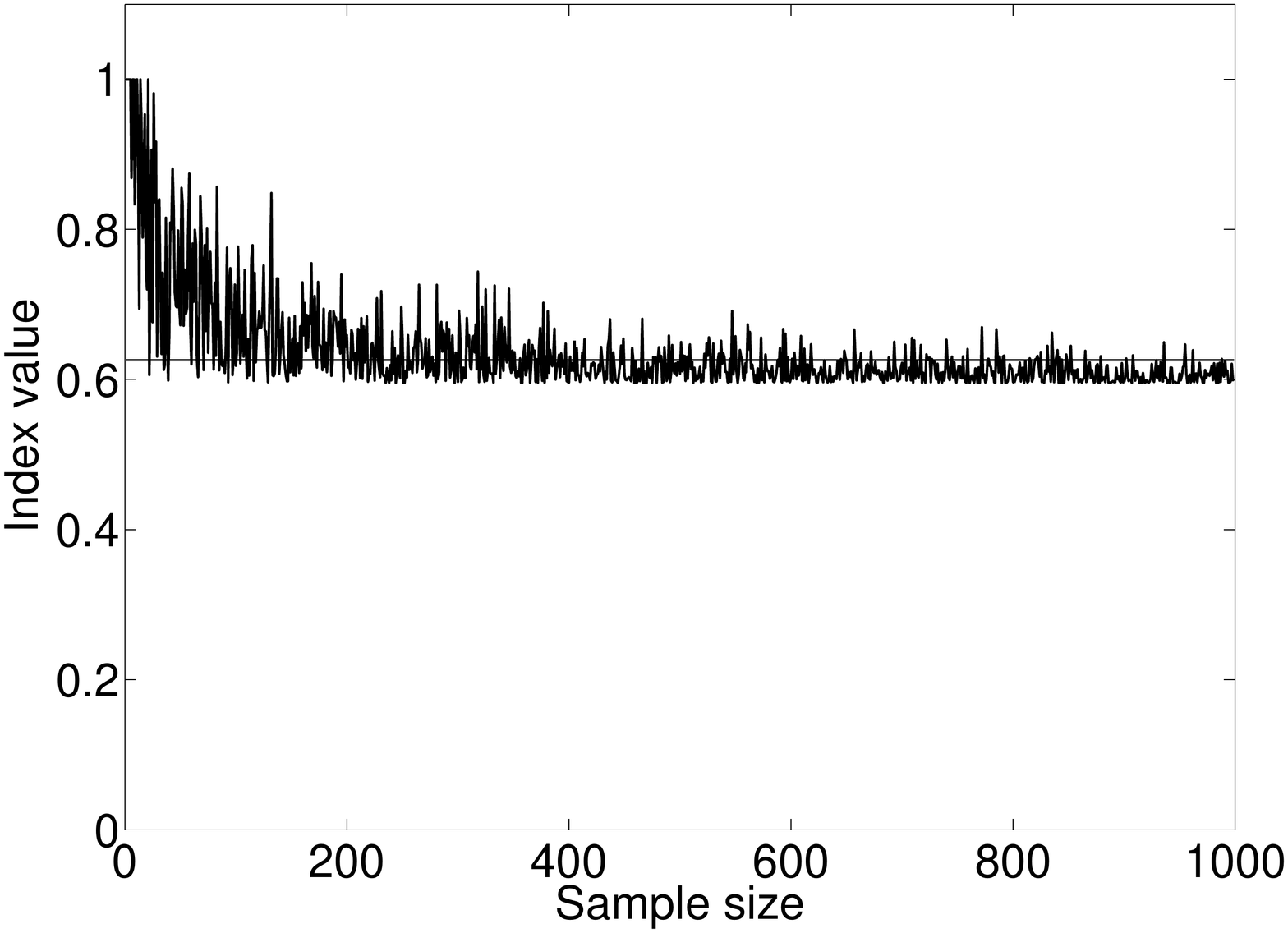}}
\caption{Performance of $\mathrm{I}_n$ in the window $(0,20]$, with the horizontal lines at the heights of the respective values of $\mathrm{I}_F(h)$; note the improving convergence from top to bottom.}
\label{fig-win20}
\end{figure}

We now look at Figure~\ref{fig-win20}
where we have depicted the performance of the estimator $\mathrm{I}_n$ and also plotted the horizontal lines at the heights of $\mathrm{I}_F(h)$. Consider first the case $\varrho=0$ when the transfer function $h(x)$ is continuous,  though not monotonic. Hence, the theoretical index $\mathrm{I}_F(h)$ is smaller than $1$, with its value calculated as follows:
\[
\mathrm{I}_F(h)
={\int_{0}^{20} (h')_{+}(x)\mathrm{d} x \over \int_{0}^{20}|h'|(x)\mathrm{d} x}
\approx {1.1178 \over 1.5151}\approx 0.7378.
\]
The two middle panels of Figure~\ref{fig-win20} correspond to this case. Note that when  $\alpha=1.5$, the index $\mathrm{I}_n$ converges, though very slowly, to the true value $\mathrm{I}_F(h)\approx 0.7378$, but there is no visible convergence in the case $\alpha=5$. To understand why this is so, we again look at Figure \ref{pdf-020}. We see from the left-hand panels of the figure that when $\alpha=1.5$, the random variable $X$ keeps producing considerable data until $x=10$ or so, whereas in the case $\alpha=5$, it virtually stops producing data already near $x=3$ or so. This scarcity of data in the case $\alpha=5$ transfers into the index's inability to capture a considerable portion of the increasing part of the transfer function $h(x)$ and thus, inevitably, tends to a value that is markedly below the theoretical one, which is
\[
\mathrm{I}_F(h)\approx 0.7378.
\]
The function  $H(u)$ on the right-hand panel of Figure \ref{pdf-020} corroborates this assessment. We next calculate the index $\mathrm{I}_F(h)$ when $\varrho=\pm 0.5$.

When $\varrho=0.5$, the transfer function $h(x)$ has a jump at $x=10$, whose size is $\Delta h_1=\varrho e^{9/10} 11/100\approx 0.1353$, a positive number. Hence,
\[
\mathrm{I}_F(h)\approx {0.1353+1.3427 \over 0.1353+1.7400 }\approx 0.7881.
\]
When $\varrho=-0.5$, the transfer function $h(x)$ has a drop at the point $x=10$ and thus $\Delta h_1=\varrho e^{9/10} 11/100\approx -0.1353$, a negative number. Consequently,
\[
\mathrm{I}_F(h)\approx {0.8928 \over 0.1353+1.2901}\approx 0.6264.
\]
In these two cases $\varrho=\pm 0.5$, just like when $\varrho=0$, we would conclude from the corresponding panels of Figure~\ref{fig-win20} that $\mathrm{I}_n$ slowly approaches the corresponding true values $\mathrm{I}_F(h)$ when $\alpha=1.5$, but hovers below the true value when $\alpha=5$. Our explanation of this phenomenon is the same as in the case $\varrho=0$ above: namely, it is due to the lack of inputs entering the filter on the right-hand side of the window $(0,20]$ and thus, in turn, results in very different values of the function $H(u)$ on different regions of its domain of definition $(0,1)$.

\section{Proofs}
\label{proofs}

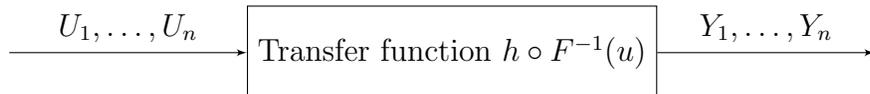
\begin{figure}[h!]\bigskip
\centering
\begin{tikzpicture}
\sbEntree{E1}
\sbBloc[8]{Bloc1}{Transfer function $h\circ F^{-1}(u)$}{E1}
\sbRelier[$U_{1},\dots , U_{n}$]{E1}{Bloc1}
\sbSortie[7]{S1}{Bloc1}
\sbRelier[$Y_{1}, \dots , Y_{n}$]{Bloc1}{S1}
\end{tikzpicture}
\caption{A counterpart of Figure~\ref{fig-2} with uniform inputs.}
	\label{fig-2b}
\end{figure}

The appearance of the composition $h\circ F^{-1}(u)$ of the transfer function and the quantile function in our main theorems is natural: it is due to the fact that $Y_i\stackrel{d}{=}h\circ F^{-1}(U_i)$, with $\stackrel{d}{=}$ denoting equality `in distribution,' where $U_i:=F(X_i)$, $i=1,\dots , n$, are iid uniform on $[0,1]$ random variables. The corresponding order statistics are $U_{1:n}< \cdots < U_{n:n}$ (we can assume without loss of generality that all the $U_i$'s are different) and thus $Y_{i,n}\stackrel{d}{=}h\circ F^{-1}(U_{i:n})$ for $i=1,\dots , n$. Consequently, Figure~\ref{fig-2} turns into Figure~\ref{fig-2b},
where the inputs are now uniformly distributed on $[0,1]$. This approach is very much in parallel to the pre-whitening technique in the time-series area that deals with  transfer functions and filtering (e.g., Box et al., 2015), except that in that area it is frequently more beneficial to reduce the input series to a white-noise time series.

\subsection{Proof of Theorem~\ref{Thm-2}}
\label{proof-2}

Since $X_{i:n}\stackrel{d}{=} F^{-1}(U_{i:n})$ for all $i=1,\dots, n$, statement \eqref{ii-main} follows from
\begin{equation}
\label{p_1}
\sum_{i=2}^n \big (h\circ F^{-1}(U_{i:n})-h\circ F^{-1}(U_{i-1:n})\big )_{+} \stackrel{\mathbf{P}}{\longrightarrow} \sum_{k=1}^m(\Delta h_k)_{+}+\int_{0}^{1} H_{+}(u)\mathrm{d}u
\end{equation}
and
\begin{equation}
\label{p_2}
\sum_{i=2}^n \big |h\circ F^{-1}(U_{i:n})-h\circ F^{-1}(U_{i-1:n})\big | \stackrel{\mathbf{P}}{\longrightarrow} \sum_{k=1}^m|\Delta h_k|+\int_{0}^{1}|H|(u)\mathrm{d}u.
\end{equation}
The proofs of these statements are virtually identical, and we thus prove only the first one.

For every $\tau_k$, we define the random variable
\[
N_k=\max_{1\leq i\leq n} \{ i: \,U_{i:n}\leq \tau_k \},
\]
which follows the binomial distribution with the parameters $\tau_k$ and $n$, because it is the number of those $ U_i$'s that do not exceed $\tau_k$. With $N_0:=0$ and $N_{m+1}:=n$, we decompose the sum on the left-hand side of statement \eqref{p_1} as follows:
\begin{equation}
\label{p_3}
\sum_{k=1}^m \big( h\circ F^{-1}(U_{N_k+1:n})-h\circ F^{-1}(U_{N_k:n})\big )_{+}
+ \sum_{k=0}^m \sum_{i=N_k+2}^{N_{k+1}} \big (h\circ F^{-1}(U_{i:n})-h\circ F^{-1}(U_{i-1:n})\big)_{+} .
\end{equation}
Since for each $k=1,\dots,m$, the variables $U_{N_k:n}$ and $U_{N_k+1:n}$ tend to $\tau_k$ almost surely from the left- and right-hand sides, respectively, the summands of the first sum of  \eqref{p_3} tend to $(\Delta h_k)_{+}$ almost surely. We next tackle the double sum of \eqref{p_3}.

By the mean-value theorem, we can find $\mathcal{U}_{i,n}\in (U_{i-1:n},U_{i:n})$ such that
\begin{align}
\label{p_4a}
\big(h\circ F^{-1}(U_{i:n})-h\circ F^{-1}(U_{i-1:n})\big)_{+}
&=H_{+}(\mathcal{U}_{i,n})(U_{i:n}-U_{i-1:n})
\notag
\\
&= H_{+}(p_{i,n})(U_{i:n}-U_{i-1:n}) + r_{i,n},
\end{align}
where the function $H(u)$ is given by equation (\ref{funct-g}), and the remainder term $r_{i,n}$ is
\[
r_{i,n}= \big( H_{+}(\mathcal{U}_{i,n})-  H_{+} (p_{i,n}) \big) (U_{i:n}-U_{i-1:n})
\]
with $p_{i,n}$ denoting the mean of the order statistic $U_{i:n}$, that is,
\[
p_{i,n}={i\over n+1}.
\]
We next show that
\begin{equation}
\label{p_3ii}
\sum_{k=0}^m \sum_{i=N_k+2}^{N_{k+1}} r_{i,n}\stackrel{\mathbf{P}}{\to} 0
 \quad \textrm{when} \quad n\to \infty .
\end{equation}
Since $m$ is fixed, in order to prove statement \eqref{p_3ii}, we need to show that, for every $k=0,1,\dots, m$, the inner sum of (\ref{p_3ii}) converges to $0$ in probability. To begin proving this fact, we first note that since $\mathcal{U}_{i,n}\in (U_{i-1:n},U_{i:n})$, we have $\mathcal{U}_{i,n}\in (\tau_{k},\tau_{k+1})$ for every $i\in [N_k+2, N_{k+1}]$. When $p_{i,n}\in (\tau_{k},\tau_{k+1})$ for $i\in [N_k+2, N_{k+1}]$, we can utilize the uniform continuity of $H_{+}(u)$ on the interval $(\tau_{k},\tau_{k+1})$, but $p_{i,n}$ may or may not be in the interval $ (\tau_{k},\tau_{k+1})$. This explains the necessity of somewhat more involved arguments that follow.

We start with the decomposition
\begin{equation}
\sum_{i=N_k+2}^{N_{k+1}} |r_{i,n}|
= \sum_{i\in \Delta_{k,1}} |r_{i,n}|
+
\sum_{i\in \Delta_{k,2}} |r_{i,n}|
+\sum_{i\in \Delta_{k,3}} |r_{i,n}|,
\label{p_3iii}
\end{equation}
where, with the notation $m_k^0=\max \{ i~:~ p_{i,n}\le \tau_{k} \}$, the three $\Delta $'s are
\begin{align*}
  \Delta_{k,1} & =\Big \{ i~:~ N_k+2\le i \le \max\{N_k+2,m_k^0+1\}-1 \Big \}, \\
  \Delta_{k,2} & =\Big \{ i~:~ \max\{N_k+2,m_k^0+1\}\le i \le \min\{N_{k+1}+1,m_{k+1}^0\}-1 \Big \}, \\
  \Delta_{k,3} & =\left \{ i~:~ \min\{N_{k+1}+1,m_{k+1}^0\}\le i \le N_{k+1} \right \}.
\end{align*}
The idea behind decomposition (\ref{p_3iii}) is the fact that, for every $i\in \Delta_{k,2}$, both $\mathcal{U}_{i,n}$ and $p_{i,n}$ belong to the open interval $ (\tau_{k},\tau_{k+1})$. The sets $\Delta_{k,1}$ and $\Delta_{k,3}$, whose definitions can of course be simplified at the cost of some transparency (e.g., $\Delta_{k,1}=\{ i~:~ N_k+2\le i \le m_k^0\}$), contain all the remaining indices $i\in [N_k+2, N_{k+1}]$. Some of the three sets (e.g., $\Delta_{k,1}$) can sometimes be empty.

We next show that the three sums on the right-hand side of equation (\ref{p_3iii}) converge to $0$ in probability when $n\to \infty $. First, we tackle the middle sum and start with the bound
\begin{multline}
\mathbf{P}\bigg( \sum_{i\in \Delta_{k,2}} |r_{i,n}|>\delta \bigg)
\le  \mathbf{P}\bigg(\max_{i\in \Delta_{k,2}} \big| H_{+}(\mathcal{U}_{i,n})- H_{+} (p_{i,n})\big|  >\delta  , \max_{2\leq i \leq n}  |\mathcal{U}_{i:n}-p_{i,n}|  \le \lambda \bigg )
\\
+ \mathbf{P}\bigg(\max_{2\leq i \leq n} |\mathcal{U}_{i:n}-p_{i,n}|  >\lambda \bigg )
\label{p_5i}
\end{multline}
for any $\delta >0$ and $\lambda>0$. As already noted above, the quantities $\mathcal{U}_{i,n}$ and $p_{i,n}$ are in the interval $ (\tau_{k},\tau_{k+1})$. Furthermore, it follows from conditions (C\ref{cond-1})--(C\ref{cond-3}) of Section~\ref{main} that the function $H_{+}(u)$ is uniformly continuous on every interval $(\tau_{k},\tau_{k+1})$, $k=0,1,\dots , m$, with finite right- and left-hand limits at the ends of these intervals.  Hence, with a sufficiently small $\lambda>0$ depending on $\delta >0$, we see that the first probability on the right-hand side of bound \eqref{p_5i} vanishes. To show that the right-most probability converges to $0$ when $n\to \infty $, we write the inequality
\begin{equation}
\label{p_6}
\mathbf{P}\Big(\max_{2\leq i \leq n}|\mathcal{U}_{i:n}-p_{i,n}| >\lambda \Big )
\leq 2 \sum_{i=1}^n \mathbf{P}\big ( |U_{i:n}-p_{i,n}| >\lambda^* \big )
\end{equation}
with $\lambda^*= \lambda -n^{-1}$ that can always be made larger than, say, $\lambda/2$ for all sufficiently large $n$. Using an exponential bound for the uniform order statistics (e.g., Shorack and Wellner, 1986, proof of Corollary 2, pp.~456--457), we have for every $i=1,\dots,n$ and all $s >1/\sqrt{n}$,
\begin{align*}
\mathbf{P}\big( |U_{i:n}-p_{i,n}|  >s \big)
&=\mathbf{P}\big( \sqrt{n}|U_{i:n}-p_{i,n}|  >\sqrt{n}s \big)
\\
&\leq 2\exp\{-\sqrt{n}s/10\}.
\end{align*}
Hence, the right-hand side of bound \eqref{p_6} and thus, in turn, the middle sum on the right-hand side of bound (\ref{p_3iii}) converge to $0$ when $n\to \infty $.

Now we are left to show that the first and third sums on the right-hand side of equation (\ref{p_3iii}) converges to $0$ when $n\to \infty $, but since the proofs are very similar, we tackle only the first sum. We need an auxiliary result. Namely, let
\[
M_n= \max_{1\leq i \leq n+1} \big \{U_{i:n}-U_{i-1:n} \big \},
\]
which is the maximal spacing, where $U_{0:n}:=0$ and $U_{n+1:n}:=1$. Then the limiting distribution of $nM_n-\log n$ is the standard Gumbel distribution
(e.g., del Barrio et al., 2007, p.~140), that is,
\begin{equation}
\label{p3_10s}
\lim_{n\to\infty}\mathbf{P}\big(nM_n-\log n\leq x\big)=\exp\{-e^{-x}\}  \quad \textrm{for} \quad x \in \mathbf{R}.
\end{equation}
With $\lambda >1$ denoting any constant, we set $x=(\lambda-1)\log n$ in statement (\ref{p3_10s}) and have
\begin{equation}
\label{p3_10s2}
\mathbf{P}\bigg(M_n\ge {\lambda \log n \over n} \bigg)=o(1)  \quad \textrm{when} \quad n\to \infty .
\end{equation}
Consequently, for every $\delta >0$ and $\lambda>1$,
\begin{align}
\label{p_5ij}
\mathbf{P}\bigg( \sum_{i\in \Delta_{k,1}}  |r_{i,n}|>\delta \bigg)
&\le  \mathbf{P}\bigg(\sum_{i\in \Delta_{k,1}} \big| H_{+}(\mathcal{U}_{i,n})- H_{+} (p_{i,n})\big|(U_{i:n}-U_{i-1:n})  >\delta  , M_n \le {\lambda \log n \over n}  \bigg )
+ o(1)
\notag
\\
&\le  \mathbf{P}\bigg( \#\{\Delta_{k,1}\}  \Vert H_{+}\Vert {\lambda \log n \over n}  >\delta  , M_n \le {\lambda \log n \over n}  \bigg )
+ o(1)
\notag
\\
&\le  \mathbf{P}\bigg( \#\{\Delta_{k,1}\}{\lambda \log n \over n} >\delta^*  \bigg )
+ o(1)
\end{align}
when $n\to \infty $, where $\delta ^*= \delta / \Vert H_{+}\Vert $ and
\[
\Vert H_{+}\Vert =\sup_{0<u<1}H_{+}(u)<\infty .
\]
To show that the probability on the right-hand side of bound (\ref{p_5ij}) converges to $0$, we first note that the cardinality of the set $\Delta_{k,1}=\{ i~:~ N_k+2\le i \le m_k^0\}$ does not exceed $N_k+1-m_k^0$. Since $N_k$ follows the binomial distribution with the parameters $\tau_k$ and $n$, its variance is of the order $O(n)$. Consequently, Chebyshev's inequality implies
\[
\mathbf{P}\bigg( \#\{\Delta_{k,1}\}{\lambda \log n \over n} >\delta^*  \bigg )
=O\bigg ( {(\log n)^2 \over n}\bigg ),
\]
and so the first sum on the right-hand side of bound (\ref{p_3iii}) converge to $0$ when $n\to \infty $. We analogously arrive at the same conclusion with respect to the third sum on the right-hand side of bound (\ref{p_3iii}), thus concluding the proof of statement (\ref{p_3ii}).

Equipped with the above asymptotic statements and after a slight rearrangement of terms, we  see that the sum on the left-hand side of statement (\ref{p_1}) can be written as follows:
\begin{align}
\sum_{i=2}^n \Big (h\circ F^{-1}(U_{i:n})-h\circ F^{-1}(U_{i-1:n})\Big )_{+}
= \sum_{k=1}^m(\Delta h_k)_{+}
&+ \sum_{i=2}^n  H_{+}(p_{i,n})\Big (U_{i:n}-U_{i-1:n}-\frac{1}{n+1}\Big )
\notag
\\
&-\sum_{k=1}^{m-1} H_{+}\Big(\frac{N_k+1}{n+1}\Big)\big(U_{N_k+1:n}-U_{N_k:n}\big)
\notag
\\
&+ \frac{1}{n+1}\sum_{i=2}^n  H_{+}(p_{i,n})
+o_{\mathbf{P}}(1).
\label{eq-03}
\end{align}
Since the right-most sum of equation (\ref{eq-03}) converges to $\int_{0}^{1} H_{+}(u)\mathrm{d}u$ when $n\to\infty$, statement (\ref{p_1}) and thus Theorem~\ref{Thm-2} follow if the two middle sums on the right-hand side of equation (\ref{eq-03}) are of the order $o_{\mathbf{P}}(1)$ when $n\to\infty$. To prove this, we start with the bound
\begin{equation}
\label{p_8}
\sum_{k=1}^{m-1} H_{+}\bigg(\frac{N_k+1}{n+1}\bigg )\big(U_{N_k+1:n}-U_{N_k:n}\big)
\leq \Vert H_{+}\Vert \sum_{k=1}^{m-1} \big(U_{N_k+1:n}-U_{N_k:n}\big)
\end{equation}
whose right-hand side converges to $0$ because both $U_{N_k+1:n}$ and $U_{N_k:n}$ converge to $\tau_k$ when $n\to \infty $. To prove that the last sum on the right-hand side of equation (\ref{eq-03}) converges to $0$ in probability, we show that its second moment converges to $0$. With $S_i:=U_{i:n}-U_{i-1:n}$ denoting the $i^{\textrm{th}}$ uniform spacing, we employ the following three formulas (e.g., Shorack and Wellner, 1986, p.~721):
\begin{gather*}
\mathbf{E}(S_i)=\frac 1{n+1},
\\
\mathbf{Var}(S_i)=\frac n{(n+1)^2(n+2)},
\\
\mathbf{Cov}(S_i,S_j)={ -1 \over (n+1)^2(n+2)}  \quad \textrm{when} \quad i\ne j .
\end{gather*}
Consequently,
\begin{align}
\label{p_11}
\mathbf{E}\bigg[\bigg \{ \sum_{i=2}^n  H_{+}(p_{i,n})\bigg (U_{i:n}-U_{i-1:n}-\frac{1}{n+1}\bigg )\bigg \}^2\bigg ]
&=\frac{n}{(n+1)^2(n+2)}\sum_{i=2}^{n} H_{+}^2(p_{i,n})
\notag
\\
&\quad -\frac{2}{(n+1)^2(n+2)}\sum_{2\leq i<j\leq n}H_{+}(p_{i,n})H_{+}(p_{j,n})
\notag
\\
&\leq \frac 1n \Vert H_{+}\Vert^2.
\end{align}
This finishes the proof of statement (\ref{p_1}) and thus, in turn, the proof of Theorem~\ref{Thm-2}.  $\hfill \square $

\subsection{Proof of Theorem~\ref{Thm-3}}
\label{proof-3}

With $\theta >0$ defined in the formulation of Theorem~\ref{Thm-3} and when $m\neq 0$, we assume without loss of generality that $\theta <\tau_1$ and $\tau_m < 1-\theta$. We need to prove statements \eqref{p_1} and \eqref{p_2} under the conditions of Theorem~\ref{Thm-3}, but we shall prove only statement \eqref{p_1} because the other one can be established analogously.

We split the sum on the left-hand side of statement \eqref{p_1} into three parts:
\begin{align}
\label{p3_1}
\bigg(\sum_{i=2}^{\lfloor n\theta \rfloor } + \sum_{\lfloor n\theta \rfloor +1}^{n-\lfloor n\theta \rfloor} +\sum_{n-\lfloor n\theta \rfloor +1}^{n}\bigg) \big (h\circ F^{-1}(U_{i:n})-h\circ F^{-1}(U_{i-1:n})\big )_{+}.
\end{align}
The conditions of Theorem~\ref{Thm-3} imply that the function $H(u)$ and thus $H_{+}(u)$ are   uniformly continuous on the interval $[\theta,1-\theta]$. Therefore, following the same arguments used in the proof of Theorem~\ref{Thm-2}, we conclude that the middle sum in \eqref{p3_1} converges to $\sum_{k=1}^m(\Delta h_k)_{+}+\int_{\theta}^{1-\theta} H_{+}(u)\mathrm{d}u$ when $n\to \infty $. It remains to show that the first and third sums in \eqref{p3_1} converge to $\int_{0}^{\theta} H_{+}(u)\mathrm{d}u$ and $\int_{1-\theta}^{1} H_{+}(u)\mathrm{d}u$, respectively, but we shall prove this for the first sum only because the third sum can be treated analogously.

Similarly to the proof of Theorem~\ref{Thm-2}, using the notation
$p_{i,n}=i/(n+1)$ and the mean-value theorem, we have
\begin{equation}
\label{p3_2}
\sum_{i=2}^{\lfloor n\theta \rfloor }\big (h\circ F^{-1}(U_{i:n})-h\circ F^{-1}(U_{i-1:n})\big )_{+}
=
\frac 1{n+1}\sum_{i=2}^{\lfloor n\theta \rfloor } H_{+}(p_{i,n})  + r_{1,n}+ r_{2,n}
\end{equation}
with the remainder terms
\begin{equation}
\label{p3_3}
r_{1,n}
= \sum_{i=2}^{\lfloor n\theta \rfloor } H_{+}(p_{i,n}) \bigg ( U_{i:n}-U_{i-1:n}-\frac 1{n+1} \bigg )
\end{equation}
and
\begin{equation}
\label{p3_4}
r_{2,n}
=\sum_{i=2}^{\lfloor n\theta \rfloor }\big( H_{+}(\mathcal{U}_{i,n})-H_{+}(p_{i,n})\big )
(U_{i:n}-U_{i-1:n}).
\end{equation}
Since the main term on the right-hand side of equation (\ref{p3_2}) converges to $\int_{0}^{\theta} H_{+}(u)\mathrm{d}u$ when $n\to \infty $, we are left to show that both $r_{1,n}$ and $r_{2,n}$ converge to zero in probability.

We start with the remainder term $r_{1,n}$ and prove that its second moment converges to $0$. For this we employ the moment-type formulas for the spacings  $S_i=U_{i:n}-U_{i-1:n}$ given at the end of the proof of the previous theorem, and arrive at the bound
\begin{align*}
\mathbf{E}\bigg[\bigg\{\sum_{i=2}^{\lfloor n\theta \rfloor } H_{+}(p_{i,n})  \bigg(U_{i:n}-U_{i-1:n}-\frac 1{n+1} \bigg)\bigg\}^2\bigg ]
=&\frac{n}{(n+1)^2(n+2)}\sum_{i=2}^{\lfloor n\theta \rfloor } H_{+}^2(p_{i,n})
\\
&-\frac{2}{(n+1)^2(n+2)}\sum_{2\leq i<j\leq \lfloor n\theta \rfloor }H_{+}(p_{i,n})H_{+}(p_{j,n})
\\
\leq &\frac{1}{n^2}\sum_{i=2}^{\lfloor n\theta \rfloor } H_{+}^2(p_{i,n}).
\end{align*}
With the help of assumption \eqref{iib}, we continue the above bound and have
\begin{align*}
\frac{1}{n^2}\sum_{i=2}^{\lfloor n\theta \rfloor } H_{+}^2(p_{i,n})
\leq  \frac{c}{n^2}\sum_{i=2}^{\lfloor n\theta \rfloor } p_{i,n}^{-2+2b}
&\leq \frac{c}{n^2}\sum_{i=2}^{\lfloor n\theta \rfloor } \Big(\frac in\Big)^{-2+2b}
\\
& \leq \frac{c}{n^{2b}} \int_{1}^{n}{1\over x^{2-2b}}\mathrm{d}x
\leq c{\log n \over \min \{n^{2b},n\} }
\end{align*}
with right-hand side converging to $0$ when $n\to \infty$. Hence, $r_{1,n}$ converges to $0$ when $n\to \infty$.

We now tackle $r_{2,n}$, whose definition is given by equation \eqref{p3_4}. Using the bound $|H_{+}(u)-H_{+}(v)|\leq |H(u)-H(v)|$ and the notation $I_{i,n}$ for the interval
\[
I_{i,n}=\big [\min\{U_{i-1:n}, p_{i,n}\},~\max\{U_{i:n}, p_{i,n}\} \big ],
\]
we have
\begin{align}
\label{p3_6}
|r_{2,n}| &\leq \sum_{i=2}^{\lfloor n\theta \rfloor }| H(\mathcal{U}_{i,n})-H(p_{i,n})|
(U_{i:n}-U_{i-1:n})
\notag
\\
&\leq \sum_{i=2}^{\lfloor n\theta \rfloor }\sup_{u\in I_{i,n}} \{ |H'(u)|\}
\max\big\{|U_{i:n}-p_{i,n}|,|U_{i-1:n}-p_{i,n}|\big\}(U_{i:n}-U_{i-1:n})
\notag
\\
&\leq \sum_{i=2}^{\lfloor n\theta \rfloor }\max_{u\in I_{i,n}} \big \{ u(1-u) \big \}^{-2+b}|U_{i:n}-p_{i,n}|(U_{i:n}-U_{i-1:n})
\notag
\\
&\quad + \sum_{i=2}^{\lfloor n\theta \rfloor }\max_{u\in I_{i,n}} \big \{ u(1-u) \big \}^{-2+b}|U_{i-1:n}-p_{i,n}|(U_{i:n}-U_{i-1:n}).
\end{align}
Without loss of generality we assume $b\leq 1$, since otherwise (i.e., when $b> 1$) the function $H(u)$ is bounded near the endpoints $0$ and $1$, and this takes us back to the case already treated in Theorem~\ref{Thm-2}. Hence, when $b\leq 1$, the function
\[
D(u)=\big ( u(1-u) \big )^{-2+b}
\]
is convex, and we have the bound $\max_{u\in I_{i,n}}D(u)\leq D(p_{i,n})+D(U_{i:n})+D(U_{i-1:n})$. Furthermore,  $|U_{i-1:n}-p_{i,n}|$ does not exceed $ |U_{i-1:n}-\frac{i-1}{n+1}|+\frac 1{n+1}$, and so estimation of the second sum on the right-hand side of equation \eqref{p3_6} is analogous to that of the first sum. Thus, in order to prove that $r_{2,n}$ tends to zero in probability when $n\to \infty $, it suffices to show that the following three sums converge to zero in probability:
\begin{align}
\label{p3_7}
&r_{2,n}^{(1)}:=\sum_{i=2}^{\lfloor n\theta \rfloor } D(p_{i,n})|U_{i:n}-p_{i,n}|(U_{i:n}-U_{i-1:n}),
\notag
\\
&r_{2,n}^{(2)}:=\sum_{i=2}^{\lfloor n\theta \rfloor } D(U_{i:n})|U_{i:n}-p_{i,n}|(U_{i:n}-U_{i-1:n}),
\notag
\\
&r_{2,n}^{(3)}:=\sum_{i=2}^{\lfloor n\theta \rfloor } D(U_{i-1:n})|U_{i:n}-p_{i,n}|(U_{i:n}-U_{i-1:n}).
\end{align}

First we consider $r_{2,n}^{(1)}$ and show that its first moment converges to $0$, which in turn implies its convergence in probability. For this, we write the bound
\begin{equation}
\label{p3_8}
\mathbf{E}\big[r_{2,n}^{(1)}\big]\leq c \sum_{i=2}^{\lfloor n\theta \rfloor } \Big(\frac i{n+1}\Big)^{-2+b}\mathbf{E}\big[|U_{i:n}-p_{i,n}|(U_{i:n}-U_{i-1:n})\big ]
\end{equation}
and then apply the Cauchy-Bunyakovsky-Schwarz inequality together with the bounds
\[
\mathbf{Var}[U_{i:n}]=\frac i{n+1} \bigg(1-\frac i{n+1}\bigg )\frac{1}{n+2}\le \frac i{n^2}
\]
and
\[
\mathbf{E}\big[(U_{i:n}-U_{i-1:n})^2\big ]=\frac 2{(n+1)(n+2)}\le \frac 2{n^2}.
\]
We have
\[
\mathbf{E}\big[r_{2,n}^{(1)}\big]
\le  \frac{c}{n^{3/2}}\sum_{i=2}^{\lfloor n\theta \rfloor }
 \Big(\frac in\Big)^{-3/2+b}
\le \frac{c}{n^b} \int_{1}^{n} {1\over x^{3/2-b}}\mathrm{d}x
\le c {\log n \over \min\{n^{b},n^{1/2}\} }
\]
with the right-hand side converging to $0$  when $n\to \infty $ because $b>0$. This proves that the remainder term $r_{2,n}^{(1)}$ converges to $0$ in probability.

Now we consider the remainder term $r_{2,n}^{(2)}$.  We again use $M_n$ to denote the maximal spacing and apply statement (\ref{p3_10s}). For every fixed $\delta>0$, we have
\begin{equation}
\label{p3_10}
\mathbf{P}\big(r_{2,n}^{(2)} \geq \delta\big) \leq  \mathbf{P}\bigg( \frac{\log n +A}{n}\sum_{i=2}^{\lfloor n\theta \rfloor } D(U_{i:n})|U_{i:n}-p_{i,n}|\geq \delta\bigg) +\mathbf{P}\big( nM_n -\log n\geq A \big).
\end{equation}
By statement (\ref{p3_10s}) and setting a sufficiently large $A$, the right-most probability can be made as small as desired for all sufficiently large $n$. This implies that in order to prove convergence of $\mathbf{P}\big(r_{2,n}^{(2)} \geq \delta\big)$ to $0$ when $n\to \infty $, we need to show that, for every fixed $A>0$, the first probability on the right-hand side of equation (\ref{p3_10}) converges to $0$ when $n\to \infty $. To this end, it is sufficient to prove the convergence when the function $D(u)$ is replaced by $D_0(u)=u^{-2+b}$. For this, we first estimate $U_{i:n}$ from above by $U_{[n\theta]:n}$ and then estimate the latter one by $G_n^{-1}(\theta )$, which is the inverse (i.e., quantile) function at the point $\theta $ arising from the empirical cdf based on the independent and uniformly on $[0,1]$ distributed random variables $U_1,\dots,U_n$. We have the bound
\begin{equation}
\label{p3_10u}
D(U_{i:n})\le {1\over U_{i:n}^{2-b} (1-G_n^{-1}(\theta ))^{2-b}}.
\end{equation}
We now split the first probability on the right-hand side of equation (\ref{p3_10}) into two parts: when  $1-G_n^{-1}(\theta )\ge 1-\gamma $ and when $1-G_n^{-1}(\theta )< 1-\gamma $. The probability of the latter event is equal to the probability of $\theta >G_n(\gamma )$, which, by the central limit theorem converges to $0$ whenever $\theta < \gamma$, and this is how we choose $\theta $. (Recall also the note at the beginning of this proof that $\theta$ must be smaller than $\tau_1$ and $1-\tau_m$.) Hence, we have reduced the problem to showing that, for every $\delta >0$, the probability
\begin{equation}
\label{p3_11}
\mathbf{P}\bigg( \frac{c\log n}{n}\sum_{i=2}^{\lfloor n\theta \rfloor } {1\over U_{i:n}^{2-b}}|U_{i:n}-p_{i,n}|\geq \delta\bigg)
\end{equation}
converges to $0$ when $n\to \infty $. To estimate this probability, we first apply bound (\ref{p3_10s}) and, with the notation $q_{i,n} =1-p_{i,n}$, obtain (e.g., Shorack and Wellner, 1986, pp.~453--455)
\[
\mathbf{P}\Big(\sqrt{n}|U_{i:n}-p_{i,n}|\geq \sqrt{p_{i,n}q_{i,n}}\,\lambda\Big)\le 2\exp\{-\lambda/10\}
\]
for all $i=1,\dots,n$ and all $\lambda \geq 1$.  This inequality implies that, for every $0<r<1/2$, which we shall specify later, the inequality
\[
\mathbf{P}\Big(|U_{i:n}-p_{i,n}|\geq n^{-1/2+r}\sqrt{p_{i,n}q_{i,n}}\lambda\Big)\le 2\exp\{-\lambda n^r/10\}
\]
holds and, in turn, implies
\begin{equation*}
\label{p3_12}
\mathbf{P}\bigg (\bigcup_{2\leq i\leq \lfloor n\theta \rfloor } \Big \{ |U_{i:n}-p_{i,n}|\geq n^{-1/2+r}\sqrt{p_{i,n}q_{i,n}}\lambda \Big \} \bigg ) \leq c n \exp\{-\lambda n^r/10\},
\end{equation*}
which converges to $0$ when $n\to\infty$.  Hence, probability \eqref{p3_11} is
\begin{equation}
\label{p3_13}
\mathbf{P}\bigg( \frac{c\log n}{n^{3/2-r}}\sum_{i=2}^{\lfloor n\theta \rfloor } {\sqrt{p_{i,n}}\over U_{i:n}^{2-b}}\geq \delta\bigg)+o(1)
\end{equation}
with a constant $c$ that does not depend on $n$. We need to prove that,  for every $\delta>0$, probability (\ref{p3_13}) converges to $0$ when $n\to \infty $. To this end, we first write
\[
 {1\over U_{i:n}}={i/n\over U_{i:n}}{n\over i}={G_n(U_{i:n})\over U_{i:n}}{n\over i}
\le {n\over i}\sup_{t\in [U_{1:n},1]}{G_n(t)\over t}=  {n\over i} \Gamma_n(U_{1:n})
\]
with the notation
\[
\Gamma_n(a):= \sup_{t\in [a,1]}{G_n(t)\over t}.
\]
This gives us the following bounds
\begin{align}
\mathbf{P}\bigg( \frac{c\log n}{n^{3/2-r}}\sum_{i=2}^{\lfloor n\theta \rfloor } {1\over U_{i:n}^{2-b}}(p_{i,n})^{1/2}\geq \delta\bigg)
&\le
\mathbf{P}\bigg( \Gamma_n^{2-b}(U_{1:n})
\frac{c\log n}{n^{3/2-r}}\sum_{i=1}^{n} \Big({n\over i}\Big)^{3/2-b}\geq \delta\bigg )
\notag
\\
&\le
\mathbf{P}\bigg( \Gamma_n^{2-b}(U_{1:n})
\frac{c\log n}{n^{b-r}}\int_{1}^{n} {1\over x^{3/2-b}}\mathrm{d}x\geq \delta\bigg )
\notag
\\
&\le \mathbf{P}\bigg( \Gamma_n^{2-b}(U_{1:n})\frac{c(\log n)^2}{n^{\min\{b,1/2\}-r}}\geq \delta\bigg ).
\label{p3_13b}
\end{align}
Of course, throughout the calculations, the value of the constant $c$ might have changed from line to line, but it never depends on $n$. Now, since we always have $r<1/2$ and can choose $r$ so that $r<b$, the probability on the right-hand side of bound (\ref{p3_13b}) converges to $0$ if, for any (sufficiently small) $\alpha>0$,
\begin{equation}
\label{p3_13c}
\mathbf{P}\big( \Gamma_n(U_{1:n})\geq n^{\alpha }\big )\to 0   \quad \textrm{when} \quad n\to \infty .
\end{equation}
We write
\begin{equation}
\label{p3_13d}
\mathbf{P}\big( \Gamma_n(U_{1:n})\geq n^{\alpha }\big )
\le \mathbf{P}\big( \Gamma_n(a_n)\geq n^{\alpha }\big )
+ \mathbf{P}\big( U_{1:n}\le  a_n\big ) ,
\end{equation}
where $a_n$ can be any sequence of positive real numbers, but we shall soon choose it in a special way. Bound (5) on page 415 of Shorack and Wellner (1986) says that, for every $a\in [0,1]$ and $x\ge 1$,
\[
\mathbf{P}\big( \Gamma_n(a)\geq x \big )
\le \exp \big \{ -n a \big(x(\log x-1)+1\big) \big \} .
\]
Applying this result on the right-hand side of bound (\ref{p3_13d}), we obtain
\begin{equation}
\label{p3_13m}
\mathbf{P}\big( \Gamma_n(U_{1:n})\geq n^{\alpha }\big )
\le \exp \big \{ -n a_n \big(n^{\alpha }(\alpha \log n-1)+1\big) \big \}
+ \mathbf{P}\big( U_{1:n}\le a_n\big ) .
\end{equation}
Setting $a_n=1/n^{1+\alpha }$ makes the right-hand side of bound (\ref{p3_13m}) converge to $0$ when $n\to \infty $, because $\alpha>0$. This proves statement (\ref{p3_13c}) and, in turn, statement (\ref{p3_13}), thus completing the proof that $r_{2,n}^{(2)}$ converges to $0$ in probability.

The treatment of the remainder term $r_{2,n}^{(3)}$ is analogous to that of $r_{2,n}^{(2)}$, and so we omit it. This concludes the proof of Theorem~\ref{Thm-3}.  $\hfill \square $

\subsection{Proof of Theorem~\ref{Thm-1}}
\label{proof-1}

With $A_n:=n^{-1/2} \sum_{i=2}^n (Y_{i,n}-Y_{i-1,n})_{+}$ and $ B_n:=n^{-1/2} \sum_{i=2}^n |Y_{i,n}-Y_{i-1,n}|$, we have the equation $\mathrm{I}_n=A_n/B_n$.  Denote $\mathbf{X}=(X_1,\dots,X_n)$ and write
\begin{equation}\label{p1}
\mathrm{I}_n=
{\big( A_n - \mathbf{E}(A_n\mid \mathbf{X})\big)+ \mathbf{E}(A_n\mid \mathbf{X})
\over
\big( B_n - \mathbf{E}(B_n\mid \mathbf{X})\big)+ \mathbf{E}(B_n\mid \mathbf{X})
}.
\end{equation}
We shall next prove the statement
\begin{equation}\label{p3}
A_n - \mathbf{E}(A_n\mid \mathbf{X})=O_{\mathbf{P}}(1).
\end{equation}
The proof of an analogous statement for $B_n$ is almost identical, and we shall therefore only give a few cursory remarks related to it.

We start the proof of statement \eqref{p3} by splitting $A_n$ into the sum $A_n^{(1)} +A_n^{(2)}$, where
\[
A_n^{(1)} ={1\over \sqrt{n}} \sum_{k=1}^{[n/2]} (Y_{2k,n}-Y_{2k-1,n})_{+}
\]
and
\[
A_n^{(2)} ={1\over \sqrt{n}} \sum_{k=1}^{[n/2]-\mathbf{1}_{ev}(n)} (Y_{2k+1,n}-Y_{2k,n})_{+}
\]
with $\mathbf{1}_{ev}(n)$ equal to $1$ if $n$ is even and $0$ otherwise. The idea of  splitting is to make the index-pairs of the summands inside each of the two sums disjoint, which will enable us to evoke independence arguments. Obviously, statement \eqref{p3} follows if
\begin{equation}\label{p4}
A_n^{(j)} - \mathbf{E}(A_n^{(j)}\mid \mathbf{X})=O_{\mathbf{P}}(1)
\end{equation}
for $j=1$ and $2$, which we prove only when $j=1$ because the case $j=2$ is virtually identical.

For every fixed $\varepsilon>0$, we have
\begin{equation}\label{p5}
\mathbf{P}\big( |A_n^{(1)} - \mathbf{E}(A_n^{(1)}\mid \mathbf{X})|>\varepsilon\big)
=\mathbf{E}\Big[
\mathbf{P}\big(|A_n^{(1)} - \mathbf{E}(A_n^{(1)}\mid \mathbf{X})|>\varepsilon \mid \mathbf{X}\big)\Big].
\end{equation}
To estimate the conditional probability inside the expectation, we use Chebyshev's inequality and the fact (Bhattacharya, 1974; Lemma 1) that, conditionally on $X_{1:n}, \dots, X_{n:n}$, the concomitants $Y_{1,n}, \dots , Y_{n,n}$ are independent and follow the cdf's $G(y\mid X_{1:n}), \dots , G(y\mid X_{n:n})$, respectively, where
$G(y\mid X_{i:n} )=\textbf{P}(Y\leq y\mid X_{i:n} )$. Hence, conditionally on $\mathbf{X}$, the summands $(Y_{2k,n}-Y_{2k-1,n})_{+}$, $k=1,\dots,[n/2]$, are independent.  Chebyshev's inequality implies
\begin{align}\label{p6}
\mathbf{P}\big ( |A_n^{(1)} - \mathbf{E}(A_n^{(1)}\mid \mathbf{X})|>\varepsilon\,\mid  \mathbf{X}\big )
\leq & \frac 1{n\varepsilon^2}\sum_{k=1}^{[n/2]} \mathbf{Var}\big((Y_{2k,n}-Y_{2k-1,n})_{+}\mid  \mathbf{X}\big)
\notag
\\
\leq & \frac c{n\varepsilon^2}\sum_{k=1}^{[n/2]} \mathbf{E}\big(Y_{[2k:n]}^2+Y_{[2k-1:n]}^2\mid  \mathbf{X}\big)
\notag
\\
\le &\frac c{n\varepsilon^2}\sum_{i=1}^{n} \mathbf{E}\big(Y_{i,n}^2\mid  \mathbf{X}\big)
\notag
\\
= &\frac c{n\varepsilon^2}\sum_{i=1}^{n} \mathbf{E}\big(Y_i^2\mid   X_i\big).
\end{align}
Plugging in this estimate on the right-hand side of equation \eqref{p5}, we have \begin{equation}\label{p7}
\mathbf{P}\big(|A_n^{(1)} - \mathbf{E}(A_n^{(1)}\mid \mathbf{X})|>\varepsilon\big)
\le \frac c{\varepsilon^2}\mathbf{E}[Y^2],
\end{equation}
which implies statement \eqref{p4} when $j=1$. With minor modifications, the same proof leads to statement \eqref{p4} when $j=2$, thus concluding the proof of statement \eqref{p3}. With some minor modifications, the proof of statement \eqref{p3} also leads to the statement
\begin{equation}\label{p3-b}
B_n - \mathbf{E}(B_n\mid \mathbf{X})=O_{\mathbf{P}}(1) .
\end{equation}
Hence, we conclude that, due to statements \eqref{p3} and \eqref{p3-b}, equation (\ref{p1}) implies
\begin{equation}\label{p1-p1}
\mathrm{I}_n=
{\mathbf{E}(A_n\mid \mathbf{X})+O_{\mathbf{P}}(1)
\over
\mathbf{E}(B_n\mid \mathbf{X})+O_{\mathbf{P}}(1)
}.
\end{equation}

We continue the proof of Theorem~\ref{Thm-1} with the equations
\begin{align}
\mathbf{E}(A_n\mid \mathbf{X})
&= {1\over \sqrt{n}} \sum_{i=2}^n \mathbf{E}\big((Y_{i,n}-Y_{i-1,n})_{+}\mid  \mathbf{X}\big)
\notag
\\
&= {1\over \sqrt{n}} \sum_{i=2}^n\mathbf{E}\big((Y_{i,n}-Y_{i-1,n})_{+}\mid  X_{i:n},X_{i-1:n}\big),
\label{p8}
\end{align}
where the last equation holds because only the pair $(X_{i:n},X_{i-1:n})$ among the coordinates of $\mathbf{X}$ is relevant for the distribution of the pair $(Y_{i,n},Y_{i-1,n})$. Furthermore, according to Bhattacharya (1974), the concomitants $Y_{i,n}$ and $Y_{i-1,n}$ are conditionally independent and follow the cdf's $G_i(y):=G(y\mid X_{i:n})$ and $G_{i-1}(y):=G(y\mid X_{i-1:n})$, respectively. Hence,
\begin{align*}
\mathbf{E}\big((Y_{i,n}-Y_{i-1,n})_{+}\mid X_{i:n},X_{i-1:n}\big)
=& \int_{-\infty}^{\infty}\int_{-\infty}^{x} (x-y)\,\mathrm{d}G_i(x)\,\mathrm{d}G_{i-1}(y)
\\
=& \int_{-\infty}^{\infty} xG_{i-1}(x) \, \mathrm{d}G_i(x) - \int_{-\infty}^{\infty}x (1-G_i(x))\,\mathrm{d}G_{i-1}(x).
\end{align*}
Plugging in this formula into the right-hand side of equation (\ref{p8}), we have
\begin{equation}\label{sum-1}
\mathbf{E}(A_n\mid \mathbf{X})=S_n^{(1)} +S_n^{(2)}+S_n^{(3)},
\end{equation}
where
\begin{align*}
S_n^{(1)}=&{1\over \sqrt{n}} \sum_{i=2}^n \int_{-\infty}^{\infty} xG_{i-1}(x) \, \mathrm{d}G_i(x) ,
\\
S_n^{(2)}=&-{1\over \sqrt{n}} \sum_{i=2}^n \int_{-\infty}^{\infty}x \,\mathrm{d}G_{i-1}(x),
\\
S_n^{(3)}=&{1\over \sqrt{n}} \sum_{i=2}^n \int_{-\infty}^{\infty}xG_i(x) \,\mathrm{d}G_{i-1}(x).
\end{align*}

Following the same arguments that had led us to equation \eqref{p8}, we now arrive at
\begin{align}
\mathbf{E}(B_n\mid \mathbf{X})
& = {1\over \sqrt{n}} \sum_{i=2}^n\mathbf{E}\big(|Y_{i,n}-Y_{i-1,n}|\mid X_{i:n},X_{i-1:n}\big)
\notag
\\
& = {1\over \sqrt{n}} \sum_{i=2}^n \int_{-\infty}^{\infty}\int_{-\infty}^{\infty} |x-y\mid \mathrm{d}G_i(x)\,\mathrm{d}G_{i-1}(y)
\notag
\\
& = S_n^{(1)} +S_n^{(2)}+S_n^{(3)}
+ {1\over \sqrt{n}} \sum_{i=2}^n \int_{-\infty}^{\infty}\int_{x}^{\infty} (y-x)\,\mathrm{d}G_i(x)\,\mathrm{d}G_{i-1}(y).
\label{p13}
\end{align}
Note that by the Fubini theorem, the right-most double integral is equal to the penultimate one if we interchange $G_i$ and $G_{i-1}$, and since the penultimate one is equal to $\mathbf{E}(A_n\mid \mathbf{X})$, from equation (\ref{sum-1}) we conclude that the double integral on the right-hand side of equation (\ref{p13}) is equal to $S_n^{(1)} +S_n^{(2)}+S_n^{(3)}  +O_{\mathbf{P}}(n^{-1/2}) $, with a remainder term $O_{\mathbf{P}}(n^{-1/2})$ added because $(-1)n^{-1/2} \sum_{i=2}^n \int_{-\infty}^{\infty}x \,\mathrm{d}G_{i}(x)$ is equal to $S_n^{(2)}+O_{\mathbf{P}}(n^{-1/2}) $. Hence, in summary, we have
\begin{align*}
\mathbf{E}(B_n\mid \mathbf{X})
&= 2\big (S_n^{(1)} +S_n^{(2)}+S_n^{(3)}\big ) +O_{\mathbf{P}}(n^{-1/2})
\\
&= 2\mathbf{E}(A_n\mid \mathbf{X})  +O_{\mathbf{P}}(n^{-1/2}) .
\end{align*}
Using this equation on the right-hand side of equation (\ref{p1-p1}), we obtain
\begin{equation}\label{p1-5a}
\mathrm{I}_n
 = {\mathbf{E}(B_n\mid \mathbf{X})/2+O_{\mathbf{P}}(1)
\over
\mathbf{E}(B_n\mid \mathbf{X})+O_{\mathbf{P}}(1)}
= {B_n/2+O_{\mathbf{P}}(1) \over B_n+O_{\mathbf{P}}(1)},
\end{equation}
where the right-most equation is due to statement \eqref{p3-b}. Since $B_n\stackrel{\mathbf{P}}{\to} \infty $ when $n\to \infty $, the right-hand side of equation (\ref{p1-5a}) converges to $1/2$ in probability. This concludes the proof of Theorem~\ref{Thm-1}. $\hfill \square $

\section{A summary and concluding notes}
\label{conclude}

We have shown how to assess monotonicity, or lack of it, of transfer functions in any window of interest based on the knowledge of data consisting of random inputs and their corresponding outputs. This enables researchers and decision makers to compare the current status of filters with their original status, and in this way helps to identify and assess potential structural changes and other abberations. The results also enable researchers to compare several filters functioning at the same time. We have also discussed potential difficulties arising in situations when outputs are contaminated by measurement errors. To aid well-informed uses of our results and to facilitate their extensions, if necessary, we have carefully specified assumptions and provided detailed proofs of the main results; they also delineate the applicability boundaries of the herein proposed methodology. Finally, we have provided numerical and graphical illustrations of the methodology, and in this way demonstrated the feasibility of its practical implementation.

\section*{Acknowledgement}

Research of the second author has been supported by the Natural Sciences and Engineering Research Council of Canada.

\end{document}